\documentstyle[12pt]{article}
\newtheorem{LE}{Lemma}
\newtheorem{TH}{Theorem}
\newtheorem{PR}{Proposition}
\newtheorem{CO}{Corollary}
\newtheorem{EX}{Example}

\newcommand{\pf}{\medskip\noindent{\sc Proof: }}
\newcommand{\la}{\longrightarrow}
\newcommand{\rTo}{\,\longrightarrow\,}
\newcommand{\qed}{$\Box$}
\newcommand{\BF}[1]{\mbox{\boldmath $#1$}}

\title{Oriented Quantum Algebras and Coalgebras, Invariants of
Oriented  $1$--$1$ Tangles, Knots and Links}
\author{Louis Kauffman\thanks{Research supported in part by NSF Grant
DMS 920-5227} $\;\;$ and $\;\;$ David E. Radford\thanks{Research
supported in part by NSF Grant
DMS 980 2178}\\  \\
Department of Mathematics, Statistics \\
and Computer Science (m/c 249)    \\
851 South Morgan Street   \\
University of Illinois at Chicago\\
Chicago, Illinois 60607-7045}
\begin{document}
\maketitle
This paper is the third in a series on oriented quantum algebras,
structures related to them, and regular isotopy invariants
associated with them. There is always a regular isotopy invariant
of oriented $1$--$1$ tangles associated to an oriented quantum
algebra. Regular isotopy invariants of oriented knots and links
can be constructed from oriented quantum algebras with a bit more
structure. These are the twist oriented quantum algebras and they
account for a very large number of the known regular isotopy
invariants of oriented knots and links.

In this paper we study oriented quantum coalgebras which are
structures closely related to oriented quantum algebras. We study
the relationship between oriented quantum coalgebras and oriented
quantum algebras and the relationship between oriented quantum
coalgebras and quantum coalgebras. We show that there are regular
isotopy invariants of oriented $1$--$1$ tangles and of oriented
knots and links associated to oriented and twist oriented quantum
coalgebras respectively. There are many parallels between the
theory of oriented quantum coalgebras and the theory of quantum
coalgebras; the latter are introduced and studied in \cite{QC}.

In the first paper \cite{RKO} of this series the notion of
oriented quantum algebra is introduced in the context of a very
natural diagrammatic formalism. In the second \cite{KRoalg2} basic
properties of oriented quantum algebras are described. Several
examples of oriented quantum algebras are given, one of which is a
paramaterized family which accounts for the Jones and HOMFLY
polynomials.


This paper is organized as follows. In Section \ref{SecPrelim} we
review most of the coalgebra prerequisites for this paper. Not
much is required. The theory of coalgebras needed for this paper
is more than adequately covered in any of \cite{LRBOOK,Mont,SBK}.
In Section \ref{SecDefExa} we review the notions of quantum
algebra, quantum coalgebra and their oriented counterparts. We
also recall examples of oriented quantum algebras described in
\cite{RKO,RK2DIM}. We explore duality relationships between
oriented algebra and coalgebra structures.

Section \ref{SecUOAandOA} is devoted to the relationship between
oriented quantum algebras and quantum algebras. We have shown that
a quantum algebra has an oriented quantum algebra structure. Here
we show how to associate a quantum algebra to an oriented quantum
algebra in a very natural way. In Section \ref{SecGenThry} we
prove some general results on oriented quantum coalgebras. Section
\ref{SecUOcandOC} is the coalgebra version of Section
\ref{SecUOAandOA}.

In Section \ref{SecInvar} we define a function from the set of
oriented $1$--$1$ tangle diagrams with respect to a vertical to
the dual algebra of an oriented quantum coalgebra and prove that
this function determines a regular isotopy invariant of oriented
$1$--$1$ tangles. In Section \ref{Sec11TINV} we show that the
invariant of Section \ref{SecInvar} is no better than the writhe
when the oriented quantum coalgebra is cocommutative. One would
expect this to be the case since the invariant of oriented
$1$--$1$ diagrams constructed from a commutative oriented quantum
algebra has the same property.

In Section \ref{SecKnotINV} the construction of the invariant of
oriented $1$--$1$ tangles described in Section \ref{SecInvar} is
used to give an invariant of oriented knots and links when the
oriented quantum coalgebra is a twist oriented quantum coalgebra.
The invariant for knots and links is a scalar.

This paper and some of its results were described in the survey
paper \cite{MSRI}. Throughout $k$ is a field and $k^\star$ will
denote the set of non-zero elements of $k$.

\section{Preliminaries}\label{SecPrelim}
For vector spaces $U$ and $V$ over $k$ we will denote the tensor
product $U{\otimes}_kV$ by $U {\otimes} V$, the identity map of
$V$ by $1_V$ and the linear dual ${\rm Hom}_k(V, k)$ of $V$ by
$V^*$. If $T$ is a linear endomorphism of $V$ then an element $v
\in V$ is $T$-{\em invariant} if $T(v) = v$. If $A$ is an algebra
over $k$ we shall let $1_A$ also denote the unit of $k$. Then
meaning $1_V$ should always be clear from context.

We will usually denote a coalgebra $(C, \Delta, \epsilon )$ over
$k$ by $C$ and we will follow the convention of writing the
coproduct $\Delta (c)$ symbolically as $\Delta (c) = c_{(1)}
{\otimes} c_{(2)}$ for all $c \in C$. This way of writing $\Delta
(c)$ is a variation of the Heyneman--Sweedler notation. The {\em
opposite coalgebra}, which we denote by $C^{cop}$, is $(C,
\Delta^{cop}, \epsilon )$, where $\Delta^{cop}(c) =
c_{(2)}{\otimes} c_{(1)}$ for all $c \in C$. An element $c \in C$
is said to be {\em cocommutative} if $\Delta (c) = c_{(2)}
{\otimes} c_{(1)} = \Delta^{cop}(c)$. The coalgebra $C$ is said to
be a {\em cocommutative coalgebra} if all of its elements are
cocommutative, or equivalently if $C = C^{cop}$.

Set $\Delta^{(1)} = \Delta$ and define $\Delta^{(n)} : C \rTo C
\otimes \cdots \otimes C$ ($n+1$ summands) for $n > 1$ inductively
by $ \Delta^{(n)} = (\Delta \otimes 1_C \otimes \cdots \otimes
1_C){\circ}\Delta^{(n-1)}. $ We generalize our notation for the
coproduct and write $\Delta^{(n-1)}(c) = c_{(1)} \otimes c_{(2)}
{\otimes} \cdots {\otimes} c_{(n)}$ for all $c \in C$.

A coalgebra which the reader will encounter several times in this
paper is the {\em comatrix coalgebra} ${\rm C}_n(k)$ which is
defined for all $n \geq 1$. As a $k$-vector space ${\rm C}_n(k)$
has basis $\{ e^\imath_\jmath \}_{1 \leq \imath, \jmath \leq n}$.
The coproduct and the counit for ${\rm C}_n(k)$ are determined by
$$
\Delta (e^\imath_\jmath) = \sum_{\ell = 1}^n e^\imath_\ell
{\otimes} e^\ell_\jmath \quad \mbox{and} \quad \epsilon
(e^\imath_\jmath ) = \delta^\imath_\jmath
$$
respectively for all $1 \leq \imath, \jmath \leq n$.

We usually denote an algebra $(A, m, \eta )$ over $k$ by $A$. The
{\em opposite algebra} is the $k$-algebra $(A, m^{op}, \eta )$
whose product is defined by $m^{op}(a {\otimes} b) = m (b
{\otimes} a) = ba$ for all $a, b \in A$. We denote the opposite
algebra by $A^{op}$. For $n \geq 1$ let ${\rm M}_n(k)$ the algebra
of all $n {\times} n$ matrices with entries in $k$ and let $\{
E^\imath_\jmath \}_{1 \leq \imath, \jmath \leq n}$ be the standard
basis for ${\rm M}_n(k)$. In our notation $E^\imath_\jmath
E^\ell_m = \delta^\ell_\jmath E^\imath_m$ for all $1 \leq \imath,
\jmath, \ell, m \leq n$.

Let $C$ be a coalgebra over $k$. Then $C^*$ is an algebra over
$k$, called the {\em dual algebra}, or {\em algebra dual to} $C$,
whose product  is determined by
$$
c^*d^*(c) = c^*(c_{(1)})d^*(c_{(2)})
$$
for all $c^*, d^* \in C^*$ and $c \in C$ and whose unit is
$\epsilon$. Note that $C$ is a $C^*$-bimodule under the left and
right actions
$$
c^*{\rightharpoonup}c = c_{(1)}c^*(c_{(2)}) \quad \mbox{and} \quad
c{\leftharpoonup} c^* = c^*(c_{(1)})c_{(2)}
$$
for all $c^* \in C^*$ and $c \in C$.

Now suppose that $A$ is an algebra over $k$. Then the subspace
$A^o$ of $A^*$ consisting of all functionals which vanish on a
cofinite ideal of $A$ is a coalgebra over $k$. For $a^o \in A^o$
the coproduct $\Delta (a^o) = \sum_{\imath = 1}^r a^o_\imath
{\otimes} b^o_\imath$ is determined by
$$
a^o (ab) = \sum_{\imath = 1}^r a^o_\imath (a)b^o_\imath (b)
$$
for all $a, b \in A$ and the counit is given by $\epsilon (a^o) =
a^o(1)$. If $f : A \la B$ is an algebra map then the restriction
$f^o$ of the transpose map $f^* : B^* \la A^*$ determines a
coalgebra map $f^o : B^o \la A^o$. Note that $A^o = A^*$ when $A$
is finite-dimensional. Observe that ${\rm C}_n(k) \simeq {\rm
M}_n(k)^*$ as coalgebras and that $\{ e^\imath_\jmath \}_{1 \leq
\imath, \jmath \leq n}$ can be identified with the basis for ${\rm
C}_n(k)$ dual to the standard basis $\{ E^\imath_\jmath \}_{1 \leq
\imath, \jmath \leq n}$ for ${\rm M}_n(k)$.

Let $V$ be a vector space over $k$ and suppose that $b : V
{\times} V \la k$ is a bilinear form. We define $b_{(\ell )},
b_{(r)} : V \la V^*$ by $b_{(\ell )}(u)(v) = b(u, v) =
b_{(r)}(v)(u)$ for all $u, v \in V$. For $\rho \in V {\otimes} V$
we define a bilinear form $b_\rho : V^* {\otimes} V^* \la k$ by
$b_\rho (u^*, v^*) = (u^* {\otimes} v^*)(\rho)$ for all $u^*, v^*
\in V^*$.

Now suppose that $C, D$ are coalgebras over $k$ and let $b, b' : C
{\times} D \la k$ be bilinear forms. Then $b'$ is an {\em inverse
for} $b$ if
$$
b'(c_{(1)}, d_{(1)})b(c_{(2)}, d_{(2)}) = \epsilon (c)\epsilon (d)
= b(c_{(1)}, d_{(1)})b'(c_{(2)}, d_{(2)})
$$
for all $c \in C$ and $d \in D$. The bilinear form $b$ has at most
one inverse which we denote $b^{-1}$ when it exists.
\section{Oriented Quantum Algebras and Coalgebras, Definitions and Examples}\label{SecDefExa}
In this section we recall the definition of quantum algebra,
oriented quantum algebra, quantum coalgebra and oriented quantum
coalgebra and list some examples of these structures which are
found in \cite{RKO} and \cite{RK2DIM}. We consider duality
relations between these algebra and coalgebra structures.

An important component of the definition of quantum algebra or
oriented quantum algebra is a solution to a the quantum
Yang--Baxter equation. Let $A$ be an algebra over the field $k$,
$\rho \in A {\otimes} A$ and write $\rho = \sum_{\imath = 1}^r
a_\imath {\otimes} b_\imath$. For $1 \leq \imath < \jmath \leq 3$
let $\rho_{\imath \, \jmath} \in A {\otimes} A {\otimes} A$ be
defined by
$$
\rho_{1\, 2} = \sum_{\imath = 1}^r a_\imath {\otimes} b_\imath
{\otimes} 1, \;\;\; \rho_{1\, 3} = \sum_{\imath = 1}^r a_\imath
{\otimes} 1 {\otimes} b_\imath \;\;\; \mbox{and} \;\;\; \rho_{2\,
3} = \sum_{\imath = 1}^r 1 {\otimes} a_\imath {\otimes} b_\imath.
$$
The quantum Yang--Baxter equation for $\rho$ is $\rho_{1\, 2}
\rho_{1\, 3} \rho_{2\, 3}  = \rho_{2\, 3} \rho_{1\, 3}\rho_{1\,
2}$.

The notion of quantum algebra is defined in \cite{KNOTS}. A {\em
quantum algebra} over the field $k$ is a triple $(A, \rho, s)$,
where $A$ is an algebra over $k$, $\rho \in A \otimes A$ is
invertible and $s : A \rTo A^{op}$ is an algebra isomorphism,
such that
\begin{flushleft}
{\rm (QA.1)} $\rho^{-1} = (s \otimes 1_A)(\rho)$, \vskip1\jot {\rm
(QA.2)} $\rho = (s \otimes s)(\rho)$ and \vskip1\jot {\rm (QA.3)}
$\rho_{1\, 2} \rho_{1\, 3} \rho_{2\, 3}  = \rho_{2\, 3} \rho_{1\,
3}\rho_{1\, 2}$.
\end{flushleft}
Suppose that $(A', \rho', s')$ is a quantum algebra over $k$ also.
A {\em morphism of quantum algebras} $f : (A, \rho, s) \la (A',
\rho', s')$ is an algebra map $ f : A \la A'$ which satisfies
$\rho' = (f {\otimes} f)(\rho)$ and $s' {\circ} f = f {\circ} s$.
Quantum algebras together with their morphisms under composition
form a monoidal category. The reader is referred to \cite[Section
3]{QC} at this point.

Our first example of a quantum algebra accounts for the Jones
polynomial when $k = {\rm C}\!\!\!\rule[.012in]{.01in}{.08in}\;$
is the field of complex numbers. See \cite[page 580]{KNOTS} and
also \cite[Section 3]{RK2DIM}.
\begin{EX}\label{ExJones}
Let $q \in k^\star$. Then $({\rm M}_2(k), \rho, s)$ is a quantum
algebra over the field $k$, where
$$
\rho = q^{-1}(E^1_1 {\otimes} E^1_1 + E^2_2 {\otimes} E^2_2) +
q(E^1_1 {\otimes} E^2_2 + E^2_2 {\otimes} E^1_1) + (q^{-1} - q^3)
E^1_2 {\otimes} E^2_1
$$
and $s(x) = Mx^tM^{-1}$ for all $x \in {\rm M}_2(k)$, where $M =
\left(\begin{array}{rr} 0 & q \\ -q^{-1} & 0 \end{array}\right)$.
\end{EX}

Finite-dimensional quasitriangular Hopf algebras account for a
large class of quantum algebras.
\begin{EX}\label{ExQTHA}
Let $(A, \rho)$ be a quasitriangular Hopf algebra with antipode
$s$ over the field $k$. Then $(A, \rho, s)$ is a quantum algebra
over $k$.
\end{EX}

The notion of oriented quantum algebra is introduced in
\cite[Section 1]{RKO}. An {\em oriented quantum algebra} over the
field $k$ is a quadruple $(A, \rho, t_{\sf d}, t_{\sf u})$, where
$A$ is an algebra over $k$, $\rho \in A {\otimes}A$ is invertible
and $t_{\sf d}, t_{\sf u}$ are commuting algebra automorphisms of
$A$, such that
\begin{flushleft}
{\rm (qa.1)} $(t_{\sf d}{\otimes}1_A)(\rho^{-1})$ and
$(1_A{\otimes} t_{\sf u})(\rho)$ are inverses in
$A{\otimes}A^{op}$, \vskip1\jot {\rm (qa.2)} $\rho = (t_{\sf d}
\otimes t_{\sf d})(\rho) = (t_{\sf u} \otimes t_{\sf u})(\rho)$
and \vskip1\jot {\rm (qa.3)} $\rho_{1\, 2} \rho_{1\, 3} \rho_{2\,
3}  = \rho_{2\, 3} \rho_{1\, 3}\rho_{1\, 2}.$
\end{flushleft}

An oriented quantum algebra $(A, \rho, t_{\sf d}, t_{\sf u})$ is
{\em standard} if $t_{\sf d} = 1_A$ and is {\em balanced} if
$t_{\sf d} = t_{\sf u}$. In the balanced case we write $(A, \rho,
t)$ for $(A, \rho, t, t)$.

Suppose that $(A, \rho, t_{\sf d}, t_{\sf u})$ and $(A, \rho',
t'_{\sf d}, t'_{\sf u})$ are oriented quantum algebras over $k$. A
{\em morphism of oriented quantum algebras} $f : (A, \rho, t_{\sf
d}, t_{\sf u}) \la (A', \rho', t'_{\sf d}, t'_{\sf u})$ is an
algebra map $ f : A \la A'$ which satisfies $\rho' = (f {\otimes}
f)(\rho)$, $t'_{\sf d} {\circ} f = f {\circ} t_{\sf d}$ and
$t'_{\sf u} {\circ} f = f {\circ} t_{\sf u}$. Oriented quantum
algebras together with their morphisms under composition form a
monoidal category.

Theorem 1 of \cite{RKO} accounts for an extensive family of
examples of balanced oriented quantum algebras.
\begin{EX}\label{ExALXHOM}
Let $n \geq 2$ and $x, b\!c \in k^\star$. Then $({\rm M}_n(k),
\rho, t)$ is a balanced oriented quantum algebra over $k$ where
$\rho = \sum_{\imath, \jmath, \ell, m = 1}^n \rho^{\imath \,
\ell}_{\jmath \, m} E^\imath_\jmath {\otimes} E^\ell_m$ satisfies
\begin{enumerate}
\item[{\rm a)}] $\rho^{\imath \, \ell}_{\jmath \, m} = 0$ unless
$\{ \imath, \ell \} = \{ \jmath, m\}$, \item[{\rm b)}]
$\rho^{\imath \, \jmath}_{\imath \, \jmath} \neq 0$ for all $1
\leq \imath, \jmath \leq n$, \item[{\rm c)}] $\rho^{\imath \,
\jmath}_{\jmath \, \imath}= x =  \rho^{\imath \, \imath}_{\imath
\, \imath} - b\!c/ \rho^{\imath \, \imath}_{\imath \, \imath}$ and
$\rho^{\jmath \, \imath}_{\imath \, \jmath} = 0$ for all $1 \leq
\imath < \jmath \leq n$, \item[{\rm d)}] $\rho^{\imath \,
\jmath}_{\imath \, \jmath}\rho^{\jmath \, \imath}_{\jmath \,
\imath} = b\!c$ for all $1 \leq \imath < \jmath \leq n$,
\item[{\rm e)}] for all $1 \leq \imath, \jmath \leq n$ either
$\rho^{\imath \, \imath}_{\imath \, \imath} = \rho^{\jmath \,
\jmath}_{\jmath \, \jmath}$ or $\rho^{\imath \, \imath}_{\imath \,
\imath} \rho^{\jmath \, \jmath}_{\jmath \, \jmath} = - b\!c$
\end{enumerate}
and $t(E^\imath_\jmath ) = (\omega_\imath/\omega_\jmath)
E^\imath_\jmath$ for all $1 \leq \imath, \jmath \leq n$, where
$\omega_1, \ldots, \omega_n \in k^\star$ satisfy
$$
\omega_\imath^2 = \left(\frac{\rho^{1 \, 1}_{1 \, 1}\rho^{\imath
\, \imath}_{\imath \, \imath}}{b\!c}\right)\left(\prod_{1 < \jmath
< \imath} \frac{(\rho^{\jmath \, \jmath}_{\jmath \,
\jmath})^2}{b\!c}\right)\omega_1^2
$$
for all $1 < \imath \leq n$.
\end{EX}

Let $k = {\rm C}\!\!\!\rule[.012in]{.01in}{.08in}\;$ be the field
of complex numbers and suppose that $q \in {\rm
C}\!\!\!\rule[.012in]{.01in}{.08in}\;^*$ is transcendental over
the subfield of rational numbers. When $b\!c = q^2$, $x = q^{-1} -
q^3$, $\rho^{\imath \, \imath}_{\imath \, \imath} =  q^{-1}$ for
all $1 \leq \imath \leq n$ and $\rho^{\imath \, \jmath}_{\imath \,
\jmath} =  q^2$ whenever $1 \leq \imath, \jmath \leq n$ are
distinct, then Example \ref{ExALXHOM} accounts for the HOMFLY
polynomial.

%

A quantum algebra always has an oriented quantum algebra structure
by virtue of \cite[Propositions 1 and 2]{KRoalg2}.
\begin{EX}\label{ExQAOQA}
If $(A, \rho, s)$ is a quantum algebra over the field $k$ then
$(A, \rho, s^{-2}, 1_A)$ and $(A, \rho, 1_A, s^{-2})$ are oriented
quantum algebras over $k$.
\end{EX}

A quantum algebra $(A, \rho, s)$ over $k$ may have no oriented
quantum algebra structures of the type $(A, \rho, t_{\sf d},
t_{\sf u})$ except those mentioned in the preceding example; see
Example 4 of \cite{KRoalg2}. A balanced oriented quantum algebra
$(A, \rho, t)$ over $k$ may not have a quantum algebra structure
of the type $(A, \rho, s)$; see Example 3 of \cite{KRoalg2}.

Balanced oriented quantum algebras arise in very natural ways.
\begin{EX}\label{ExQTHAs2tm2}
Let $(A, \rho)$ be a finite-dimensional quasitriangular Hopf
algebra over the field $k$ and suppose that $t$ is a Hopf algebra
automorphism of $A$ which satisfies $\rho = (t {\otimes} t)(\rho)$
and $t^2 = s^{-2}$. Then $(A, \rho, t)$ is a balanced oriented
quantum algebra.
\end{EX}
Very important examples of a finite-dimensional quasitriangular
Hopf algebras over $k$ are the quantum doubles $(D(A), {\BF
\rho})$ of finite-dimensional Hopf algebras $A$ with antipode $s$
over $k$. We write $D(A) = A^* {\otimes} A$ as a vector space.
\begin{EX}\label{ExDA}
Let $A$ be a finite-dimensional Hopf algebra over $k$ and suppose
that $t$ is a Hopf algebra automorphism of $A$ which satisfies
$t^2 = s^{-2}$. Then $(D(A), {\BF \rho}, T)$ is a balanced
oriented quantum algebra over $k$, where $T = (t^{-1})^* {\otimes}
t$.
\end{EX}
For details concerning these two examples see \cite[Corollary
2]{RKO} and the discussion preceding it.

We now turn to quantum coalgebras and oriented quantum coalgebras.
The notion of quantum coalgebra was introduced in \cite[Section
4]{QC}. Strict quantum coalgebras form an important class of
quantum coalgebras.

A {\em strict quantum coalgebra over} $k$ is a triple $(C, b, S)$,
where $C$ is a coalgebra over $k$, $b : C {\times} C \la k$ is an
invertible bilinear form and $S : C \la C^{cop}$ is a coalgebra
isomorphism, such that
\begin{flushleft}
{\rm (QC.1)} $b^{-1}(c, d) = b(S(c), d)$, \vskip1\jot {\rm (QC.2)}
$b(c, d) = b(S(c), S(d))$ and \vskip1\jot {\rm (QC.3)} $b(c_{(1)},
d_{(1)})b(c_{(2)}, e_{(1)})b(d_{(2)}, e_{(2)}) = b(c_{(2)},
d_{(2)})b(c_{(1)}, e_{(2)})b(d_{(1)}, e_{(1)})$
\end{flushleft}
for all $c, d, e \in C$. A {\em quantum coalgebra over} $k$ is a
triple $(C, b, S)$, where $C$ is a coalgebra over $k$, $b : C
{\times} C \la k$ is an invertible bilinear form, $S : C \la
C^{cop}$ is a coalgebra isomorphism of $C$ with respect to $b$,
such that (QC.1)--(QC.3) hold. That $S$ is a coalgebra isomorphism
with respect to $b$ means $S$ is a linear isomorphism which
satisfies $\epsilon {\circ} S = \epsilon$,
$$
b(S(c_{(1)}), d)b(S(c_{(2)}), e) = b(S(c)_{(2)}, d)b(S(c)_{(1)},
e)  \;\; \mbox{and}
$$
$$
b(d, S(c_{(1)}))b(e, S(c_{(2)})) = b(d, S(c)_{(2)})b(e,
S(c)_{(1)})
$$
for all $c, d, e \in C$.

A {\em morphism of quantum coalgebras} $f : (C, b, S) \la (C', b',
S')$ is a coalgebra map $ f : C \la C'$ which satisfies $b(c, d) =
b'(f(c), f(d))$ for all $c, d \in C$ and $S' {\circ} f = f {\circ}
S$. Quantum coalgebras over $k$ together with their morphisms
under composition form a monoidal category; the strict quantum
coalgebras over $k$ form a subcategory of this category.

The notions of quantum algebra and strict quantum coalgebra are
dual as was remarked in \cite[Section 4.1]{QC}. More formally,
\begin{PR}\label{QAdualQC}
Let $A$ be a finite-dimensional algebra over $k$, let $\rho \in A
{\otimes} A$ and suppose that $s$ is a linear automorphism of $A$.
Let $A^*$ be the dual coalgebra of $A$. Then the following are
equivalent:
\begin{enumerate}
\item[{\rm a)}] $(A, \rho, s)$ is a quantum algebra over $k$.
\item[{\rm b)}] $(A^*, b_\rho, s^*)$ is a strict quantum coalgebra
over $k$.
\end{enumerate}
\end{PR}
\qed
\medskip

A little more can be squeezed from a proof of the proposition.
\begin{CO}\label{AdualQC}
Suppose that $(A, \rho, s)$ is any quantum algebra over $k$. Then
$(A^o, b, s^o)$ is a strict quantum coalgebra over $k$, where
$b(a^o, b^o) = (a^o {\otimes} b^o)(\rho)$ for all $a^o, b^o \in
A^o$.
\end{CO}
\qed
\medskip

The strict quantum coalgebra $(A^o, b, s^o)$ of Corollary
\ref{AdualQC} is called the {\em dual quantum coalgebra of} $(A,
\rho, s)$. The dual quantum coalgebra of the quantum algebra of
Example \ref{ExJones} is a basic example of a (strict) quantum
coalgebra.
\begin{EX}\label{ExdualJones}
Let $q \in k^\star$. Then $({\rm C}_2(k), b, S)$ is a quantum
coalgebra over $k$ where
$$
b(e^1_1, e^1_1) = q^{-1} = b(e^2_2, e^2_2), \quad b(e^1_1, e^2_2)
= q = b(e^2_2, e^1_1), \quad b(e^1_2, e^2_1) = q^{-1} - q^3
$$
and $b(e^\imath_\jmath, e^\ell_m) = 0$ otherwise, and
$$
S(e^1_1) = e^2_2, \quad S(e^2_2) = e^1_1, \quad S(e^1_2) =
-q^2e^1_2 \;\; \mbox{and} \;\; S(e^2_1) = -q^{-2}e^2_1.
$$
\end{EX}
Also see \cite[Section 8]{QC}.

Just as finite-dimensional quasitriangular Hopf algebras give rise
to quantum algebras, it is easy to see, following the discussion
of \cite[Section 7.3]{LRBOOK} for example, that:

\begin{EX}
Let $(A, \beta )$ be a coquasitriangular Hopf algebra with
antipode $s$ over the field $k$. Then $(A, \beta, s)$ is a quantum
coalgebra over $k$.
\end{EX}

The notion of oriented quantum coalgebra is introduced in
\cite[Section 4]{KRoalg2}. Strict oriented quantum coalgebras form
an important class of oriented quantum coalgebras. A {\em strict
oriented quantum coalgebra over} $k$ is a quadruple $(C, b, T_{\sf
d}, T_{\sf u})$, where $C$ is a coalgebra over $k$, $b : C
{\times} C \la k$ is an invertible bilinear form and $T_{\sf d},
T_{\sf u}$ are commuting coalgebra automorphisms of $C$, such that
\begin{flushleft}
{\rm (qc.1)} $b(c_{(1)}, T_{\sf u}(d_{(2)}))b^{-1}(T_{\sf
d}(c_{(2)}), d_{(1)}) =  \epsilon (c)\epsilon (d)$ and
$\phantom{aaaaaaa} b^{-1}(T_{\sf d}(c_{(1)}), d_{(2)})b(c_{(2)},
T_{\sf u}(d_{(1)}))=  \epsilon (c)\epsilon (d)$, \vskip1\jot {\rm
(qc.2)} $b(c, d) = b(T_{\sf d}(c), T_{\sf d}(d))  = b(T_{\sf
u}(c), T_{\sf u}(d))$ and \vskip1\jot {\rm (qc.3)} $b(c_{(1)},
d_{(1)})b(c_{(2)}, e_{(1)})b(d_{(2)}, e_{(2)}) = b(c_{(2)},
d_{(2)})b(c_{(1)}, e_{(2)})b(d_{(1)}, e_{(1)})$
\end{flushleft}
for all $c, d, e \in C$. An {\em oriented quantum coalgebra over}
$k$ is a quadruple $(C, b, T_{\sf d}, T_{\sf u})$, where $C$ is a
coalgebra over $k$, $b : C {\times} C \la k$ is an invertible
bilinear form and $T$ is a coalgebra automorphism of $C$ with
respect to $\{ b, b^{-1}\}$, such that (qc.1)--(qc.3) hold.
Generally if $C, D$ are coalgebras over $k$ and ${\cal S}$ is a
set of bilinear forms $b : D {\times} D \la k$, then a linear map
(respectively isomorphism) $T : C \la D$ is a {\em coalgebra map}
(respectively {\em isomorphism}) {\em with respect to} ${\cal S}$
if
$$
b(T(c_{(1)}), d)b'(T(c_{(2)}), e) = b(T(c)_{(1)}, d)b'(T(c)_{(2)},
e)
$$
and
$$
b(d, T(c_{(1)}))b'(e, T(c_{(2)})) = b(d, T(c)_{(1)})b'(e,
T(c)_{(2)})
$$
for all $b, b' \in {\cal S}$,  $c \in C$ and $d, e \in D$. When $C
= D$ and $T$ is a coalgebra isomorphism with respect to ${\cal S}$
then $T$ is called a {\em coalgebra automorphism of} $C$ {\em with
respect to} ${\cal S}$.

An oriented quantum coalgebra $(C, b, T_{\sf d}, T_{\sf u})$ is
{\em standard} if $T_{\sf d} = 1_C$ and is {\em balanced} if
$T_{\sf d} = T_{\sf u}$. In the balanced case we write $(C, b, T)$
for $(C, b, T, T)$. A {\em morphism of oriented quantum
coalgebras} $f : (C, b, T_{\sf d}, T_{\sf u}) \la (C', b', T'_{\sf
d}, T'_{\sf u})$ is a coalgebra map $f : C \la C'$ which satisfies
$b(c, d) = b'((f(c), f(d))$ for all $c, d \in C$ and $T'_{\sf d}
{\circ} f = f {\circ} T_{\sf d}$, $T'_{\sf u} {\circ} f = f
{\circ} T_{\sf u}$. Oriented quantum coalgebras together with
their morphisms under composition form a monoidal category.

As remarked in \cite[Section 3]{RKO}, the notions of oriented
quantum algebra and strict oriented quantum coalgebra are dual. We
state here more formally:
\begin{PR}\label{OQCdualOQAf}
Suppose that $A$ is a finite-dimensional algebra over $k$, $\rho
\in A {\otimes} A$ and $t_{\sf d}, t_{\sf u}$ are commuting linear
automorphisms of $A$. Let $A^*$ be the dual coalgebra of $A$. Then
the following are equivalent:
\begin{enumerate}
\item[{\rm a)}] $(A, \rho, t_{\sf d}, t_{\sf u})$ is an oriented
quantum algebra over $k$. \item[{\rm b)}] $(A^*, b_\rho, t_{\sf
d}^*, t_{\sf u}^*)$ is a strict oriented quantum coalgebra over
$k$.
\end{enumerate}
\end{PR}
\qed
\medskip

Moreover:
\begin{CO}\label{OQCdualOQA}
Suppose that $(A, \rho, t_{\sf d}, t_{\sf u})$ is any oriented
quantum algebra over $k$. Then $(A^o, b,t_{\sf d}^o, t_{\sf u}^o)$
is a strict oriented quantum coalgebra over $k$, where $b(a^o,
b^o) = (a^o {\otimes} b^o)(\rho)$ for all $a^o, b^o \in A^o$.
\end{CO}
\qed
\medskip

The strict oriented quantum coalgebra $(A^o, b, t_{\sf d}^o,
t_{\sf d}^o)$ of Corollary \ref{OQCdualOQA} is called the {\em
dual oriented quantum coalgebra of} $(A, \rho, t_{\sf d}, t_{\sf
u})$. The duals of the strict quantum algebras of Example
\ref{ExALXHOM} form a rather extensive family of balanced strict
oriented quantum coalgebras.
\begin{EX}\label{ExdualHOMFLY}
Suppose $n \geq 2$ and $b\!c, x \in k^\star$. Suppose that $\{
\rho^{\imath \, \ell}_{\jmath \, m} \}_{1 \leq \imath, \jmath \leq
n} \subseteq k$ and $\{ \omega_\imath \}_{1 \leq \imath \leq n}
\subseteq k^\star$ satisfy conditions a)--d) and e) respectively
of Example \ref{ExALXHOM}. Then $({\rm C}_n(k), b, T)$ is a strict
balanced oriented quantum coalgebra over $k$, where
$b(e^\imath_\jmath, e^\ell_m) = \rho^{\imath \, \ell}_{\jmath \,
m}$ for all $1 \leq \imath, \jmath, \ell, m \leq n$ and
$T(e^\imath_\jmath) =
(\omega_\imath/\omega_\jmath)e^\imath_\jmath$ for all $1 \leq
\imath, \jmath \leq n$.
\end{EX}

Observe that the quantum coalgebra $({\rm C}_2(k), b, S)$ of
Example \ref{ExdualJones} has a strict oriented balanced quantum
coalgebra structure $({\rm C}_2(k), b, T)$ which is a special case
of the previous example with $\omega_1 = q^{-1}$ and $\omega_2 =
-q$. Note that $S{\circ}T = T{\circ}S$ and $T^2 = S^{-2}$.

We end this section with a result on coalgebra automorphisms with
respect to a set of bilinear forms which will be useful for the
proof of Theorem \ref{fCRegIso} of Section \ref{SecORINVThm}. Let
$C$ be a coalgebra over $k$ and suppose that ${\cal S}$ is a set
of bilinear forms $b : C {\times} C \la k$. The set of linear
automorphisms $T$ of  $C$ which satisfy (qc.2) for all $b \in
{\cal S}$ is easily seen to be a subgroup of the multiplicative
group of all linear automorphisms of $C$.
\begin{LE}\label{TmAuto}
Let $C$ be a coalgebra over the field $k$ and suppose that ${\cal
S}$ is a set of bilinear forms $b : C {\times} C \la k$.
\begin{enumerate}
\item[{\rm a)}] The set of coalgebra automorphisms $T$ of $C$ with
respect to ${\cal S}$ which satisfy (qc.2) for all $b \in {\cal
S}$ form a subgroup ${\cal G}(C, {\cal S})$ of the group of linear
automorphisms of $C$ under composition. \item[{\rm b)}] The
equations
$$
b(T^{u + \ell}(c_{(1)}), d) b'(T^{v + \ell}(c_{(2)}), e) =
b(T^u(T^\ell(c)_{(1)}), d) b'(T^v(T^\ell(c)_{(2)}), e)
$$
and
$$
b(d, T^{u + \ell}(c_{(1)})) b'(e, T^{v + \ell}(c_{(2)})) = b(d,
T^u(T^\ell(c)_{(1)})) b'(e, T^v(T^\ell(c)_{(2)}))
$$
hold for all $b, b' \in {\cal S}$, for all $T \in {\cal G}(C,
{\cal S})$, for all integers $u, v, \ell$ and $c, d, e \in C$.
\end{enumerate}
\end{LE}

\pf It is clear that the identity map of $C$ lies in ${\cal G}(C,
{\cal S})$. Suppose that $T, U \in {\cal G}(C, {\cal S})$. To
complete the proof of part a) we need only show that
$T^{-1}{\circ}U \in {\cal G}(C, {\cal S})$. For all $b, b' \in S$
and $c, d, e \in C$ observe that
\begin{eqnarray*}
\lefteqn{b(T^{-1}(U(c_{(1)})), d) b'(T^{-1}(U(c_{(2)})), e) } \\
& = & b(U(c_{(1)}), T(d)) b'(U(c_{(2)}), T(e))  \\
& = & b(U(c)_{(1)}, T(d)) b'(U(c)_{(2)}, T(e))  \\
& = & b(T(T^{-1}(U(c)))_{(1)}, T(d)) b'(T(T^{-1}(U(c)))_{(2)}, T(e))  \\
& = & b(T(T^{-1}(U(c))_{(1)}), T(d)) b'(T(T^{-1}(U(c))_{(2)}), T(e))  \\
& = & b(T^{-1}(U(c))_{(1)}, d) b'(T^{-1}(U(c)_{(2)}, e)
\end{eqnarray*}
and likewise
$$
b(d, T^{-1}(U(c_{(1)}))) b'(e, T^{-1}(U(c_{(2)}))) = b(d,
T^{-1}(U(c))_{(1)}) b'(e, T^{-1}(U(c))_{(2)}).
$$
Therefore $T^{-1}{\circ}U \in {\cal G}(C, {\cal S})$. Since
$b(T^{u + \ell}(c), d) = b(T^\ell (c), T^{-u}(d))$ and $b(d, T^{u
+ \ell}(c)) = b(T^{-u}(d), T^\ell (c))$ and for all integers $u,
\ell$ and $c, d, e \in C$, part b) follows from part a). \qed
\section{A Basic Relationship Between Oriented and Unoriented
Quantum Algebra Structures}\label{SecUOAandOA}
Suppose that $(A, \rho, s)$ is a quantum algebra over $k$. Then
$(A, \rho, 1_A, s^{-2})$ is a standard oriented quantum algebra
over $k$ by virtue of Example \ref{ExQAOQA}. We have seen that not
every oriented quantum algebra is of this form in our discussion
following Example \ref{ExQAOQA}. The main purpose of this section
is to show that a quantum algebra can be associated with an
oriented quantum algebra in a natural way.

Let $(A, \rho, 1_A, t)$ be a standard oriented quantum algebra
over $k$, let ${\cal A} = A {\oplus} A^{op}$ be the direct product
of $A$ and $A^{op}$, and let $\pi : {\cal A} \longrightarrow A$ be
the projection onto the first factor. We will construct a quantum
algebra $({\cal A}, {\BF \rho}, {\bf s})$ on ${\cal A}$ such that
$\pi : ({\cal A}, {\bf \rho}, 1_{\cal A}, {\bf s}^{-2})
\longrightarrow (A, \rho, 1_A, t)$ is a morphism of oriented
quantum algebras.

Let $\overline{(\;\;)}$ denote the linear involution of ${\cal A}$
which exchanges the direct summands of ${\cal A}$. Thus
$\overline{a {\oplus} b} = b {\oplus} a$ for all $a, b \in A$. We
regard $A$ as a subspace of ${\cal A}$ by the identification $a =
a {\oplus} 0$ for all $a \in A$. Therefore $\overline{a} = 0
{\oplus} a$ and every element of ${\cal A}$ has a unique
decomposition of the form $a + \overline{b}$ for some $a, b \in
A$. Observe that
\begin{equation}\label{EqA}
\overline{(\overline{a})} = a, \;\; \overline{ab} =
\overline{b}\overline{a} \;\; \mbox{and} \;\; a \overline{b} = 0 =
\overline{a}b
\end{equation}
for all $a, b \in A$. The main result of this section is:
\begin{TH}\label{ThmOQtoOQ}
Let $(A, \rho, t_{\sf d}, t_{\sf u})$ be an oriented quantum
algebra over the field $k$, let ${\cal A} = A {\oplus} A^{op}$ be
the direct product of $A$ and $A^{op}$ and write $\rho =
\sum_{\imath = 1}^r a_\imath {\otimes} b_\imath$, $\rho^{-1} =
\sum_{\jmath = 1}^s \alpha_\jmath {\otimes} \beta_\jmath$. Then:
\begin{enumerate}
\item[{\rm a)}] $({\cal A}, \BF{\rho}, \BF{s})$ is a quantum
algebra over $k$, where
$$
\BF{\rho} = \sum_{\imath = 1}^r (a_\imath {\otimes} b_\imath +
\overline{a_\imath} {\otimes} \overline{b_\imath}) + \sum_{\jmath
= 1}^s (\overline{\alpha_\jmath} {\otimes} \beta_\jmath +
\alpha_\jmath {\otimes} \overline{t_{\sf d}^{-1} {\circ}t_{\sf
u}^{-1} (\beta_\jmath)})
$$
and $\BF{s}(a {\oplus} b) = b {\oplus} t_{\sf d}^{-1}
{\circ}t_{\sf u}^{-1}(a)$ for all $a, b \in A$. \item[{\rm b)}]
$({\cal A}, \BF{\rho}, \BF{t_{\sf d}}, \BF{t_{\sf u}})$ is an
oriented quantum algebra over $k$, $\BF{t_{\sf d}}, \BF{t_{\sf
u}}$ commute with $\BF{s}$ and $\BF{t_{\sf d}} {\circ} \BF{t_{\sf
u}} = \BF{s}^{-2}$, where $\BF{t_{\sf d}}(a {\oplus} b) = t_{\sf
d}(a) {\oplus} t_{\sf d}(b)$ and $\BF{t_{\sf u}}(a {\oplus} b) =
t_{\sf u}(a) {\oplus} t_{\sf u}(b)$ for all $a, b \in A$.
\item[{\rm c)}] The projection $\pi : {\cal A} \la A$ onto the
first factor determines a morphism $\pi : ({\cal A}, \BF{\rho},
\BF{t_{\sf d}}, \BF{t_{\sf u}}) \la (A, \rho, t_{\sf d}, t_{\sf
u})$ of oriented quantum algebras.
\end{enumerate}
\end{TH}

\pf This result was announced as \cite[Theorem 2]{MSRI} and the
proof here was also given in that paper. We repeat the proof here
for the reader's convenience and to connect it to a proof of
Theorem \ref{ThmCdouble}.

Part b) is a straightforward calculation which is left to the
reader and part c) follows by definitions. As for part a) we may
assume that $(A, \rho, t_{\sf d}, t_{\sf u}) = (A, \rho, 1_A, t)$
is standard. In this case
$$
\BF{\rho} = \sum_{\imath = 1}^r (a_\imath {\otimes} b_\imath +
\overline{a_\imath} {\otimes} \overline{b_\imath}) + \sum_{\jmath
= 1}^s (\overline{\alpha_\jmath} {\otimes} \beta_\jmath +
\alpha_\jmath {\otimes} \overline{t^{-1}(\beta_\jmath)}) \quad
\mbox{and} \quad \BF{s}(a {\oplus} b) = b {\oplus} t^{-1}(a)
$$
for all $a, b \in A$. Since $t$ is an algebra automorphism of $A$
it follows that $t^{-1}$ is also. Thus $\BF{s} : {\cal A} \la
{\cal A}^{op}$ is an algebra isomorphism. By definition $\BF{s}(a)
= \overline{t^{-1}(a)}$ and $\BF{s}(\overline{a}) = a$ for all $a
\in A$. Since $\rho = (t {\otimes} t)(\rho)$ it follows that
$\rho^{-1} = (t^{-1} {\otimes} t^{-1})(\rho^{-1})$. At this point
it is easy to see that $\BF{\rho} = (\BF{s} {\otimes}
\BF{s})(\BF{\rho})$, or (QA.2) is satisfied for $\BF{\rho}$ and
$\BF{s}$. Using the equation $\rho^{-1} = (t^{-1} {\otimes}
t^{-1})(\rho^{-1})$ we calculate
$$
(s {\otimes}1_{\cal A})(\BF{\rho}) =  \sum_{\imath = 1}^r
\left(\overline{t^{-1}(a_\imath)} {\otimes} b_\imath + a_\imath
{\otimes} \overline{b_\imath}\right) + \sum_{\jmath = 1}^s
\left(\alpha_\jmath {\otimes} \beta_\jmath +
\overline{\alpha_\jmath} {\otimes}\overline{\beta_\jmath}\right).
$$
Using (\ref{EqA}), the equation $(t^{-1} {\otimes} 1_A)(\rho) =
(1_A {\otimes} t)(\rho)$, which follows by (qa.2), we see that
$$
\BF{\rho}((s {\otimes}1_{\cal A})(\BF{\rho})) = 1 {\otimes} 1 +
\overline{1} {\otimes} \overline{1} + \overline{1} {\otimes} 1 + 1
{\otimes} \overline{1} = 1_{\cal A} {\otimes} 1_{\cal A} = ((s
{\otimes}1_{\cal A})(\BF{\rho}))\BF{\rho}.
$$
Therefore $\BF{\rho}$ is invertible and $\BF{\rho}^{-1} = (\BF{s}
{\otimes} 1_A)(\BF{\rho})$. We have shown that (QA.1) holds for
$\BF{\rho}$ and $\BF{s}$.

That $\BF{\rho}$ satisfies (QA.3) is a rather lengthy and
interesting calculation. It is a straightforward exercise to see
that (QA.3) for $\BF{\rho}$ is equivalent to a set of eight
equations. With the notation convention $(\rho^{-1})_{\imath \,
\jmath} = \rho^{-1}_{\imath \, \jmath}$ for $1 \leq \imath <
\jmath \leq 3$, this set of eight equations can be rewritten as
set of six equations which are:
\begin{equation}\label{EqB}
\rho_{1 \, 2} \rho_{1 \, 3} \rho_{2 \, 3} = \rho_{2 \, 3} \rho_{1
\, 3} \rho_{1 \, 2},
\end{equation}
\begin{equation}\label{EqC}
\rho_{1 \, 2} \rho_{2 \, 3}^{-1} \rho_{1 \, 3}^{-1} = \rho_{1 \,
3}^{-1} \rho_{2 \, 3}^{-1} \rho_{1 \, 2},
\end{equation}
\begin{equation}\label{EqD}
\rho_{1 \, 3}^{-1} \rho_{1 \, 2}^{-1} \rho_{2 \, 3} = \rho_{2 \,
3} \rho_{1 \, 2}^{-1} \rho_{1 \, 3}^{-1},
\end{equation}
\begin{equation}\label{EqE}
\sum_{\ell = 1}^r \sum_{\jmath, m = 1}^s a_\ell \alpha_\jmath
{\otimes} \beta_\jmath \alpha_m  {\otimes}t^{-1}(\beta_m) b_\ell =
\sum_{\jmath, m = 1}^s \sum_{\ell = 1}^r \alpha_\jmath a_\ell
{\otimes} \alpha_m \beta_\jmath {\otimes}b_\ell t^{-1}(\beta_m),
\end{equation}
\begin{equation}\label{EqF}
\sum_{\jmath, m = 1}^s \sum_{\ell = 1}^r \alpha_\jmath a_\ell
{\otimes} \alpha_m t^{-1}(\beta_\jmath ) {\otimes}b_\ell \beta_m =
\sum_{\ell = 1}^r \sum_{\jmath, m = 1}^s a_\ell \alpha_\jmath
{\otimes} t^{-1}(\beta_\jmath ) \alpha_m  {\otimes}\beta_m b_\ell
\end{equation}
and
\begin{equation}\label{EqG}
\sum_{\jmath, \ell = 1}^s \sum_{\ell = 1}^r \alpha_\jmath
\alpha_\ell{\otimes} a_mt^{-1}(\beta_\jmath ) {\otimes}b_m
t^{-1}(\beta_\ell) = \sum_{\ell, \jmath = 1}^s \sum_{m = 1}^r
\alpha_\ell \alpha_\jmath {\otimes} t^{-1}(\beta_\jmath) a_m
{\otimes} t^{-1}(\beta_\ell) b_m.
\end{equation}
By assumption (\ref{EqB}) holds. Since $\rho_{\imath \, \jmath}$
is invertible and $(\rho_{\imath \, \jmath})^{-1} =
(\rho^{-1})_{\imath \, \jmath} = \rho^{-1}_{\imath \, \jmath}$,
equations (\ref{EqC})--(\ref{EqD}) hold since (\ref{EqB}) does.

We note that $t^{-1}$ is an algebra automorphism of $A$ and
$\rho^{-1} = (t^{-1} {\otimes} t^{-1})(\rho^{-1})$. Thus applying
$1_A {\otimes} t^{-1} {\otimes} 1_A$ to both sides of the equation
of (\ref{EqE}) we see that (\ref{EqE}) and (\ref{EqF}) are
equivalent; applying $t^{-1} {\otimes} 1_A {\otimes} 1$ to both
sides of (\ref{EqG}) we see that (\ref{EqG}) is equivalent to
$\rho_{2 \,3} \rho^{-1}_{1 \, 2} \rho^{-1}_{1 \, 3} = \rho^{-1}_{1
\, 3} \rho^{-1}_{1 \, 2} \rho_{2 \, 3}$, a consequence of
(\ref{EqB}). To complete the proof of part a) we need only show
that (\ref{EqE}) holds.

By assumption $(1_A {\otimes} t)(\rho)$ and $\rho^{-1}$ are
inverses in $A {\otimes} A^{op}$. Thus $\rho$ and $(1_A {\otimes}
t^{-1})(\rho)$ are inverses in $A {\otimes} A^{op}$ as $1_A
{\otimes} t^{-1}$  is an algebra endomorphism of $A {\otimes}
A^{op}$. Recall that $\rho^{-1}$ satisfies (QA.3). Therefore
\begin{eqnarray*}
\lefteqn{\sum_{\jmath, m = 1}^s \sum_{\ell = 1}^r   \alpha_\jmath a_\ell
{\otimes}\alpha_m \beta_\jmath {\otimes} b_\ell t^{-1}(\beta_m)} \\
& = &
\sum_{v, \ell = 1}^r \sum_{u, \jmath, m = 1}^s (a_v\alpha_u) \alpha_\jmath a_\ell
{\otimes} \alpha_m  \beta_\jmath {\otimes} b_\ell t^{-1}(\beta_m)(t^{-1}(\beta_u)b_v) \\
& = &
\sum_{v, \ell = 1}^r \sum_{u, \jmath, m = 1}^s a_v(\alpha_u \alpha_\jmath) a_\ell
{\otimes} \alpha_m  \beta_\jmath {\otimes} b_\ell t^{-1}(\beta_m\beta_u)b_v \\
& = &
\sum_{v, \ell = 1}^r \sum_{u, \jmath, m = 1}^s a_v(\alpha_\jmath \alpha_u) a_\ell
{\otimes} \beta_\jmath  \alpha_m {\otimes} b_\ell t^{-1}(\beta_u\beta_m)b_v \\
& = &
\sum_{v, \ell = 1}^r \sum_{u, \jmath, m = 1}^s a_v \alpha_\jmath (\alpha_u a_\ell)
{\otimes} \beta_\jmath  \alpha_m {\otimes} (b_\ell t^{-1}(\beta_u))t^{-1}(\beta_m)b_v \\
& = & \sum_{v = 1}^r \sum_{\jmath, m = 1}^s a_v \alpha_\jmath
{\otimes} \beta_\jmath  \alpha_m {\otimes} t^{-1}(\beta_m)b_v.
\end{eqnarray*}
which establishes (\ref{EqE}). \qed
\medskip

Denote by ${\cal C}_q$ the category whose objects are quintuples
$(A, \rho, s, t_{\sf d}, t_{\sf u})$, where $(A, \rho, s)$ is a
quantum algebra over $k$ and $(A, \rho, t_{\sf d}, t_{\sf u})$ is
an oriented quantum algebra over $k$ such that $t_{\sf d}, t_{\sf
u}$ commute with $s$ and $t_{\sf d} {\circ} t_{\sf u} = s^{-2}$,
and whose morphisms $f : (A, \rho, s, t_{\sf d}, t_{\sf u}) \la
(A', \rho', s', t'_{\sf d}, t'_{\sf u})$ are algebra maps $f : A
\la A'$ which determine morphisms $f : (A, \rho, s) \la (A',
\rho', s')$  and $f : (A, \rho, t_{\sf d}, t_{\sf u}) \la (A',
\rho', t'_{\sf d}, t'_{\sf u})$. The construction $({\cal A},
\BF{\rho}, \BF{s}, \BF{t_{\sf d}}, \BF{t_{\sf d}})$ of Theorem
\ref{ThmOQtoOQ} is a cofree object of ${\cal C}_q$. Let $\pi :
{\cal A} \la A$ be the projection onto the first factor.
\begin{PR}\label{UMPqaoqa}
Let $(A, \rho, t_{\sf d}, t_{\sf u})$ be an oriented quantum
algebra over the field $k$. Then the pair $(({\cal A}, \BF{\rho},
\BF{s}, \BF{t_{\sf d}}, \BF{t_{\sf u}}), \pi)$ satisfies the
following properties:
\begin{enumerate}
\item[{\rm a)}] $({\cal A}, \BF{\rho}, \BF{s}, \BF{t_{\sf d}},
\BF{t_{\sf u}})$  is an object of ${\cal C}_q$ and $\pi : ({\cal
A}, \BF{\rho}, \BF{t_{\sf d}}, \BF{t_{\sf u}}) \la (A, \rho,
t_{\sf d}, t_{\sf u})$ is a morphism of oriented quantum algebras.
\item[{\rm b)}] Suppose that  $(A', \rho', s', t'_{\sf d}, t'_{\sf
u})$ is an object of ${\cal C}_q$ and suppose that $f : (A',
\rho', t'_{\sf d}, t'_{\sf u}) \la (A, \rho, t_{\sf d}, t_{\sf
u})$ is a morphism of oriented quantum algebras. Then there is a
morphism $F : (A', \rho', s', t'_{\sf d}, t'_{\sf u}) \la ({\cal
A}, \BF{\rho}, \BF{s}, \BF{t_{\sf d}}, \BF{t_{\sf u}})$ uniquely
determined by $\pi {\circ} F = f$.
\end{enumerate}
\end{PR}

\pf We have shown part a). Let $f : (A', \rho, t'_{\sf d}, t'_{\sf
u}) \la (A, \rho, t_{\sf d}, t_{\sf u})$ be a morphism of oriented
quantum algebras. To show part b) we first suppose that $F : (A',
\rho', s', t'_{\sf d}, t'_{\sf u}) \la ({\cal A}, \BF{\rho},
\BF{s}, \BF{t_{\sf d}}, \BF{t_{\sf u}})$ is a morphism which
satisfies $\pi {\circ} F = f$. Now there are linear maps $g, h :
A' \la A$ such that $F(x) = g(x) {\oplus} h(x)$ for all $x \in
A'$. Since $\pi {\circ} F = f$ it follows that $g = f$. Since
$s{\circ}F = F{\circ} s'$ it follows that $h = g{\circ}s' =
f{\circ}s'$. Thus $F(x) = f(x) {\oplus} f(s'(x))$ for all $x \in
A'$ which establishes the uniqueness assertion of part b).

To establish the existence assertion of part b), we consider the
algebra homomorphism $F : A' \la A$ defined by $F(x) = f(x)
{\oplus} f(s'(x))$ for all $x \in A'$. Since $f : (A', \rho',
t'_{\sf d}, t'_{\sf u}) \la (A, \rho, t_{\sf d}, t_{\sf u})$ is a
morphism $f{\circ}t'_{\sf d} = t_{\sf d}{\circ}f$ and
$f{\circ}t'_{\sf u} = t_{\sf u}{\circ}f$. Therefore
\begin{eqnarray*}
{\BF s}{\circ}F(x)  & = & f(s'(x)) {\oplus} (t_{\sf d}^{-1} {\circ} t_{\sf u}^{-1} {\circ}f)(x) \\
&  =  & f(s'(x)) {\oplus} (f {\circ} t_{\sf d}'^{-1} {\circ}t_{\sf u}'^{-1}) (x) \\
&  = & f(s'(x)) {\oplus} f ((s'^2(x)) \\
&  = & F{\circ}s'(x)
\end{eqnarray*}
for all $x \in A'$ which shows that ${\BF s}{\circ} F = F {\circ}
s'$. Since $t_{\sf d}'$ and $s'$ commute  the calculation
\begin{eqnarray*}
{\BF t}_d{\circ}F(x) & = & t_{\sf d}(f(x)) {\oplus} t_{\sf d}(f(s'(x))) \\
& =  & f(t'_{\sf d}(x)) {\oplus} f(t'_{\sf d}(s'(x))) \\
&  = & f(t'_{\sf d}(x)) {\oplus} f(s'(t'_{\sf d}(x))) \\
& =  & F{\circ}t'_{\sf d}(x)
\end{eqnarray*}
for all $x \in A'$ establishes ${\BF t}_d {\circ} F = F {\circ}
t'_{\sf d}$. Likewise ${\BF t}_u {\circ} F = F {\circ} t'_{\sf
u}$.

To complete the proof that $F : (A', \rho', s', t'_{\sf d},
t'_{\sf u}) \la ({\cal A}, {\BF \rho}, {\BF s}, {\BF T}_d, {\BF
T}_d)$ is a morphism, and thus to complete the proof of the
proposition, we need only show that ${\BF \rho} =
(F{\oplus}F)({\BF \rho})$. Since $f : (A', \rho', t'_{\sf d},
t'_{\sf u}) \la (A, \rho, t_{\sf d}, t_{\sf u})$ is a morphism of
oriented quantum algebras $\rho = (f {\otimes} f)(\rho')$ and thus
$\rho^{-1} = (f {\oplus} f)(\rho'^{-1})$. Now $f{\circ}s'^2 =
t_{\sf d}^{-1}{\circ} t_{\sf u}^{-1} {\circ} f $ follows from the
hypothesis of part b). Write $\rho =  \sum_{\imath = 1}^r a_\imath
{\otimes} b_\imath$ and $\rho' =  \sum_{\imath' = 1}^{r'}
a'_{\imath'} {\otimes} b_{\imath'}$. Using the fact that $(A',
\rho', s')$ is a quantum algebra we can now calculate
\begin{eqnarray*}
(F {\otimes} F)(\rho') & = & \sum_{\imath' = 1}^{r'} F(a'_{\imath'}) {\otimes} F(b'_{\imath'}) \\
& = &  \sum_{\imath' = 1}^{r'} \left(f(a'_{\imath'}) {\oplus} f(s'(a'_{\imath'}))\right) {\otimes} \left(f(b'_{\imath'}) {\oplus} f(s'(b'_{\imath'}))\right) \\
& = &  \sum_{\imath' = 1}^{r'} \left(f(a'_{\imath'}) + \overline{ f(s'(a'_{\imath'}))}\right) {\otimes} \left(f(b'_{\imath'}) + \overline{f(s'(b'_{\imath'}))}\right) \\
& = &  \sum_{\imath' = 1}^{r'} \left(f(a'_{\imath'}) {\otimes}f(b'_{\imath'}) + \overline{ f(s'(a'_{\imath'}))} {\otimes} \overline{f(s'(b'_{\imath'}))}\right) \\
& & \quad  + \sum_{\imath' = 1}^{r'}
\left(\overline{f(s(a'_{\imath'}))} {\otimes} f(b'_{\imath'})
  + f(a'_{\imath'}) {\otimes} \overline{f(s'(b'_{\imath'}))}\right)  \\
& = &  \sum_{\imath' = 1}^{r'} \left(f(a'_{\imath'}) {\otimes}f(b'_{\imath'}) + \overline{ f(a'_{\imath'})} {\otimes} \overline{f(b'_{\imath'})}\right)  \\
& &  \quad  + \sum_{\imath' = 1}^{r'}
\left(\overline{f(s(a'_{\imath'}))} {\otimes} f(b'_{\imath'})
 + f(s'(a'_{\imath'})) {\otimes} \overline{f(s'^2(b'_{\imath'}))}\right)  \\
& = &  \sum_{\imath' = 1}^{r'} \left(f(a'_{\imath'}) {\otimes}f(b'_{\imath'}) + \overline{ f(a'_{\imath'})} {\otimes} \overline{f(b'_{\imath'})}\right)  \\
& &  \quad  + \sum_{\imath' = 1}^{r'}
\left(\overline{f(s(a'_{\imath'}))} {\otimes} f(b'_{\imath'})
 + f(s'(a'_{\imath'})) {\otimes} \overline{t_{\sf d}^{-1}{\circ} t_{\sf u}^{-1}(f(b'_{\imath'}))}\right)  \\
& = &  \sum_{\imath = 1}^r \left(a_\imath {\otimes} b_\imath + \overline{a_\imath} {\otimes} \overline{b_\imath}\right) + \sum_{\jmath = 1}^s \left(\overline{\alpha_\jmath} {\otimes} \beta_\jmath + \alpha_\jmath {\otimes}\overline{t_{\sf d}^{-1}{\circ} t_{\sf u}^{-1}(\beta_\jmath)}\right) \\
& = & {\BF \rho}.
\end{eqnarray*}
\qed
\section{General Results for Oriented Quantum Coalgebras}\label{SecGenThry}
In this section we develop some of the basic theory of oriented
quantum coalgebras. Our discussion parallels that of
\cite[Sections 4 and 5]{QC} to a good extent. Proofs of most
assertions made in this section can be obtained by modifying the
proofs of corresponding statements in \cite{QC} about quantum
coalgebras. Thus we shall tend to omit many details here.

Let $(C, b, T_{\sf d}, T_{\sf u})$ be an oriented quantum
coalgebra over the field $k$ and let $b', b'' : C {\times} C^{cop}
\la k$ be the bilinear forms defined by $b'(c, d) = b(c, T_{\sf
u}(d))$ and $b''(c, d) = b^{-1}(T_{\sf d}(c), d)$ for all $c, d
\in C$. Then the two equations of (qc.1) may be viewed as
technical formulations of the statements $b''$ is a left inverse
for $b'$ and $b''$ is a right inverse for $b'$. When $C$ is
finite-dimensional the two equations of (qc.1) are equivalent.
Thus in the finite-dimensional case, axiom (qc.1) for oriented
quantum coalgebra can be simplified.

A straightforward calculation shows that $(C^{cop}, b, T_{\sf d},
T_{\sf u})$ is an oriented quantum coalgebra over $k$. Let $T =
T_{\sf d}$ or $T = T_{\sf u}$. Since $T^{-1}$ is a coalgebra
automorphism of $C$ with respect to $\{ b, b^{-1}\}$, by part a)
of Lemma \ref{TmAuto} and the equation $b^{-1}(T(c), T(d)) =
b^{-1}(c, d)$ for all $c, d \in C$, it follows that $(C, b^{-1},
T_{\sf d}^{-1}, T_{\sf u}^{-1})$ is an oriented quantum coalgebra
over $k$. See the proof of the corresponding statement for quantum
coalgebras given in \cite[Section 4.2]{QC}. Let $b^{op} : C
{\times} C \la k$ be the bilinear form defined by $b^{op}(c, d) =
b(d, c)$ for all $c, d \in C$. Then $(C, b^{op}, T_{\sf u}, T_{\sf
d})$ is an oriented quantum coalgebra over $k$ as well.

Let $K$ be a field extension of $k$. Then $(C {\otimes} K,
b{\otimes}1_{K{\otimes} K}, T _d{\otimes} 1_K, T_{\sf u} {\otimes}
1_K)$ is an oriented quantum coalgebra over $K$, where $b
{\otimes} 1_{K {\otimes} K} (c {\otimes} \alpha, d {\otimes} \beta
) = \alpha \beta b(c, d)$ for all $c, d \in C$ and $\alpha, \beta
\in K$. Suppose that $(C', b', T'_{\sf d}, T'_{\sf u})$ is also an
oriented quantum coalgebra over $k$. Then $(C {\otimes} C', b'',
T_{\sf d}{\otimes} T'_{\sf d}, T_{\sf u}{\otimes} T'_{\sf u})$ is
an oriented quantum coalgebra over $k$, called the {\em tensor
product of oriented quantum coalgebras over} $k$, where $b'' (c
{\otimes} c', d {\otimes} d') = b(c, d)b'(c', d')$ for all $c, d
\in C$ and $c', d' \in C'$.

An {\em oriented quantum subcoalgebra of} $(C, b, T_{\sf d},
T_{\sf u})$ is an oriented quantum coalgebra $(D, b', T'_{\sf d},
T'_{\sf u})$, where $D$ is a subcoalgebra of $C$ and the inclusion
$\imath : D \la C$ determines a morphism $\imath : (D, b', T'_{\sf
d}, T'_{\sf u}) \la (C, b, T_{\sf d}, T_{\sf u})$. In this case
$b' = b|_{D {\times} D}$,  $T_{\sf d}(D) = D = T_{\sf u}(D)$ and
$T'_{\sf d} = T_{\sf d}|_D$, $T'_{\sf u} = T_{\sf u}|_D$.
Conversely, if $D$ is a subcoalgebra of $C$ and $T_{\sf d}(D) = D
= T_{\sf u}(D)$, then $(D, b_{D {\times} D}, T_{\sf d}|_D, T_{\sf
u}|_D)$ is an oriented quantum subcoalgebra of $(C, b, T_{\sf d},
T_{\sf u})$.

Let $I$ be a coideal of $C$ which satisfies $T_{\sf d}(I) = I =
T_{\sf u}(I)$ and $b(I, C) = (0) = b(C, I)$. Then the quotient
$C/I$ has a unique oriented quantum coalgebra structure $(C/I,
\overline{b}, \overline{T}_{\sf d}, \overline{T}_{\sf u})$, which
we refer to as a {\em quotient oriented quantum coalgebra
structure}, such that the projection $\pi : C \la C/I$ induces a
morphism $\pi : (C, b, T_{\sf d}, T_{\sf u}) \la (C/I,
\overline{b}, \overline{T}_{\sf d}, \overline{T}_{\sf u})$.

Let $I$ be the sum of the coideals $J$ of $C$ such that $T_{\sf
d}(J) = J = T_{\sf u}(J)$ and $b(J, C) = (0) = b(C, J)$. Then $I$
is a coideal of $C$ which satisfies $T_{\sf d}(I) = I = T_{\sf
u}(I)$ and $b(I, C) = (0) = b(C, I)$. Set $C_r = C/I$. The
quotient oriented quantum coalgebras structure  $(C_r, b_r, T_{r
\, {\sf d}}, T_{r \, {\sf u}})$ is the coalgebra counterpart of
the minimal oriented quantum subalgebra $(A_\rho, \rho, t_{\sf
d}|_{A_\rho}, t_{\sf u}|_{A_\rho})$ of an oriented quantum algebra
$(A, \rho, t_{\sf d}, t_{\sf u})$ over $k$. Observe that if $J$ is
a coideal of $C_r$ such that $T_{r \, {\sf d}}(J) = J =  T_{r \,
{\sf u}}(J)$ and $b_r(J, C_r) = (0) = b_r(C_r, J)$ then $J = (0)$.

Recall that a bilinear form $\beta : V {\times} V \la k$ define
for a vector space over $k$ is {\em left} (respectively {\em
right}) {\em non-singular} if $\beta_{(\ell )}$ (respectively
$\beta_{(r )}$) is one-one. If $b$ is either left non-singular or
right non-singular then $(C, b, T_{\sf d}, T_{\sf u})$ is strict.
\begin{LE}\label{NoNSing}
Let $(C, b, T_{\sf d}, T_{\sf u})$ be an oriented quantum
coalgebra over the field $k$ and suppose that $b$ is either left
non-singular or is right non-singular. Then $T_{\sf d}$ and
$T_{\sf u}$ are coalgebra automorphisms of $C$.
\end{LE}

\pf Consider the linear map $b_{(\ell )} {\otimes} b_{(\ell )} : C
{\otimes} C \la C^* {\otimes} C^*$ and regard $C^* {\otimes} C^*$
as a subspace of $(C {\otimes} C)^*$ by $(c^* {\otimes} d^*)(c
{\otimes} d) = c^*(c) d^*(d)$ for all $c^*, d^* \in C^*$ and $c, d
\in C$. Let $T = T_{\sf d}$ or $T = T_{\sf u}$. Then
$$
b(d, T(c_{(1)}))b(e, T(c_{(2)})) = b(d, T(c)_{(1)})b(e,
T(c)_{(2)})
$$
for all $d, c, e \in C$ which holds if and only if
$$
b_{(\ell )} {\otimes} b_{(\ell )}( T(c_{(1)}){\otimes}T(c_{(2)}))
= b_{(\ell )} {\otimes} b_{(\ell )}(
T(c)_{(1)}{\otimes}T(c)_{(2)})
$$
for all $c \in C$. Since $\epsilon {\circ} T = \epsilon$, it
follows that $T$ is a coalgebra automorphism of $C$ if $b_{(\ell
)}$ is one-one. Likewise $T$ is a coalgebra automorphism of $C$ if
$b_{(r)}$ is one-one. \qed
\medskip

Oriented quantum coalgebra structures can be pulled back just as
quantum coalgebra structures can be pulled back.
\begin{TH}\label{pullback}
Suppose that $\pi : C \la C'$ is an onto map of coalgebras over
the field $k$ and suppose that $(C', b', T_{\sf d}', T_{\sf u}')$
is an oriented quantum coalgebra structure on $C$. Then there
exists an oriented quantum coalgebra structure $(C, b, T_{\sf d},
T_{\sf u})$ on $C$ such that $\pi : (C, b, T_{\sf d}, T_{\sf u})
\la (C', b', T'_{\sf d}, T'_{\sf u})$ is a morphism.
\end{TH}

\pf The proof boils down to finding commuting linear automorphisms
$T_{\sf d}, T_{\sf u}$ of $C$ which satisfies $T'_{\sf d} {\circ}
\pi = \pi {\circ} T_{\sf d}$ and $T'_{\sf u} {\circ} \pi = \pi
{\circ} T_{\sf u}$. This is easy enough to do. The reader may want
to refer to the proof of the corresponding result for quantum
coalgebras \cite[Theorem 2]{QC}. \qed
\medskip

By \cite[Proposition 1.1.1]{LRBOOK}, for example, every
finite-dimensional coalgebra over $k$ is the quotient of  ${\rm
C}_n(k)$ for some $n \geq 1$. Thus as a corollary to Theorem
\ref{pullback}:
\begin{CO}\label{CnkC}
Every finite-dimensional oriented quantum coalgebra over the field
$k$ is the quotient of an oriented quantum coalgebra structure on
${\rm C}_n(k)$ for some $n \geq 1$.
\end{CO}
\qed
\medskip

There is an analog of \cite[Proposition 2]{KRoalg2} for quantum
coalgebras.
\begin{PR}\label{QCOQC}
If $(C, b, S)$ is a quantum (respectively strict quantum)
coalgebra over $k$ then $(C, b, 1_C, S^{-2})$ is an oriented
(respectively strict oriented) quantum coalgebra over $k$.
\end{PR}

\pf Let $(C, b, S)$ be a quantum coalgebra over $k$. To prove the
proposition we need only show that $(C, b, 1_C, S^{-2})$ is an
oriented quantum coalgebra over $k$. Since $S : C \la C^{cop}$ is
a coalgebra isomorphism with respect to $b$ and (QC.1), (QC.2)
hold for $S$ it follows that $S^{-2}$ is a coalgebra automorphism
of $C$ with respect to $\{b, b^{-1}\}$. See the proof of Lemma
\ref{TmAuto}. Since (QC.2) holds for $S$ it follows that (qc.2)
holds for $S^{-2}$. Now (QC.3) is (qc.3). Thus to complete the
proof we need only show that (qc.1) holds for $T_{\sf d} = 1_C$
and $T_{\sf d} = S^{-2}$. The calculation
\begin{eqnarray*}
b(c_{(1)}, S^{-2}(d_{(2)}))b^{-1}(c_{(2)}, d_{(1)})
& = & b(S(c_{(1)}), S^{-1}(d_{(2)}))b(S(c_{(2)}), d_{(1)}) \\
& = & b(S(c)_{(2)}, S^{-1}(d_{(2)}))b(S(c)_{(1)}, d_{(1)})  \\
& = & b^{-1}(S(c)_{(2)}, d_{(2)})b(S(c)_{(1)}, d_{(1)})  \\
& = & \epsilon(S(c)) \epsilon(d) \\
& = &  \epsilon(c) \epsilon(d)
\end{eqnarray*}
for all $c, d \in C$ shows that $b(c_{(1)},
S^{-2}(d_{(2)}))b^{-1}(c_{(2)}, d_{(1)}) =  \epsilon(c)
\epsilon(d)$ for all $c, d \in C$. Likewise $b^{-1}(c_{(1)},
d_{(2)})b(c_{(2)}, S^{-2}(d_{(1)})) =  \epsilon(c) \epsilon(d)$
for all $c, d \in C$. We have established (qc.1) for $T_{\sf d} =
1_C$ and $T_{\sf u} = S^{-2}$. \qed
\medskip

By virtue of the preceding proposition every quantum coalgebra
over $k$ has the structure of a standard oriented quantum
coalgebra. Every oriented quantum coalgebra over $k$ does also by
the analog of \cite[Proposition 1]{KRoalg2}.
\begin{PR}\label{OQCLRShift}
If $(C, b, T_{\sf d}, T_{\sf u})$ is an oriented quantum coalgebra
over $k$ then $(C, b, T_{\sf d} {\circ} T_{\sf u}, 1_C)$ and $(C,
b, 1_C, T_{\sf d} {\circ} T_{\sf u})$ are also.
\end{PR}

\pf Let $(C, b, T_{\sf d}, T_{\sf u})$ be an oriented quantum
coalgebra over $k$. By Lemma \ref{TmAuto} the composition $T_{\sf
d}{\circ} T_{\sf u}$ is a coalgebra automorphism of $C$ with
respect to $\{ b, b^{-1}\}$. The proof boils down to showing that
(qc.1) holds for $(C, b, T_{\sf d} {\circ} T_{\sf u}, 1_C)$ and
$(C, b, 1_C, T_{\sf d} {\circ} T_{\sf u})$. To show that (qc.1)
holds for $(C, b, T_{\sf d} {\circ} T_{\sf u}, 1_C)$ we note that
\begin{eqnarray*}
\lefteqn{b(c_{(1)}, d_{(2)})b^{-1}(T_{\sf d}{\circ} T_{\sf u}(c_{(2)}), d_{(1)}) } \\
& = & b(T_{\sf u}(c_{(1)}), T_{\sf u}(d_{(2)}))b^{-1}(T_{\sf u}(c_{(2)}), T_{\sf d}^{-1}(d_{(1)})) \\
& = & b(T_{\sf u}(c)_{(1)}, T_{\sf u}(d_{(2)}))b^{-1}(T_{\sf u}(c)_{(2)}, T_{\sf d}^{-1}(d_{(1)})) \\
& = & b(T_{\sf u}(c)_{(1)}, T_{\sf u}(d_{(2)}))b^{-1}(T_{\sf d}(T_{\sf u}(c)_{(2)}), d_{(1)}) \\
& = & \epsilon (T_{\sf u}(c))  \epsilon (d) \\
& = & \epsilon (c)  \epsilon (d)
\end{eqnarray*}
for all $c, d \in C$ and likewise $b^{-1}(T_{\sf d}{\circ} T_{\sf
u}(c_{(1)}), d_{(2)}) b(c_{(2)}, d_{(1)})
 = \epsilon (c)  \epsilon (d)$ for all $c, d \in C$. Similar calculations show that (qc.1) holds for $(C, b, 1_C, T_{\sf d} {\circ} T_{\sf u})$ also. The fact that $T_{\sf d}$ and $T_{\sf u}$ commute is used in the latter.
\qed
\section{A Basic Relationship Between Oriented and Unoriented Quantum Coalgebra Structures}\label{SecUOcandOC}
Let $(C, b, T_{\sf d}, T_{\sf u})$ be an oriented quantum
coalgebra over the field $k$ and let ${\cal C} = C {\oplus}
C^{cop}$ be the direct sum of the coalgebras $C$ and $C^{cop}$.
There is a quantum coalgebra structure $({\cal C}, {\BF \beta},
{\bf S})$ on ${\cal C}$ which is accounted for by \cite[Theorem
1]{RKQCS} and there is a coalgebra counterpart of Theorem
\ref{ThmOQtoOQ}.

Let $\imath : C \la {\cal C}$ be the one-one map defined by
$\imath (c) = c {\oplus} 0$ for all $c \in C$, make the
identification $c = \imath (c)$ for all $c \in C$ and define
$\overline{c{\oplus}d} = d {\oplus}c$ for all $c, d \in C$.
\begin{TH}\label{ThmCdouble}
Let $(C, b, T_{\sf d}, T_{\sf u})$ be an oriented quantum
coalgebra over the field $k$ and let ${\cal C} = C {\oplus}
C^{cop}$ be the direct sum of $C$ and $C^{cop}$. Then:
\begin{enumerate}
\item[{\rm a)}]
 $({\cal C}, {\BF \beta}, {\bf S})$ is a quantum coalgebra over $k$, where ${\BF \beta}$ is determined by
$$
{\BF \beta}(c, d) = b(c, d) =  {\BF \beta}(\overline{c},
\overline{d}), \quad {\BF \beta}(\overline{c}, d) = b^{-1}(c, d),
$$
$$
{\BF \beta}(c, \overline{d}) = b^{-1}(c, T^{-2}(d)) \;\;
\mbox{and} \;\; {\bf S}(c {\oplus} d) = T_{\sf d}^{-1} {\circ}
T_{\sf u}^{-1}(d) {\oplus}c
$$
for all $c, d \in C$. \item[{\rm b)}] $({\cal C}, {\BF \beta},
{\bf T}_{\sf d}, {\bf T}_{\sf u})$ is an oriented quantum
coalgebra over $k$ and ${\bf T}_{\sf d}, {\bf T}_{\sf u}$ commute
with ${\bf S}$, where ${\bf T}_{\sf d} (c {\oplus} d) = T_{\sf
d}(c) {\oplus} T_{\sf d}(d)$ and  ${\bf T}_{\sf u} (c {\oplus} d)
= T_{\sf u}(c) {\oplus} T_{\sf u}(d)$ for all $c, d \in C$.
\item[{\rm c)}] The inclusion $\imath : C \la {\cal C}$ induces a
morphism of oriented quantum coalgebras $\imath : (C, b, T_{\sf
d}, T_{\sf u}) \la ({\cal C}, {\BF \beta}, {\bf T}_{\sf d}, {\bf
T}_{\sf u})$.
\end{enumerate}
\end{TH}

\pf The proofs of parts b) and c) are straightforward and are left
to the reader. As for part a), we first note that $(C, b, 1_C,
T_{\sf d} {\circ} T_{\sf u})$ is an oriented quantum coalgebra
over $k$ by Proposition \ref{OQCLRShift} and that $(C, b, (T_{\sf
d} {\circ} T_{\sf u})^{-1})$  is a $(T_{\sf d} {\circ} T_{\sf
u})^{-1}$-form structure \cite[Section 3]{RKQCS}. Thus part a)
follows by \cite[Theorem 1]{RKQCS}. \qed
\medskip

The proof of \cite[Theorem 1]{RKQCS} is conceptually far more
difficult than the proof of Theorem \ref{ThmOQtoOQ}. The
formulation of \cite[Theorem 1]{RKQCS} preceded the definitions of
oriented quantum algebra and oriented quantum coalgebra. The
motivation for this theorem was to simplify calculation of
invariants of $1$--$1$ tangles which arise from certain quantum
coalgebras.

Let ${\cal C}_{cq}$ be the category whose objects are quintuples
$(C, b, S, T_{\sf d}, T_{\sf u})$, where $(C, b, S)$ is a quantum
coalgebra over $k$, $(C, b, T_{\sf d}, T_{\sf u})$ is an oriented
quantum coalgebra over $k$ and $T_{\sf d}$, $T_{\sf u}$ commute
with $S$, and whose morphisms $f : (C, b, S, T_{\sf d}, T_{\sf u})
\la (C', b', S', T'_{\sf d}, T'_{\sf u})$ are morphisms of quantum
coalgebras $f : (C, b, S) \la (C', b', S')$ and morphisms of
oriented quantum coalgebras $f : (C, b, T_{\sf d}, T_{\sf u}) \la
(C', b', T'_{\sf d}, T'_{\sf u})$. Our construction gives rise to
a free object of ${\cal C}_{cq}$. The following result, whose
proof is left to the reader, is a coalgebra counterpart of
Proposition \ref{UMPqaoqa}.
\begin{PR}
Let $(C, b, T_{\sf d}, T_{\sf u})$ be an oriented quantum
coalgebra over the field $k$. Then the pair $(\imath, ({\cal C},
\BF{\beta}, \BF{S}, \BF{T}_{\sf d}, \BF{T}_{\sf u}))$ satisfies
the following properties:
\begin{enumerate}
\item[{\rm a)}] $({\cal C}, \BF{\beta}, \BF{S}, \BF{T}_{\sf d},
\BF{T}_{\sf u})$ is an object of ${\cal C}_{cq}$ and $\imath : (C,
b, T_{\sf d}, T_{\sf u}) \la ({\cal C}, \BF{\beta}, \BF{T}_{\sf
d}, \BF{T}_{\sf u})$ is a morphism of oriented quantum coalgebras
over $k$. \item[{\rm b)}] Suppose that $(C', b', S', T'_{\sf d},
T'_{\sf u})$ is an object of ${\cal C}_{cq}$ and $f : (C, b,
T_{\sf d}, T_{\sf u}) \la (C', b', T'_{\sf d}, T'_{\sf u})$  is a
morphism of oriented quantum coalgebras over $k$. There is a
morphism $F :  ({\cal C}, \BF{\beta}, \BF{S}, \BF{T}_{\sf d},
\BF{T}_{\sf u}) \la (C', b', S', T'_{\sf d}, T'_{\sf u})$ uniquely
determined by $F {\circ} \imath = f$.\end{enumerate} \qed
\end{PR}
\section{A Regular Isotopy Invariant of Oriented $1$--$1$ Tangles
Which Arises from an Oriented Quantum Coalgebra}\label{SecInvar}
In this section we construct a regular isotopy invariant ${\bf
Inv}_C$ of oriented $1$--$1$ tangle diagrams from an oriented
quantum coalgebra $C$ over $k$ in much the same manner that we
constructed an invariant of $1$--$1$ tangle diagrams from a
quantum coalgebra over $k$ in \cite[Section 6.1]{QC}. The
invariant we construct can be considered the coalgebra version of
the invariant of oriented $1$--$1$ tangle diagrams described in
\cite[Section 1]{RKO} and \cite[Section 3]{KRoalg2} which arises
from an oriented quantum algebra. In Section \ref{SecORINV} we
describe ${\bf Inv}_C$ and in Section \ref{SecORINVThm} we prove
that ${\bf Inv}_C$ is a regular isotopy invariant of oriented
$1$--$1$ tangle diagrams (and thus determines a regular isotopy
invariant of oriented $1$--$1$ tangles).
\subsection{Invariants of Oriented $1$--$1$ Tangle Diagrams Arising
from Oriented Quantum Coalgebras}\label{SecORINV}
We represent oriented $1$--$1$ tangles as diagrams in the plane
with respect to the vertical direction. Simple examples of these
diagrams are
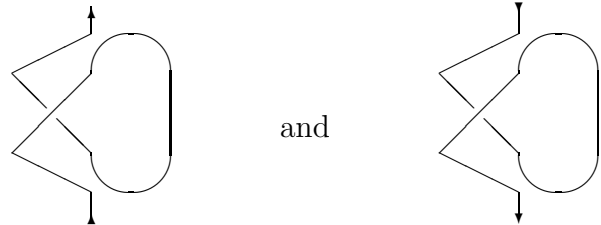
\begin{center}
\mbox{
\begin{picture}(80,90)(-10,-10)
\put(45,15){\oval(30,30)[b]} \put(45,45){\oval(30,30)[t]}
\put(60,15){\line(0,1){30}} \put(30,0){\line(0,-1){10}}
\put(30,-8){\vector(0,1){0}} \put(30,60){\line(0,1){10}}
\put(30,70){\vector(0,1){0}}
%
%
\put(0,15){\line(1,1){30}} \put(30,15){\line(-1,1){13}}
\put(0,45){\line(1,-1){13}}
%
%
\put(0,15){\line(2,-1){30}} \put(0,45){\line(2,1){30}}
\end{picture}
} \qquad \raisebox{6ex}{and} \qquad
\mbox{
\begin{picture}(80,90)(-10,-10)
\put(45,15){\oval(30,30)[b]} \put(45,45){\oval(30,30)[t]}
\put(60,15){\line(0,1){30}} \put(30,0){\line(0,-1){10}}
\put(30,-12){\vector(0,-1){0}} \put(30,60){\line(0,1){10}}
\put(30,68){\vector(0,-1){0}}
%
%
\put(0,15){\line(1,1){30}} \put(30,15){\line(-1,1){13}}
\put(0,45){\line(1,-1){13}}
%
%
\put(0,15){\line(2,-1){30}} \put(0,45){\line(2,1){30}}
\end{picture}
}
\end{center}
\noindent where the arrow heads indicate orientation. We shall
refer to the tangle diagrams above as ${\bf T}_{\rm curl}$ and
${\bf T}^{op}_{\rm curl}$ respectively. We let ${\bf Tang}$ denote
the set of all oriented $1$--$1$ tangle diagrams. If ${\bf T} \in
{\bf Tang}$ then ${\bf T}^{op}$ is the underlying diagram of ${\bf
T}$ with the opposite orientation.

The point on the tangle diagram at which a traversal of the
diagram in the direction of the orientation begins is called the
{\em base point of the diagram} and the point at which such a
traversal ends is called the {\em end point of the diagram}. We
require $1$--$1$ tangle diagrams to be completely contained in a
box except for two protruding line segments as indicated by the
two examples below.
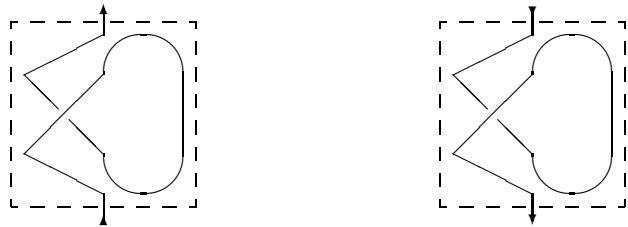
\begin{center}
\mbox{
\begin{picture}(80,90)(-10,-10)
\put(45,15){\oval(30,30)[b]} \put(45,45){\oval(30,30)[t]}
\put(60,15){\line(0,1){30}} \put(30,0){\line(0,-1){10}}
\put(30,-8){\vector(0,1){0}} \put(30,60){\line(0,1){10}}
\put(30,72){\vector(0,1){0}}
%
%
\put(0,15){\line(1,1){30}} \put(30,15){\line(-1,1){13}}
\put(0,45){\line(1,-1){13}}
%
%
\put(0,15){\line(2,-1){30}} \put(0,45){\line(2,1){30}}
%
%
\put(-5,-5){\dashbox{5}(70,70)}
\end{picture}
} \qquad \qquad \qquad \mbox{
\begin{picture}(80,90)(-10,-10)
\put(45,15){\oval(30,30)[b]} \put(45,45){\oval(30,30)[t]}
\put(60,15){\line(0,1){30}} \put(30,0){\line(0,-1){10}}
\put(30,-12){\vector(0,-1){0}} \put(30,60){\line(0,1){10}}
\put(30,67){\vector(0,-1){0}}
%
%
\put(0,15){\line(1,1){30}} \put(30,15){\line(-1,1){13}}
\put(0,45){\line(1,-1){13}}
%
%
\put(0,15){\line(2,-1){30}} \put(0,45){\line(2,1){30}}
%
%
\put(-5,-5){\dashbox{5}(70,70)}
\end{picture}
}
\end{center}

When an oriented $1$--$1$ tangle diagram ${\bf T}$ can be written
as the union of two $1$--$1$ tangle diagrams ${\bf T}_1$ and
${\bf T}_2$, where the end point of ${\bf T}_1$ is the base point
of  ${\bf T}_2$, and the horizontal line passing through this
common point otherwise separates  ${\bf T}_1$ and  ${\bf T}_2$,
then ${\bf T}$ is called the {\em product of} ${\bf T}_1$ and
${\bf T}_2$ and we write ${\bf T} = {\bf T}_1 {\star} {\bf T}_2$.
For example,
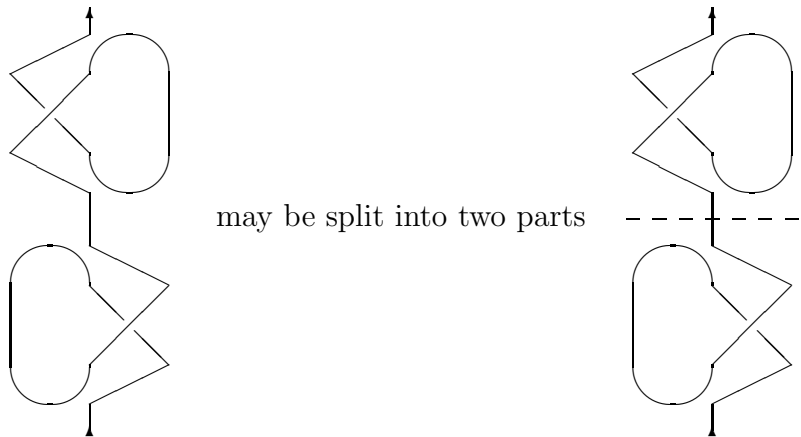
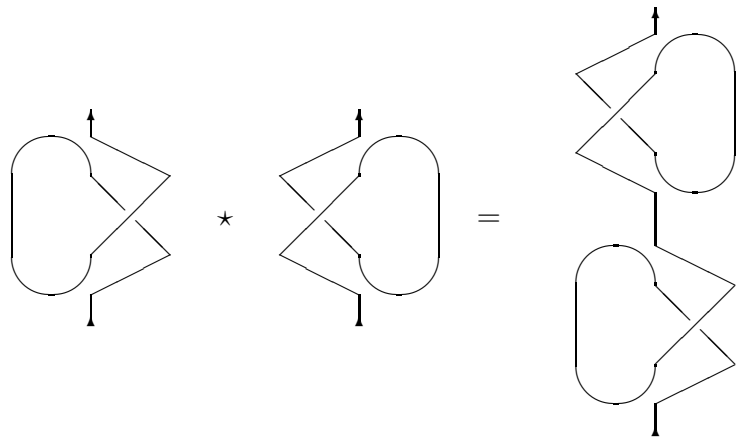
\begin{center}
\mbox{
\begin{picture}(80,160)(-10,-10)
\put(45,95){\oval(30,30)[b]} \put(45,125){\oval(30,30)[t]}
\put(60,95){\line(0,1){30}} \put(30,80){\line(0,-1){10}}
%
\put(30,140){\line(0,1){10}} \put(30,150){\vector(0,1){0}}
%
%
\put(0,95){\line(1,1){30}} \put(30,95){\line(-1,1){13}}
\put(0,125){\line(1,-1){13}}
%
%
\put(0,95){\line(2,-1){30}} \put(0,125){\line(2,1){30}}
\put(15,15){\oval(30,30)[b]} \put(15,45){\oval(30,30)[t]}
\put(0,15){\line(0,1){30}} \put(30,0){\line(0,-1){10}}
\put(30,-8){\vector(0,1){0}} \put(30,60){\line(0,1){10}}
%
%
\put(30,15){\line(1,1){30}} \put(60,15){\line(-1,1){13}}
\put(30,45){\line(1,-1){13}}
%
%
\put(30,0){\line(2,1){30}} \put(30,60){\line(2,-1){30}}
\end{picture}  }
\raisebox{15ex}{may be split into two parts} \mbox{
\begin{picture}(80,160)(-10,-10)
\put(45,95){\oval(30,30)[b]} \put(45,125){\oval(30,30)[t]}
\put(60,95){\line(0,1){30}} \put(30,80){\line(0,-1){10}}
\put(27.5,70){\line(1,0){5}} \put(37.5,70){\line(1,0){5}}
\put(47.5,70){\line(1,0){5}} \put(57.5,70){\line(1,0){5}}
\put(17.5,70){\line(1,0){5}} \put(7.5,70){\line(1,0){5}}
\put(-2.5,70){\line(1,0){5}}
\put(30,140){\line(0,1){10}} \put(30,150){\vector(0,1){0}}
%
%
\put(0,95){\line(1,1){30}} \put(30,95){\line(-1,1){13}}
\put(0,125){\line(1,-1){13}}
%
%
\put(0,95){\line(2,-1){30}} \put(0,125){\line(2,1){30}}
\put(15,15){\oval(30,30)[b]} \put(15,45){\oval(30,30)[t]}
\put(0,15){\line(0,1){30}} \put(30,0){\line(0,-1){10}}
\put(30,-8){\vector(0,1){0}} \put(30,60){\line(0,1){10}}
%
%
\put(30,15){\line(1,1){30}} \put(60,15){\line(-1,1){13}}
\put(30,45){\line(1,-1){13}}
%
%
\put(30,0){\line(2,1){30}} \put(30,60){\line(2,-1){30}}
\end{picture}  }
\end{center}
and thus
\begin{center}
\raisebox{8ex}{ \mbox{
\begin{picture}(80,100)(-10,-10)
\put(15,15){\oval(30,30)[b]} \put(15,45){\oval(30,30)[t]}
\put(0,15){\line(0,1){30}} \put(30,0){\line(0,-1){10}}
\put(30,-8){\vector(0,1){0}} \put(30,60){\line(0,1){10}}
\put(30,70){\vector(0,1){0}}
%
%
\put(30,15){\line(1,1){30}} \put(60,15){\line(-1,1){13}}
\put(30,45){\line(1,-1){13}}
%
%
\put(30,0){\line(2,1){30}} \put(30,60){\line(2,-1){30}}
%
%
\end{picture} }
\raisebox{7ex}{$\star$} \mbox{
\begin{picture}(70,100)(-10,-10)
\put(45,15){\oval(30,30)[b]} \put(45,45){\oval(30,30)[t]}
\put(60,15){\line(0,1){30}} \put(30,0){\line(0,-1){10}}
\put(30,-8){\vector(0,1){0}} \put(30,60){\line(0,1){10}}
\put(30,70){\vector(0,1){0}}
%
%
\put(0,15){\line(1,1){30}} \put(30,15){\line(-1,1){13}}
\put(0,45){\line(1,-1){13}}
%
%
\put(0,15){\line(2,-1){30}} \put(0,45){\line(2,1){30}}
%
%
\end{picture} }
\raisebox{7ex}{$\;\; = \;\;$} } \mbox{
\begin{picture}(80,200)(-10,-10)
\put(45,95){\oval(30,30)[b]} \put(45,125){\oval(30,30)[t]}
\put(60,95){\line(0,1){30}} \put(30,80){\line(0,-1){10}}
%
\put(30,140){\line(0,1){10}} \put(30,150){\vector(0,1){0}}
%
%
\put(0,95){\line(1,1){30}} \put(30,95){\line(-1,1){13}}
\put(0,125){\line(1,-1){13}}
%
%
\put(0,95){\line(2,-1){30}} \put(0,125){\line(2,1){30}}
\put(15,15){\oval(30,30)[b]} \put(15,45){\oval(30,30)[t]}
\put(0,15){\line(0,1){30}} \put(30,0){\line(0,-1){10}}
\put(30,-8){\vector(0,1){0}} \put(30,60){\line(0,1){10}}
%
%
\put(30,15){\line(1,1){30}} \put(60,15){\line(-1,1){13}}
\put(30,45){\line(1,-1){13}}
%
%
\put(30,0){\line(2,1){30}} \put(30,60){\line(2,-1){30}}
\end{picture}  }
\end{center}
Multiplication is clearly an associative operation.

Oriented $1$--$1$ tangle diagrams consist of some or all of the
following components:
\begin{enumerate}
\item[$\bullet$] oriented crossings;
\begin{enumerate}
\item[{\rm \phantom{a}}] {\em under crossings}

%
%
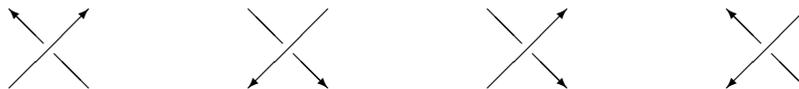
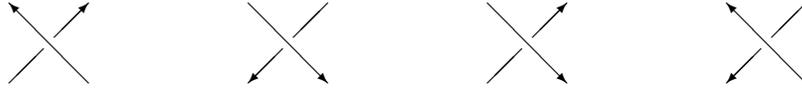
\mbox{
\begin{picture}(80,50)(-10,-10)
\put(0,0){\line(1,1){30}} \put(0,30){\line(1,-1){13}}
\put(30,0){\line(-1,1){13}} \put(0, 30){\vector(-1, 1){0}}
\put(30, 30){\vector(1, 1){0}}
\end{picture} }
%
%
\mbox{
\begin{picture}(80,50)(-10,-10)
\put(0,0){\line(1,1){30}} \put(0,30){\line(1,-1){13}}
\put(30,0){\line(-1,1){13}} \put(0, 0){\vector(-1, -1){0}}
\put(30, 0){\vector(1, -1){0}}
\end{picture} }
%
%
\mbox{
\begin{picture}(80,50)(-10,-10)
\put(0,0){\line(1,1){30}} \put(0,30){\line(1,-1){13}}
\put(30,0){\line(-1,1){13}} \put(30, 0){\vector(1, -1){0}}
\put(30, 30){\vector(1, 1){0}}
\end{picture} }
%
%
\mbox{
\begin{picture}(80,50)(-10,-10)
\put(0,0){\line(1,1){30}} \put(0,30){\line(1,-1){13}}
\put(30,0){\line(-1,1){13}} \put(0, 0){\vector(-1, -1){0}} \put(0,
30){\vector(-1, 1){0}}
\end{picture} }
\item[{\rm \phantom{a}}] {\em over crossings}

%
%
\mbox{
\begin{picture}(80,50)(-10,-10)
\put(0,0){\line(1,1){13}} \put(30,0){\line(-1,1){30}}
\put(30,30){\line(-1,-1){13}} \put(0, 30){\vector(-1, 1){0}}
\put(30, 30){\vector(1, 1){0}}
\end{picture} }
%
%
\mbox{
\begin{picture}(80,50)(-10,-10)
\put(0,0){\line(1,1){13}} \put(30,0){\line(-1,1){30}}
\put(30,30){\line(-1,-1){13}} \put(0, 0){\vector(-1, -1){0}}
\put(30, 0){\vector(1, -1){0}}
\end{picture} }
%
%
\mbox{
\begin{picture}(80,50)(-10,-10)
\put(0,0){\line(1,1){13}} \put(0,30){\line(1,-1){30}}
\put(30,30){\line(-1,-1){13}} \put(30, 0){\vector(1, -1){0}}
\put(30, 30){\vector(1, 1){0}}
\end{picture} }
%
%
\mbox{
\begin{picture}(80,50)(-10,-10)
\put(0,0){\line(1,1){13}} \put(30,30){\line(-1,-1){13}}
\put(30,0){\line(-1,1){30}} \put(0, 0){\vector(-1, -1){0}} \put(0,
30){\vector(-1, 1){0}}
\end{picture} }
\end{enumerate}
\item[$\bullet$] oriented local extrema;
\begin{enumerate}
\item[{\rm \phantom{a}}] {\em local maxima} \mbox
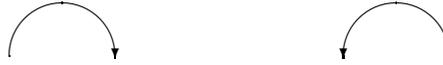
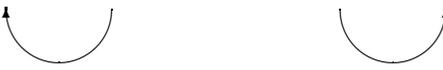
{
\begin{picture}(75,20)(-25,0)
\put(40,0){\vector(0, -1){}} \put(20,0){\oval(40, 40)[t]}
\end{picture}}
\qquad \mbox{
\begin{picture}(75,20)(-45,0)
\put(0,0){\vector(0, -1){}} \put(20,0){\oval(40, 40)[t]}
\end{picture}}

\item[{\rm \phantom{a}}] {\em local minima}
\mbox{
\begin{picture}(75,40)(-25,-20)
\put(0,0){\vector(0, 1){}} \put(20,0){\oval(40, 40)[b]}
\end{picture}}
\qquad \mbox{
\begin{picture}(75,40)(-45,-20)
\put(40,0){\vector(0, 1){}} \put(20,0){\oval(40, 40)[b]}
\end{picture}}
\end{enumerate}
\end{enumerate}
and
\begin{enumerate}
\item[$\bullet$] oriented ``vertical" lines.
\end{enumerate}
The orientations of adjoining components of the tangle diagram
must be compatible.

For an oriented quantum coalgebra $(C, b, T_{\sf d}, T_{\sf u})$
over the field $k$ the invariant we describe in this section is a
function ${\bf Inv}_C : {\bf Tang} \la C^*$ which is the function
${\bf Inv}_A$ of \cite[Section 1]{RKO} and \cite[Section
3]{KRoalg2} when $C$ is finite-dimensional, strict and $A = C^*$
is the dual quantum algebra. To motivate the definition of ${\bf
Inv}_C$ we first review how ${\bf Inv}_A$ is constructed for
oriented quantum algebras $A$ over $k$. The reader is encouraged
to refer to \cite[Section 1]{RKO} or \cite[Section 3]{KRoalg2} at
this point. Much of the discussion which follows parallels
\cite[Section 6.1]{QC}.

Let $(A, \rho, t_{\sf d}, t_{\sf u})$ be an oriented quantum
algebra defined over the field $k$ and suppose that ${\bf T} \in
{\bf Tang}$. We decorate each crossing of ${\bf T}$ according to
the scheme

%
%
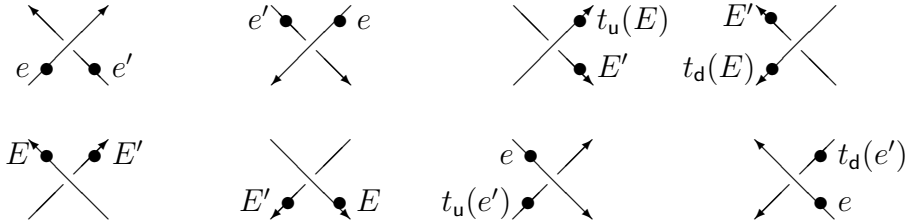
\mbox{
\begin{picture}(80,50)(-10,-10)
\put(0,0){\line(1,1){30}} \put(0,30){\line(1,-1){13}}
\put(30,0){\line(-1,1){13}} \put(-5, 3){$e \; \bullet$} \put(22,
3){$\bullet \; e'$} \put(0, 30){\vector(-1, 1){0}} \put(30,
30){\vector(1, 1){0}}
\end{picture} }
%
%
\mbox{
\begin{picture}(80,50)(-10,-10)
\put(0,0){\line(1,1){30}} \put(0,30){\line(1,-1){13}}
\put(30,0){\line(-1,1){13}} \put(23, 21){$\bullet \; e$} \put(-9,
21){$e' \; \bullet$}
\put(0, 0){\vector(-1, -1){0}} \put(30, 0){\vector(1, -1){0}}
\end{picture} }
%
%
\mbox{
\begin{picture}(80,50)(-10,-10)
\put(0,0){\line(1,1){30}} \put(0,30){\line(1,-1){13}}
\put(30,0){\line(-1,1){13}} \put(22, 21 ){$\bullet \; t_{\sf
u}(E)$} \put(22, 3){$\bullet \; E'$} \put(30, 0){\vector(1,
-1){0}} \put(30, 30){\vector(1, 1){0}}
\end{picture} }
%
%
\mbox{
\begin{picture}(80,50)(-10,-10)
\put(0,0){\line(1,1){30}} \put(0,30){\line(1,-1){13}}
\put(30,0){\line(-1,1){13}} \put(-28, 3){$ t_{\sf d}(E) \;
\bullet$} \put(-13, 22){$E' \; \bullet$} \put(0, 0){\vector(-1,
-1){0}} \put(0, 30){\vector(-1, 1){0}}
\end{picture} }

%
\mbox{
\begin{picture}(80,50)(-10,-10)
\put(0,0){\line(1,1){13}} \put(30,0){\line(-1,1){30}}
\put(30,30){\line(-1,-1){13}} \put(-9, 21){$E \; \bullet$}
\put(22, 21){$\bullet \; E'$} \put(0, 30){\vector(-1, 1){0}}
\put(30, 30){\vector(1, 1){0}}
\end{picture} }
%
%
\mbox{
\begin{picture}(80,50)(-10,-10)
\put(0,0){\line(1,1){13}} \put(30,0){\line(-1,1){30}}
\put(30,30){\line(-1,-1){13}} \put(23, 3){$\bullet \; E$}
\put(-12, 3){$E' \; \bullet$} \put(0, 0){\vector(-1, -1){0}}
\put(30, 0){\vector(1, -1){0}}
\end{picture} }
%
%
\mbox{
\begin{picture}(80,50)(-10,-10)
\put(0,0){\line(1,1){13}} \put(0,30){\line(1,-1){30}}
\put(30,30){\line(-1,-1){13}} \put(-5, 21){$e \; \bullet$}
\put(-27, 3){$t_{\sf u}(e') \; \bullet$} \put(30, 0){\vector(1,
-1){0}} \put(30, 30){\vector(1, 1){0}}
\end{picture} }
%
%
\mbox{
\begin{picture}(80,50)(-10,-10)
\put(0,0){\line(1,1){13}} \put(30,30){\line(-1,-1){13}}
\put(30,0){\line(-1,1){30}} \put(22,3){$\bullet \; e$} \put(22,
21){$\bullet \; t_{\sf d}(e')$} \put(0, 0){\vector(-1, -1){0}}
\put(0, 30){\vector(-1, 1){0}}
\end{picture}$\;\; .$ }

\noindent where we use the shorthand $\rho = e {\otimes} e'$,
$\rho^{-1} =  E{\otimes} E'$ and $(1_A {\otimes} t)(\rho) = e
{\otimes} t(e')$, $(t {\otimes} 1_A)(\rho^{-1}) = t(E) {\otimes}
E'$ for $t = t_{\sf d}, t_{\sf u}$. In practice we let $e
{\otimes} e'$, $f {\otimes} f'$, $g {\otimes} g' \ldots$ denote
copies of $\rho$ and $E {\otimes} E'$, $F {\otimes} F'$, $G
{\otimes} G' \ldots$ denote copies of $\rho^{-1}$.

Think of the oriented tangle as a rigid wire and think of the
decorations as labeled beads which slide freely around the wire.
Starting at the base point of the tangle diagram, traverse the
diagram pushing the labeled beads along the wire so that the end
result is a juxtaposition of labeled beads at the end point of the
diagram. As a labeled bead passes through a local extrema its
label is altered according to the following rules:
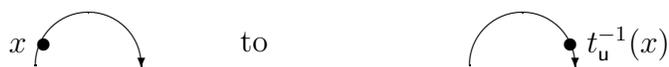
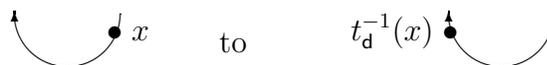
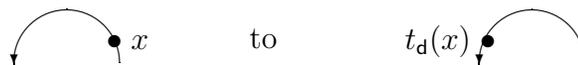
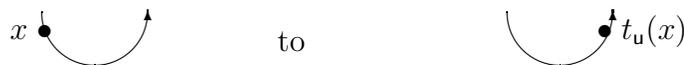
\begin{center}
\mbox{
\begin{picture}(75,40)(-25,-10)
\put(40,0){\vector(0, -1){}} \put(20,0){\oval(40, 40)[t]}
\put(-10,5){$x \; \bullet$}
\end{picture}}
\qquad \raisebox{3ex}{to} \qquad \mbox{
\begin{picture}(75,40)(-45,-10)
\put(40,0){\vector(0, -1){}} \put(20,0){\oval(40, 40)[t]}
\put(35,5){$\bullet \; t_{\sf u}^{-1}(x)$}
\end{picture}}
\end{center}

\noindent and
\begin{center}
\mbox{
\begin{picture}(75,40)(-25,-30)
\put(0,0){\vector(0, 1){}} \put(20,0){\oval(40, 40)[b]}
\put(35,-10){$\bullet \; x$}
\end{picture}}
\qquad \raisebox{3ex}{to} \qquad \mbox{
\begin{picture}(75,40)(-45,-30)
\put(0,0){\vector(0, 1){}} \put(20,0){\oval(40, 40)[b]}
\put(-37,-10){$t_{\sf d}^{-1}(x) \; \bullet$}
\end{picture}}
\end{center}
for clockwise motion;
\begin{center}
\mbox{
\begin{picture}(75,40)(-25,-10)
\put(0,0){\vector(0, -1){}} \put(20,0){\oval(40, 40)[t]}
\put(35,5){$\bullet \; x$}
\end{picture}}
\qquad  \quad \raisebox{3ex}{to} \qquad \mbox{
\begin{picture}(75,40)(-45,-10)
\put(0,0){\vector(0, -1){}} \put(20,0){\oval(40, 40)[t]}
\put(-28,5){$t_{\sf d}(x) \;\bullet$}
\end{picture}}
\end{center}

\noindent and
\begin{center}
\mbox{
\begin{picture}(75,40)(-25,-30)
\put(40,0){\vector(0, 1){}} \put(20,0){\oval(40, 40)[b]}
\put(-12,-10){$x \; \bullet$}
\end{picture}}
\qquad  \quad \raisebox{3ex}{to} \qquad \mbox{
\begin{picture}(75,40)(-45,-30)
\put(40,0){\vector(0, 1){}} \put(20,0){\oval(40, 40)[b]}
\put(34,-10){$\bullet \; t_{\sf u}(x)$}
\end{picture}}
\end{center}
for counterclockwise motion. We refer to the oriented local
extrema
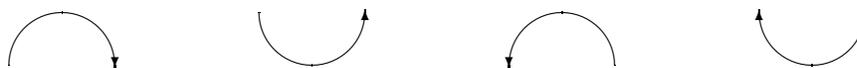
\begin{center}
\mbox{
\begin{picture}(75,40)(-25,-10)
\put(40,0){\vector(0, -1){}} \put(20,0){\oval(40, 40)[t]}
\end{picture}}
\quad \mbox{
\begin{picture}(75,40)(-25,-30)
\put(40,0){\vector(0, 1){}} \put(20,0){\oval(40, 40)[b]}
\end{picture}}
\quad \mbox{
\begin{picture}(75,40)(-25,-10)
\put(0,0){\vector(0, -1){}} \put(20,0){\oval(40, 40)[t]}
\end{picture}}
\quad \mbox{
\begin{picture}(75,40)(-25,-30)
\put(0,0){\vector(0, 1){}} \put(20,0){\oval(40, 40)[b]}
\end{picture}}
\end{center}
as having type (${\sf u}_-$), (${\sf u}_+$), (${\sf d}_+$) and
(${\sf d}_-$) respectively. Reading the juxtaposed labeled beads
in the direction of orientation results in a formal word ${\bf
W}_A({\bf T})$.

Now to define ${\bf W}_A({\bf T})$ more formally. If ${\bf T}$ has
no crossings then ${\bf W}_A({\bf T}) = 1$. Suppose that ${\bf T}$
has $n \geq 1$ crossings. Traverse ${\bf T}$ in the direction of
orientation and label the crossing lines $1, 2, \ldots, 2n$ in the
order in which they are encountered. For $1 \leq \imath \leq 2n$
let $u_{\sf d}(\imath)$ be the number of  local extrema of type
(${\sf d}_+$) minus the number of type (${\sf d}_-$) encountered
on the portion of the traversal from line $\imath$ to the end of
the traversal of ${\bf T}$. We define $u_{\sf u}(\imath)$ in the
same way where (${\sf u}_+$) and (${\sf u}_-$) replace (${\sf
d}_+$) and (${\sf d}_-$) respectively. Then
\begin{equation}\label{EQWAT}
{\bf W}_A({\bf T}) = t_{\sf d}^{u_{\sf d}(1)} {\circ} t_{\sf
u}^{u_{\sf u}(1)}(x_1) \cdots t_{\sf d}^{u_{\sf d}(2n)} {\circ}
t_{\sf u}^{u_{\sf u}(2n)}(x_{2n}),
\end{equation}
where $x_\imath$ is the decoration on the crossing line $\imath$.
Replacing the formal representations of $\rho$ and $\rho^{-1}$ in
${\bf W}_A({\bf T})$ by $\rho$ and $\rho^{-1}$ respectively we
obtain an element ${\bf Inv}_A({\bf T}) \in A$.

For example, consider the oriented $1$--$1$ tangle diagram ${\bf
T}_{\rm trefoil}$ depicted below on the left.
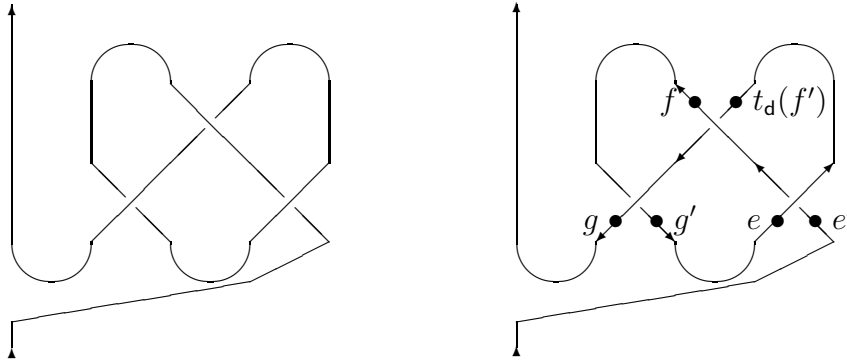
\begin{center}
\begin{picture}(140,140)(-10,-10)
\put(0,0){\line(0,-1){10}} \put(0,-10){\vector(0, 1){0}}
\put(0,30){\line(0,1){90}} \put(0,120){\vector(0, 1){0}}
\put(30,60){\line(0,1){30}} \put(120,60){\line(0,1){30}}
\put(45,90){\oval(30,30)[t]} \put(105,90){\oval(30,30)[t]}
\put(15,30){\oval(30,30)[b]} \put(75,30){\oval(30,30)[b]}
\put(30,30){\line(1,1){30}} \put(60,30){\line(-1,1){13}}
\put(30,60){\line(1,-1){13}}
\put(60,90){\line(1,-1){30}} \put(60,60){\line(1,1){13}}
\put(90,90){\line(-1,-1){13}}
\put(90,30){\line(1,1){30}} \put(90,60){\line(1,-1){13}}
\put(120,30){\line(-1,1){13}}
\put(120,30){\line(-2,-1){30}} \put(0,0){\line(6,1){90}}
%
%
\end{picture}
\qquad \qquad
\begin{picture}(140,140)(-10,-10)
\put(0,0){\line(0,-1){10}} \put(0,-9){\vector(0, 1){0}}
\put(0,30){\line(0,1){90}} \put(0,121){\vector(0, 1){0}}
\put(30,60){\line(0,1){30}} \put(120,60){\line(0,1){30}}
\put(45,90){\oval(30,30)[t]} \put(105,90){\oval(30,30)[t]}
\put(15,30){\oval(30,30)[b]} \put(75,30){\oval(30,30)[b]}
\put(30,30){\line(1,1){30}} \put(25,35){$g \; \bullet$}
\put(60,30){\line(-1,1){13}} \put(30,30){\vector(-1, -1){0}}
\put(30,60){\line(1,-1){13}} \put(60,30){\vector(1, -1){0}}
\put(50,35){$\bullet \; g'$}
\put(60,90){\line(1,-1){30}} \put(54,80){$f \; \bullet$}
\put(60,60){\vector(-1, -1){0}} \put(60,60){\line(1,1){13}}
\put(90,90){\line(-1,-1){13}} \put(60,90){\vector(-1, 1){0}}
\put(80,80){$\bullet \; t_{\sf d}(f')$}
\put(90,30){\line(1,1){30}} \put(87,35){$e \; \bullet$}
\put(90,60){\vector(-1, 1){0}} \put(90,60){\line(1,-1){13}}
\put(120,30){\line(-1,1){13}} \put(120,60){\vector(1, 1){0}}
\put(110,35){$\bullet \; e'$}
\put(120,30){\line(-2,-1){30}} \put(0,0){\line(6,1){90}}
%
%
\end{picture}
\end{center}
Traversal of the $1$--$1$ tangle diagram ${\bf T}_{\rm trefoil}$
results in the juxtaposition of labeled beads
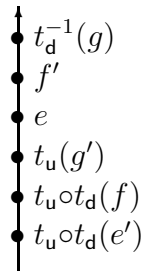
\begin{center}
\begin{picture}(40,110)(-10,-10)
\put(0,0){\line(0,-1){10}} \put(0,0){\vector(0, 1){90}}
\put(-3.5,0){$\bullet \; t_{\sf u}{\circ} t_{\sf d}(e')$}
\put(-3.5,15){$\bullet \; t_{\sf u}{\circ} t_{\sf d}(f)$}
\put(-3.5,30){$\bullet \; t_{\sf u}(g')$} \put(-3.5,45){$\bullet
\; e$} \put(-3.5,60){$\bullet \; f'$} \put(-3.5,75){$\bullet \;
t_{\sf d}^{-1}(g)$}
\end{picture}
\end{center}
Thus
$$
{\bf W}_A({\bf T}_{\rm trefoil}) = \left(\!\!\!\!\!\!\!\!
\phantom{\begin{array}{r} a  \end{array}}t_{\sf u}{\circ} t_{\sf
d}(e')\right) \left(\!\!\!\!\!\!\!\! \phantom{\begin{array}{r} a
\end{array}}t_{\sf u}{\circ} t_{\sf d}(f) \right)
\left(\!\!\!\!\!\!\!\! \phantom{\begin{array}{r} a
\end{array}}t_{\sf u}(g') \right) ef' \left(\!\!\!\!\!\!\!\!
\phantom{\begin{array}{r} a  \end{array}}t_{\sf d}^{-1}(g)\right)
$$
from which we obtain after substitution
$$
{\bf Inv}_A({\bf T}_{\rm trefoil}) = \sum_{\imath, \jmath, \ell =
1}^r \left(\!\!\!\!\!\!\!\! \phantom{\begin{array}{r} a
\end{array}}t_{\sf u}{\circ} t_{\sf d}(b_\imath) \right)
\left(\!\!\!\!\!\!\!\! \phantom{\begin{array}{r} a
\end{array}}t_{\sf u}{\circ} t_{\sf d}(a_\jmath)  \right)
\left(\!\!\!\!\!\!\!\! \phantom{\begin{array}{r} a  \end{array}}
t_{\sf u}(b_\ell) \right) a_\imath b_\jmath \left(\!\!\!\!\!\!\!\!
\phantom{\begin{array}{r} a  \end{array}}t_{\sf d}^{-1}(a_\ell)
\right),
$$
where $\rho = \sum_{\imath = 1}^r a_\imath {\otimes} b_\imath \in
A {\otimes} A$. Generally, the formal word ${\bf W}_A({\bf T})$
can be viewed as merely a device which encodes instructions for
defining an element of $A$. Since $\rho = (t {\otimes} t)(\rho)$
and $\rho^{-1} = (t {\otimes} t)(\rho^{-1})$, or symbolically $e
{\otimes} e' = t(e) {\otimes} t(e')$ and $E {\otimes} E' = t(E)
{\otimes} t(E')$ for $t = t_{\sf d}, t_{\sf u}$, we may introduce
the rules
$$
{\bf W}_A({\bf T})  =  \cdots t^p(x) \cdots t^q(y) \cdots = \cdots
t^{p+ \ell}(x) \cdots t^{q + \ell}(y) \cdots
$$
for all integers $\ell$, where $x {\otimes} y$ or $y {\otimes} x$
represents either $\rho$ or $\rho^{-1}$. Thus we may rewrite
$$
{\bf Inv}_A({\bf T}_{\rm trefoil}) = \sum_{\imath, \jmath, \ell =
1}^r \left(\!\!\!\!\!\!\!\! \phantom{\begin{array}{r} a
\end{array}}t_{\sf u}{\circ} t_{\sf d} (b_\imath) \right)
\left(\!\!\!\!\!\!\!\! \phantom{\begin{array}{r} a  \end{array}}
t_{\sf u}{\circ} t_{\sf d} (a_\jmath) \right)
\left(\!\!\!\!\!\!\!\! \phantom{\begin{array}{r} a  \end{array}}
t_{\sf u}{\circ} t_{\sf d} (b_\ell) \right) a_\imath b_\jmath
a_\ell.
$$

As a small exercise the reader is left to show that
\begin{eqnarray*}
\lefteqn{{\bf Inv}_A({\bf T}^{op}_{\rm trefoil}) } \\
 &  = & \sum_{\ell, \jmath, \imath  = 1}^r
\left(\!\!\!\!\!\!\!\! \phantom{\begin{array}{r} a  \end{array}}
t_{\sf u}^{-2}{\circ} t_{\sf d}^{-1} (a_\ell) \right)
\left(\!\!\!\!\!\!\!\! \phantom{\begin{array}{r} a  \end{array}}
t_{\sf u}^{-1}{\circ} t^{-1}_{\sf d} (b_\jmath) \right)
\left(\!\!\!\!\!\!\!\! \phantom{\begin{array}{r} a  \end{array}}
t_{\sf u}^{-1}{\circ} t^{-1}_{\sf d} (a_\imath ) \right)
\left(\!\!\!\!\!\!\!\! \phantom{\begin{array}{r} a  \end{array}}
t^{-1}_{\sf u} (b_\ell) \right)  a_\jmath b_\imath \\
& = & \sum_{\ell, \jmath, \imath  = 1}^r a_\ell b_\jmath a_\imath
\left(\!\!\!\!\!\!\!\! \phantom{\begin{array}{r} a  \end{array}}
t_{\sf u}{\circ} t_{\sf d} (b_\ell) \right) \left(\!\!\!\!\!\!\!\!
\phantom{\begin{array}{r} a  \end{array}}t_{\sf u}{\circ} t_{\sf
d}  (a_\jmath) \right) \left(\!\!\!\!\!\!\!\!
\phantom{\begin{array}{r} a  \end{array}} t_{\sf u}{\circ} t_{\sf
d} (b_\imath) \right),
\end{eqnarray*}
and also that
$$
{\bf Inv}_A({\bf T}_{\rm curl}) = \sum_{\imath = 1}^r a_\imath
\left(\!\!\!\!\!\!\!\! \phantom{\begin{array}{r} a  \end{array}}
t_{\sf u}{\circ} t_{\sf d} (b_\imath) \right),  \quad {\bf
Inv}_A({\bf T}_{\rm curl}^{op}) = \sum_{\imath = 1}^r
\left(\!\!\!\!\!\!\!\! \phantom{\begin{array}{r} a  \end{array}}
t_{\sf u}{\circ} t_{\sf d}  (b_\imath) \right) a_\imath.
$$

Assume that $A$ is a finite-dimensional and let $(C, b, T_{\sf d},
T_{\sf u}) = (A^*, b_\rho, t_{\sf d}^*, t_{\sf u}^*)$ be the dual
(strict) oriented quantum coalgebra. For all ${\bf T} \in {\bf
Tang}$ we regard ${\bf Inv}_A({\bf T}) \in A = A^{**} = C^*$ as a
functional on $C$. Here we think of $A$ as $A^{**}$ under the
identification $a(a^*) = a^*(a)$ for all $a \in A$ and $a^* \in
A^*$. We set ${\bf Inv}_C = {\bf Inv}_A$. Thus for ${\bf T} \in
{\bf Tang}$ the functional ${\bf Inv}_C({\bf T}) \in A = C^*$ is
evaluated on $c \in C$ as follows. Use (\ref{EQWAT}) to make the
formal calculation
$$
{\bf W}_A({\bf T})(c) = c({\bf W}_A({\bf T})) = c_{(1)}(t_{\sf
d}^{u_{\sf d}(1)} {\circ} t_{\sf u}^{u_{\sf u}(1)}(x_1)) \cdots
 c_{(2n)}(t_{\sf d}^{u_{\sf d}(2n)} {\circ} t_{\sf u}^{u_{\sf u}(2n)}(x_{2n}))
$$
and replace the formal copies of $\rho$ and $\rho^{-1}$ by their
actual values to obtain a scalar ${\bf Inv}_C({\bf T})(c)$.

We will evaluate ${\bf Inv}_C({\bf T}_{\rm trefoil})(c)$ to
illustrate this procedure. Recall that $b(c, d) = (c {\otimes}
d)(\rho) = \sum_{\imath = 1}^r c(a_\imath ) d(b_\imath )$ for all
$c, d \in C$. Thus we calculate, omitting the summation symbol,
\begin{eqnarray*}
\lefteqn{{\bf Inv}_C({\bf T}_{\rm trefoil})(c) } \\
& = & c \left(\!\!\!\!\!\!\!\! \phantom{\begin{array}{r} a
\end{array}} \left(\!\!\!\!\!\!\!\! \phantom{\begin{array}{r} a
\end{array}}t_{\sf u}{\circ} t_{\sf d}(b_\imath) \right)
\left(\!\!\!\!\!\!\!\! \phantom{\begin{array}{r} a
\end{array}}t_{\sf u}{\circ} t_{\sf d}(a_\jmath)  \right)
\left(\!\!\!\!\!\!\!\! \phantom{\begin{array}{r} a  \end{array}}
t_{\sf u}(b_\ell) \right) a_\imath b_\jmath
\left(\!\!\!\!\!\!\!\! \phantom{\begin{array}{r} a  \end{array}}
t_{\sf d}^{-1}(a_\ell) \right) \right) \\
& = & c_{(1)} \left(\!\!\!\!\!\!\!\! \phantom{\begin{array}{r} a
\end{array}}t_{\sf u}{\circ} t_{\sf d}(b_\imath) \right) c_{(2)}
\left(\!\!\!\!\!\!\!\! \phantom{\begin{array}{r} a
\end{array}}t_{\sf u}{\circ} t_{\sf d}(a_\jmath)  \right) c_{(3)}
\left(\!\!\!\!\!\!\!\! \phantom{\begin{array}{r} a  \end{array}}
t_{\sf u}(b_\ell) \right) c_{(4)}(a_\imath) c_{(5)}(b_\jmath)
c_{(6)}\left(\!\!\!\!\!\!\!\! \phantom{\begin{array}{r} a
\end{array}}t_{\sf d}^{-1}(a_\ell) \right)  \\
& = & \left(\!\!\!\!\!\!\!\! \phantom{\begin{array}{r} a
\end{array}}T_{\sf u}{\circ} T_{\sf d}(c_{(1)})(b_\imath) \right)
\left(\!\!\!\!\!\!\!\! \phantom{\begin{array}{r} a
\end{array}}T_{\sf u}{\circ} T_{\sf d}(c_{(2)}) (a_\jmath)
\right) \left(\!\!\!\!\!\!\!\! \phantom{\begin{array}{r} a
\end{array}} T_{\sf u}(c_{(3)})(b_\ell) \right) c_{(4)}(a_\imath)
c_{(5)}(b_\jmath)
\left(\!\!\!\!\!\!\!\! \phantom{\begin{array}{r} a  \end{array}}
T_{\sf d}^{-1}(c_{(6)})(a_\ell) \right)  \\
& = & b(c_{(4)}, T_{\sf d} {\circ} T_{\sf u} (c_{(1)})) b(T_{\sf
d} {\circ} T_{\sf u} (c_{(2)}), c_{(5)}) b(T^{-1}_d (c_{(6)}),
T_{\sf u}(c_{(3)}) )
\end{eqnarray*}
and thus
$$
{\bf Inv}_C({\bf T}_{\rm trefoil})(c)  =  b(c_{(4)}, T_{\sf d}
{\circ} T_{\sf u} (c_{(1)})) b(T_{\sf d} {\circ} T_{\sf u}
(c_{(2)}), c_{(5)}) b(T^{-1}_d (c_{(6)}), T_{\sf u}(c_{(3)}) )
$$
for all $c \in C$.

Now suppose that $(C, b, T_{\sf d}, T_{\sf u})$ is any oriented
quantum coalgebra over $k$. We shall define ${\bf Inv}_C$ in a way
which agrees with our definition when $C$ is the dual of a
finite-dimensional oriented quantum algebra over $k$.

Let ${\bf T} \in {\bf Tang}$. If ${\bf T}$ has no crossings set
${\bf Inv}_C({\bf T}) = \epsilon$. Suppose that ${\bf T}$ has $n
\geq 1$ crossings. Starting at the base point of the tangle
diagram ${\bf T}$, traverse ${\bf T}$ labeling the crossing lines
of the diagram $1, \ldots, 2n$ in the order encountered. For $1
\leq \imath \leq 2n$ let $u_{\sf d}(\imath)$ and $u_{\sf
u}(\imath)$ be as defined earlier in this section.

Let $\chi$ be a crossing and suppose that its over crossing and
under crossing lines are labeled $\imath$ and $\jmath$
respectively. For $c \in C$ the scalar ${\bf Inv}_C({\bf T})(c)$
is the sum of products, where each crossing contributes a factor
according to:
%
%
\begin{flushleft}
\mbox{
\begin{picture}(40,50)(-5,-15)
\put(30,0){\line(-1,1){13}} \put(0,0){\line(1,1){30}}
\put(0,30){\line(1,-1){13}} \put(-7,-10){$\imath$}
\put(33,-10){$\jmath$} \put(0,30){\vector(-1,1){0}}
\put(30,30){\vector(1,1){0}}
\end{picture} }
\quad \raisebox{5ex}{$\;\; {\bf Inv}_C({\bf T})(c)  =   \cdots
b(T_{\sf d}^{u_{\sf d} (\imath)} {\circ} T_{\sf u}^{u_{\sf u}
(\imath)}(c_{(\imath)}), T_{\sf d}^{u_{\sf d} (\jmath)} {\circ}
T_{\sf u}^{u_{\sf u} (\jmath)}(c_{(\jmath)})) \cdots $}
\end{flushleft}
\vskip4\jot
\begin{flushleft}
\mbox{
\begin{picture}(40,50)(-5,-15)
\put(30,0){\line(-1,1){13}} \put(0,0){\line(1,1){30}}
\put(0,30){\line(1,-1){13}} \put(-7, 35){$\jmath$}
\put(33,35){$\imath$} \put(0,0){\vector(-1,-1){0}}
\put(30,0){\vector(1,-1){0}}
\end{picture} }
\quad \raisebox{5ex}{$\;\; {\bf Inv}_C({\bf T})(c)  =   \cdots
b(T_{\sf d}^{u_{\sf d} (\imath)} {\circ} T_{\sf u}^{u_{\sf u}
(\imath)}(c_{(\imath)}), T_{\sf d}^{u_{\sf d} (\jmath)} {\circ}
T_{\sf u}^{u_{\sf u} (\jmath)}(c_{(\jmath)})) \cdots $}
\end{flushleft}
\vskip4\jot
%
%
%
\begin{flushleft}
\mbox{
\begin{picture}(40,50)(-5,-15)
\put(30,0){\line(-1,1){13}} \put(0,0){\line(1,1){30}}
\put(0,30){\line(1,-1){13}} \put(-7,-10){$\imath$} \put(-7,
35){$\jmath$} \put(30,0){\vector(1,-1){0}}
\put(30,30){\vector(1,1){0}}
\end{picture} }
\quad \raisebox{5ex}{$\;\; {\bf Inv}_C({\bf T})(c)  =   \cdots
b^{-1}(T_{\sf d}^{u_{\sf d} (\imath)} {\circ} T_{\sf u}^{u_{\sf
u}(\imath){+}1}(c_{(\imath)}), T_{\sf d}^{u_{\sf d} (\jmath)}
{\circ} T_{\sf u}^{u_{\sf u} (\jmath)}(c_{(\jmath)})) \cdots $}
\end{flushleft}
\vskip4\jot
\begin{flushleft}
\mbox{
\begin{picture}(40,50)(-5,-15)
\put(30,0){\line(-1,1){13}} \put(0,0){\line(1,1){30}}
\put(0,30){\line(1,-1){13}} \put(33,-10){$\jmath$}
\put(33,35){$\imath$} \put(0,0){\vector(-1,-1){0}}
\put(0,30){\vector(-1,1){0}}
\end{picture} }
\quad \raisebox{5ex}{$\;\; {\bf Inv}_C({\bf T})(c)  =   \cdots
b^{-1}(T_{\sf d}^{u_{\sf d}(\imath){+}1 } {\circ} T_{\sf
u}^{u_{\sf u} (\imath)}(c_{(\imath)}), T_{\sf d}^{u_{\sf d}
(\jmath)} {\circ} T_{\sf u}^{u_{\sf u} (\jmath)}(c_{(\jmath)}))
\cdots $}
\end{flushleft}
\vskip4\jot
\noindent for under crossings;
\medskip
%
%
\begin{flushleft}
\mbox{
\begin{picture}(40,50)(-5,-15)
\put(0,30){\line(1,-1){30}} \put(0,0){\line(1,1){13}}
\put(30,30){\line(-1,-1){13}} \put(-7,-10){$\jmath$}
\put(33,-10){$\imath$} \put(0,30){\vector(-1,1){0}}
\put(30,30){\vector(1,1){0}}
\end{picture} }
\quad \raisebox{5ex}{$\;\; {\bf Inv}_C({\bf T})(c)  =   \cdots
b^{-1}(T_{\sf d}^{u_{\sf d}(\imath)} {\circ} T_{\sf u}^{u_{\sf u}
(\imath)}(c_{(\imath)}), T_{\sf d}^{u_{\sf d} (\jmath)} {\circ}
T_{\sf u}^{u_{\sf u} (\jmath)}(c_{(\jmath)})) \cdots $}
\end{flushleft}
\vskip4\jot
\begin{flushleft}
\mbox{
\begin{picture}(40,50)(-5,-15)
\put(0,30){\line(1,-1){30}} \put(0,0){\line(1,1){13}}
\put(30,30){\line(-1,-1){13}} \put(-7, 35){$\imath$}
\put(33,35){$\jmath$} \put(0,0){\vector(-1,-1){0}}
\put(30,0){\vector(1,-1){0}}
\end{picture} }
\quad \raisebox{5ex}{$\;\; {\bf Inv}_C({\bf T})(c)  =   \cdots
b^{-1}(T_{\sf d}^{u_{\sf d}(\imath)} {\circ} T_{\sf u}^{u_{\sf u}
(\imath)}(c_{(\imath)}), T_{\sf d}^{u_{\sf d} (\jmath)} {\circ}
T_{\sf u}^{u_{\sf u} (\jmath)}(c_{(\jmath)})) \cdots $}
\end{flushleft}
\vskip4\jot
%
%
\begin{flushleft}
\mbox{
\begin{picture}(40,50)(-5,-15)
\put(0,30){\line(1,-1){30}} \put(0,0){\line(1,1){13}}
\put(30,30){\line(-1,-1){13}} \put(-7,-10){$\jmath$} \put(-7,
35){$\imath$} \put(30,0){\vector(1,-1){0}}
\put(30,30){\vector(1,1){0}}
\end{picture} }
\quad \raisebox{5ex}{$\;\; {\bf Inv}_C({\bf T})(c)  =   \cdots
b(T_{\sf d}^{u_{\sf d}(\imath)} {\circ} T_{\sf u}^{u_{\sf u}
(\imath)}(c_{(\imath)}), T_{\sf d}^{u_{\sf d} (\jmath)} {\circ}
T_{\sf u}^{u_{\sf u} (\jmath){+}1}(c_{(\jmath)})) \cdots $}
\end{flushleft}
\vskip4\jot
\begin{flushleft}
\mbox{
\begin{picture}(40,50)(-5,-15)
\put(0,30){\line(1,-1){30}} \put(0,0){\line(1,1){13}}
\put(30,30){\line(-1,-1){13}} \put(33,-10){$\imath$}
\put(33,35){$\jmath$} \put(0,0){\vector(-1,-1){0}}
\put(0,30){\vector(-1,1){0}}
\end{picture} }
\quad \raisebox{5ex}{$\;\; {\bf Inv}_C({\bf T})(c)  =   \cdots
b(T_{\sf d}^{u_{\sf d}(\imath)} {\circ} T_{\sf u}^{u_{\sf u}
(\imath)}(c_{(\imath)}), T_{\sf d}^{u_{\sf d} (\jmath){+}1}
{\circ} T_{\sf u}^{u_{\sf u} (\jmath) }(c_{(\jmath)})) \cdots $}
\end{flushleft}
\vskip4\jot \noindent for over crossings. In the next section we
will show that ${\bf Inv}_C$ determines a regular isotopy
invariant of $1$--$1$ tangle diagrams.

Let us reconsider the tangle diagram ${\bf T}_{\rm Trefoil}$.
Diagram traversal results in the labeling
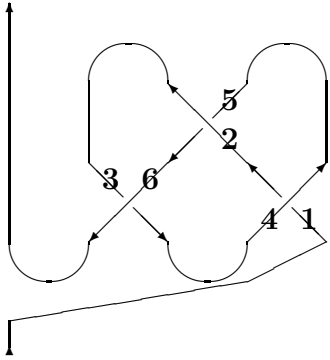
\begin{center}
\begin{picture}(140,140)(-10,-10)
\put(0,0){\line(0,-1){10}} \put(0,-9){\vector(0, 1){0}}
\put(0,30){\line(0,1){90}} \put(0,121){\vector(0, 1){0}}
\put(30,60){\line(0,1){30}} \put(120,60){\line(0,1){30}}
\put(45,90){\oval(30,30)[t]} \put(105,90){\oval(30,30)[t]}
\put(15,30){\oval(30,30)[b]} \put(75,30){\oval(30,30)[b]}
\put(30,30){\line(1,1){30}} \put(50,50){$\bf 6$}
\put(60,30){\line(-1,1){13}} \put(30,30){\vector(-1, -1){0}}
\put(30,60){\line(1,-1){13}} \put(60,30){\vector(1, -1){0}}
\put(35,50){$\bf 3$}
\put(60,90){\line(1,-1){30}} \put(80,65){$\bf 2$}
\put(60,60){\vector(-1, -1){0}} \put(60,60){\line(1,1){13}}
\put(90,90){\line(-1,-1){13}} \put(60,90){\vector(-1, 1){0}}
\put(80,80){$\bf 5$}
\put(90,30){\line(1,1){30}} \put(95,35){$\bf 4$}
\put(90,60){\vector(-1, 1){0}} \put(90,60){\line(1,-1){13}}
\put(120,30){\line(-1,1){13}} \put(120,60){\vector(1, 1){0}}
\put(110,35){$\bf 1$}
\put(120,30){\line(-2,-1){30}} \put(0,0){\line(6,1){90}}
%
%
\end{picture}
\end{center}
and thus
$$
{\bf Inv}_C({\bf T}_{\rm trefoil})(c)  =  b(c_{(4)}, T_{\sf d}
{\circ} T_{\sf u} (c_{(1)})) b(T_{\sf d} {\circ} T_{\sf u}
(c_{(2)}), c_{(5)}) b(T^{-1}_d (c_{(6)}), T_{\sf u}(c_{(3)}) )
$$
for all $c \in C$ by the algorithm described above which agrees
with our previous calculation. Observe that
$$
{\bf Inv}_C({\bf T}_{\rm curl})(c) = b(T_{\sf d} {\circ} T_{\sf
u}(c_{1}), c_{(2)}) \quad \mbox{and} \quad  {\bf Inv}_C({\bf
T}^{op}_{\rm curl})(c) = b(T_{\sf d} {\circ} T_{\sf u}(c_{2}),
c_{(1)}).
$$
Note that if ${\bf T}, {\bf T}', {\bf T}''  \in {\bf Tang}$ and
${\bf T}'' = {\bf T}{\star} {\bf T}'$ then
$$
{\bf Inv}_C({\bf T}{\star} {\bf T}') = {\bf Inv}_C({\bf T}){\bf
Inv}_C({\bf T}'),
$$
where the righthand side of the equation is the product in the
dual algebra $C^*$.
\subsection{A Proof That ${\bf Inv}_C$ Determines a Regular Isotopy
Invariant of Oriented $1$--$1$ Tangles}\label{SecORINVThm}
The sole purpose of this section is to show that the function
${\bf Inv}_C$ of Section \ref{SecORINV} determines a regular
isotopy invariant of oriented $1$--$1$ tangles. This follows by:
\begin{TH}\label{fCRegIso}
Let $(C, b, T_{\sf d}, T_{\sf u})$ be an oriented quantum
coalgebra defined over the field $k$ and suppose that ${\bf Inv}_C
: {\bf Tang} \la C^*$ is the function of  the previous section. If
${\bf T}, {\bf T}' \in {\bf Tang}$ are regularly isotopic then
${\bf Inv}_C({\bf T}) = {\bf Inv}_C({\bf T}')$.
\end{TH}

\pf The reader will find a discussion of regular isotopy, which we
assume as background material, in many references. Here we follow
the conventions of \cite{KNOTS}.

The regular isotopy equivalences are
%
\begin{flushleft}
\raisebox{6ex}{(M.1) \quad} \mbox{
\begin{picture}(60,50)(-30,-25)
\put(-30,0){\line(6,-5){30}} \put(0,-25){\line(0,-1){5}}
\put(0,25){\line(0,1){5}} \put(30,0){\line(-6,5){30}} \put(-15,
0){\oval(30,30)[t]} \put(15, 0){\oval(30,30)[b]}
\end{picture} }
\raisebox{4ex}{$\;\; \approx \;\;$}
\mbox{
\begin{picture}(10,60)(-5,0)
\put(0,-5){\line(0,1){60}}
\end{picture} }
\raisebox{4ex}{\quad and \quad}
%
%
\mbox{
\begin{picture}(60,50)(-30,-25)
\put(-30,0){\line(6,5){30}} \put(0,-25){\line(0,-1){5}}
\put(0,25){\line(0,1){5}} \put(30,0){\line(-6,-5){30}} \put(-15,
0){\oval(30,30)[b]} \put(15, 0){\oval(30,30)[t]}
\end{picture} }
\raisebox{4ex}{$\;\; \approx \;\;$}
\mbox{
\begin{picture}(10,60)(-5,0)
\put(0,-5){\line(0,1){60}}
\end{picture} }
\end{flushleft}
\vskip1\jot
%
\begin{flushleft}
\raisebox{6ex}{(M.2) \quad} \mbox{
\begin{picture}(30,70)(0,-5)
\put(0,0){\line(0,-1){5}} \put(0,60){\line(0,1){5}}
\put(30,0){\line(0,-1){5}} \put(30,60){\line(0,1){5}}
\put(0,0){\line(1,1){30}} \put(0,30){\line(1,-1){13}}
\put(30,0){\line(-1,1){13}}
\put(0,30){\line(1,1){13}} \put(30,60){\line(-1,-1){13}}
\put(0,60){\line(1,-1){30}}
\end{picture} }
\raisebox{6ex}{$\;\; \approx \;\; $}
\mbox{
\begin{picture}(30,70)(0,-5)
\put(0,-5){\line(0,1){70}} \put(30,-5){\line(0,1){70}}
\end{picture} }
\end{flushleft}
\vskip1\jot
\begin{flushleft}
\raisebox{9ex}{(M.3) \quad} \mbox{
\begin{picture}(60,100)(0,-5)
\put(0,0){\line(0,-1){5}} \put(30,0){\line(0,-1){5}}
\put(60,0){\line(0,-1){5}}
\put(0,90){\line(0,1){5}} \put(30,90){\line(0,1){5}}
\put(60,90){\line(0,1){5}}
\put(0,30){\line(0,1){30}} \put(60,0){\line(0,1){30}}
\put(60,60){\line(0,1){30}}
\put(0,0){\line(1,1){13}} \put(30,30){\line(-1,-1){13}}
\put(0,30){\line(1,-1){30}}
\put(0,60){\line(1,1){13}} \put(30,90){\line(-1,-1){13}}
\put(0,90){\line(1,-1){30}}
\put(30,30){\line(1,1){13}} \put(60,60){\line(-1,-1){13}}
\put(30,60){\line(1,-1){30}}
\end{picture} }
\raisebox{9ex}{$\;\; \approx \;\;$}
\mbox{
\begin{picture}(60,100)(0,-5)
\put(0,0){\line(0,-1){5}} \put(30,0){\line(0,-1){5}}
\put(60,0){\line(0,-1){5}}
\put(0,90){\line(0,1){5}} \put(30,90){\line(0,1){5}}
\put(60,90){\line(0,1){5}}
\put(0,0){\line(0,1){30}} \put(0,60){\line(0,1){30}}
\put(60,30){\line(0,1){30}}
\put(0,30){\line(1,1){13}} \put(30,60){\line(-1,-1){13}}
\put(30,30){\line(-1,1){30}}
\put(30,0){\line(1,1){13}} \put(60,30){\line(-1,-1){13}}
\put(60,0){\line(-1,1){30}}
\put(30,60){\line(1,1){13}} \put(60,90){\line(-1,-1){13}}
\put(60,60){\line(-1,1){30}}
\end{picture} }
\end{flushleft}
\vskip1\jot
%
%
\begin{flushleft}
\raisebox{9ex}{(M.4) \quad} \mbox{
\begin{picture}(100,100)(0,0)
\put(0,30){\line(0,-1){5}} \put(45,15){\line(0,-1){5}}
\put(90,30){\line(0,-1){5}} \put(45,90){\line(0,1){5}}
\put(0,60){\line(3,2){45}} \put(60,60){\line(1,-1){30}}
\put(30,30){\line(1,-1){15}}
\put(0,30){\line(1,1){30}} \put(0,60){\line(1,-1){13}}
\put(30,30){\line(-1,1){13}}
\put(45,60){\oval(30,30)[t]}
\end{picture} }
\raisebox{9ex}{$\;\; \approx \;\;$}
\mbox{
\begin{picture}(100,100)(0,0)
\put(0,30){\line(0,-1){5}} \put(45,15){\line(0,-1){5}}
\put(90,30){\line(0,-1){5}} \put(45,90){\line(0,1){5}}
\put(90,60){\line(-3,2){45}} \put(0,30){\line(1,1){30}}
\put(60,30){\line(-1,-1){15}}
\put(60,60){\line(1,-1){30}} \put(60,30){\line(1,1){13}}
\put(90,60){\line(-1,-1){13}}
\put(45,60){\oval(30,30)[t]}
\end{picture} }
\end{flushleft}

and
%
%
\begin{flushleft}
\raisebox{9ex}{\phantom{(M.4)} \quad} \mbox{
\begin{picture}(100,100)(0,0)
\put(0,60){\line(0,1){5}} \put(45,75){\line(0,1){5}}
\put(90,60){\line(0,1){5}} \put(45,0){\line(0,-1){5}}
\put(0,30){\line(3,-2){45}} \put(60,30){\line(1,1){30}}
\put(30,60){\line(1,1){15}}
\put(0,60){\line(1,-1){30}} \put(0,30){\line(1,1){13}}
\put(30,60){\line(-1,-1){13}}
\put(45,30){\oval(30,30)[b]}
\end{picture} }
\raisebox{6ex}{$\;\; \approx \;\;$}
\mbox{
\begin{picture}(100,100)(0,0)
\put(0,60){\line(0,1){5}} \put(45,75){\line(0,1){5}}
\put(90,60){\line(0,1){5}} \put(45,0){\line(0,-1){5}}
\put(90,30){\line(-3,-2){45}} \put(60,60){\line(-1,1){15}}
\put(0,60){\line(1,-1){30}}
\put(60,30){\line(1,1){30}} \put(60,60){\line(1,-1){13}}
\put(90,30){\line(-1,1){13}}
\put(45,30){\oval(30,30)[b]}
\end{picture} }
\end{flushleft}
\vskip1\jot and (M.2rev)--(M.4rev), which are (M.2)--(M.4)
respectively with over crossings replaced by under crossings and
vice versa. The ``twist moves"

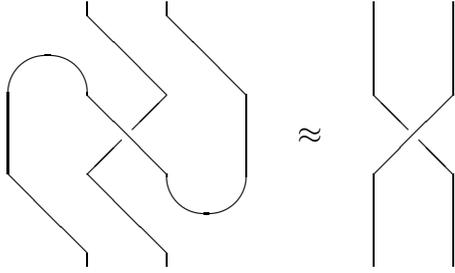
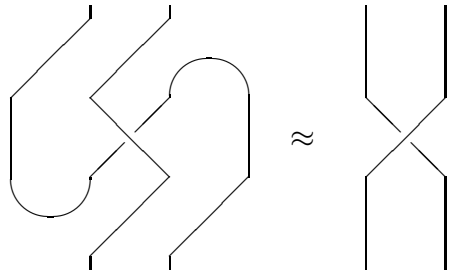
\begin{flushleft}
\mbox{
\begin{picture}(100,100)(-5,-5)
\put(30,0){\line(0,-1){5}} \put(60,0){\line(0,-1){5}}
\put(30,90){\line(0,1){5}} \put(60,90){\line(0,1){5}}
\put(0,30){\line(1,-1){30}} \put(30,30){\line(1,-1){30}}
\put(30,90){\line(1,-1){30}} \put(60,90){\line(1,-1){30}}
\put(15,60){\oval(30,30)[t]} \put(75,30){\oval(30,30)[b]}
\put(0,30){\line(0,1){30}} \put(90,30){\line(0,1){30}}
\put(30,30){\line(1,1){13}} \put(30,60){\line(1,-1){30}}
\put(60,60){\line(-1,-1){13}}
\end{picture} }
\raisebox{9ex}{$\;\; \approx \;\;$}
\mbox{
\begin{picture}(40,100)(-5,-5)
\put(0,0){\line(0,-1){5}} \put(30,0){\line(0,-1){5}}
\put(0,90){\line(0,1){5}} \put(30,90){\line(0,1){5}}
\put(0,0){\line(0,1){30}} \put(0,60){\line(0,1){30}}
\put(30,0){\line(0,1){30}} \put(30,60){\line(0,1){30}}
\put(0,30){\line(1,1){30}} \put(0,60){\line(1,-1){13}}
\put(30,30){\line(-1,1){13}}
\end{picture} }
\end{flushleft}
\raisebox{9ex}{$\qquad $ and $\;\;$}
%
\mbox{
\begin{picture}(100,100)(-5,-5)
\put(30,0){\line(0,-1){5}} \put(60,0){\line(0,-1){5}}
\put(30,90){\line(0,1){5}} \put(60,90){\line(0,1){5}}
\put(30,0){\line(1,1){30}} \put(60,0){\line(1,1){30}}
\put(0,60){\line(1,1){30}} \put(30,60){\line(1,1){30}}
\put(15,30){\oval(30,30)[b]} \put(75,60){\oval(30,30)[t]}
\put(0,30){\line(0,1){30}} \put(90,30){\line(0,1){30}}
\put(30,30){\line(1,1){13}} \put(30,60){\line(1,-1){30}}
\put(60,60){\line(-1,-1){13}}
\end{picture}
\raisebox{9ex}{$\;\; \approx \;\;$}
\mbox{
\begin{picture}(40,100)(-5,-5)
\put(0,0){\line(0,-1){5}} \put(30,0){\line(0,-1){5}}
\put(0,90){\line(0,1){5}} \put(30,90){\line(0,1){5}}
\put(0,0){\line(0,1){30}} \put(0,60){\line(0,1){30}}
\put(30,0){\line(0,1){30}} \put(30,60){\line(0,1){30}}
\put(0,30){\line(1,1){30}} \put(0,60){\line(1,-1){13}}
\put(30,30){\line(-1,1){13}}
\end{picture} }
}
\vskip1\jot \noindent are consequences of (M.1), (M.2) and (M.4).
These are important in that they allow for crossing types to be
changed.

Let ${\bf T}, {\bf T}' \in {\bf Tang}$ and suppose that a part of
${\bf T}$ is the figure on the left in one of the equivalences of
(M.1)--(M.5) or (M.2rev)--(M.5rev) and that ${\bf T}'$ is obtained
from ${\bf T}$ by replacing the figure on the left with the figure
on the right. To prove the theorem we need only show that ${\bf
Inv}_C({\bf T}) = {\bf Inv}_C({\bf T}')$. There are many cases to
consider since all possible orientations must be taken into
account. We will carefully analyze the typical cases, leaving the
remainder for the reader to work out. Let $u'_{\sf d}$ and
$u'_{\sf u}$ be the counterparts of $u_{\sf d}$ and $u_{\sf u}$
respectively for ${\bf T}'$.

Consider the first equivalence
\begin{equation}\label{EqSildeMove}
\mbox{
\begin{picture}(100,100)(0,0)
\put(0,30){\line(0,-1){5}} \put(45,15){\line(0,-1){5}}
\put(90,30){\line(0,-1){5}} \put(45,90){\line(0,1){5}}
\put(0,60){\line(3,2){45}} \put(60,60){\line(1,-1){30}}
\put(30,30){\line(1,-1){15}}
\put(0,30){\line(1,1){30}} \put(0,60){\line(1,-1){13}}
\put(30,30){\line(-1,1){13}}
\put(45,60){\oval(30,30)[t]}
\end{picture} }
\raisebox{9ex}{$\;\; \approx \;\;$}
\mbox{
\begin{picture}(100,100)(0,0)
\put(0,30){\line(0,-1){5}} \put(45,15){\line(0,-1){5}}
\put(90,30){\line(0,-1){5}} \put(45,90){\line(0,1){5}}
\put(90,60){\line(-3,2){45}} \put(0,30){\line(1,1){30}}
\put(60,30){\line(-1,-1){15}}
\put(60,60){\line(1,-1){30}} \put(60,30){\line(1,1){13}}
\put(90,60){\line(-1,-1){13}}
\put(45,60){\oval(30,30)[t]}
\end{picture} }
\end{equation}
of  (M.4). In this case
$$
\begin{picture}(40,50)(-5,-15)
\put(30,0){\line(-1,1){13}} \put(0,0){\line(1,1){30}}
\put(0,30){\line(1,-1){13}} \put(-7,-10){$\imath$}
\put(33,-10){$\jmath$}
\end{picture}
\raisebox{6ex}{\quad in ${\bf T}$ is replaced by \quad}
\begin{picture}(40,50)(-5,-15)
\put(30,30){\line(-1,-1){13}} \put(0,0){\line(1,1){13}}
\put(0,30){\line(1,-1){30}} \put(-7,-10){$\jmath$}
\put(33,-10){$\imath$}
\end{picture}
\raisebox{6ex}{\quad in ${\bf T}'$.}
$$
There are four possible orientations associated with
(\ref{EqSildeMove}).
\medskip

\noindent {\em Case} M.4.1:
$$
\mbox{
\begin{picture}(100,100)(0,0)
\put(0,30){\line(0,-1){5}} \put(45,15){\line(0,-1){5}}
\put(90,30){\line(0,-1){5}} \put(45,90){\line(0,1){5}}
\put(45,94){\vector(0,1){}} \put(90,25){\vector(0,-1){}}
\put(0,60){\line(3,2){45}} \put(60,60){\line(1,-1){30}}
\put(30,30){\line(1,-1){15}}
\put(0,30){\line(1,1){30}} \put(0,60){\line(1,-1){13}}
\put(30,30){\line(-1,1){13}}
\put(45,60){\oval(30,30)[t]}
\end{picture} }
\qquad \qquad
\mbox{
\begin{picture}(100,100)(0,0)
\put(0,30){\line(0,-1){5}} \put(45,15){\line(0,-1){5}}
\put(90,30){\line(0,-1){5}} \put(45,90){\line(0,1){5}}
\put(45,94){\vector(0,1){}} \put(90,25){\vector(0,-1){}}
\put(90,60){\line(-3,2){45}} \put(0,30){\line(1,1){30}}
\put(60,30){\line(-1,-1){15}}
\put(60,60){\line(1,-1){30}} \put(60,30){\line(1,1){13}}
\put(90,60){\line(-1,-1){13}}
\put(45,60){\oval(30,30)[t]}
\end{picture} }
$$
In this case $\quad$
\begin{picture}(40,50)(-5,-15)
\put(30,0){\line(-1,1){13}} \put(0,0){\line(1,1){30}}
\put(0,30){\line(1,-1){13}} \put(-7,-10){$\imath$}
\put(33,-10){$\jmath$} \put(0,30){\vector(-1,1){0}}
\put(30,30){\vector(1,1){0}}
\end{picture}
$\quad$ in ${\bf T}$ is replaced by $\quad$
\begin{picture}(40,50)(-5,-15)
\put(30,30){\line(-1,-1){13}} \put(0,0){\line(1,1){13}}
\put(0,30){\line(1,-1){30}} \put(-7,35){$\imath$}
\put(-7,-10){$\jmath$} \put(30,0){\vector(1,-1){0}}
\put(30,30){\vector(1,1){0}}
\end{picture}
$\quad$ in ${\bf T}'$. Observe that $u'_{\sf d}$ and $u'_{\sf u}$
agree with $u_{\sf d}$ and $u_{\sf u}$ respectively with the
exception $u'_{\sf u}(\imath) = u_{\sf u}(\imath) + 1$. Since
$T_{\sf d}, T_{\sf u}$ commute and (qc.2) holds for $b$ it follows
that
\begin{eqnarray*}
\lefteqn{b(T_{\sf d}^{u'_{\sf d}(\imath)}{\circ} T_{\sf u}^{u'_{\sf u}(\imath)}
(c_{(\imath)}), T_{\sf d}^{u'_{\sf d}(\jmath) } {\circ} T_{\sf u}^{u'_{\sf u}
(\jmath){+}1}(c_{(\jmath)}) ) } \\
& = &
b(T_{\sf d}^{u_{\sf d}(\imath)}{\circ} T_{\sf u}^{u_{\sf u}(\imath){+}1}
(c_{(\imath)}), T_{\sf d}^{u_{\sf d}(\jmath) } {\circ} T_{\sf u}^{u_{\sf u}
(\jmath){+}1}(c_{(\jmath)}) )  \\
& = & b(T_{\sf d}^{u_{\sf d}(\imath)}{\circ} T_{\sf u}^{u_{\sf
u}(\imath)} (c_{(\imath)}), T_{\sf d}^{u_{\sf d}(\jmath)} {\circ}
T_{\sf u}^{u_{\sf u}(\jmath)}(c_{(\jmath)}) )
\end{eqnarray*}
The contributions which the other crossings of ${\bf T}$ make to
${\bf Inv}_C({\bf T})(c)$ are unaffected by the replacement of the
figure on the left in (\ref{EqSildeMove}) with the right on the
right. Therefore ${\bf Inv}_C({\bf T})(c) = {\bf Inv}_C({\bf
T}')(c)$ in this case.
\medskip

\noindent {\em Case} M.4.2:
$$
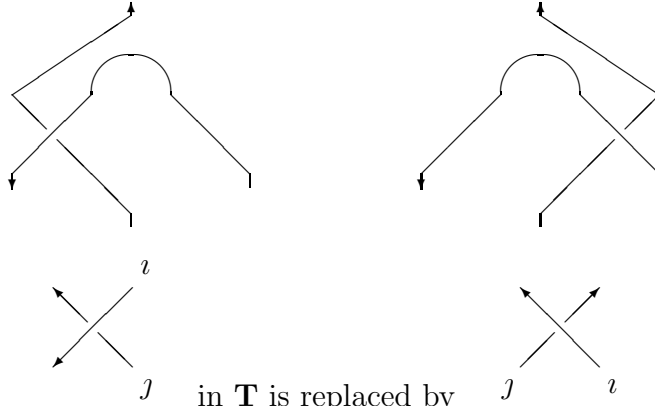
\mbox{
\begin{picture}(100,100)(0,0)
\put(0,30){\line(0,-1){5}} \put(45,15){\line(0,-1){5}}
\put(90,30){\line(0,-1){5}} \put(45,90){\line(0,1){5}}
\put(45,94){\vector(0,1){}} \put(0,25){\vector(0,-1){}}
\put(0,60){\line(3,2){45}} \put(60,60){\line(1,-1){30}}
\put(30,30){\line(1,-1){15}}
\put(0,30){\line(1,1){30}} \put(0,60){\line(1,-1){13}}
\put(30,30){\line(-1,1){13}}
\put(45,60){\oval(30,30)[t]}
\end{picture} }
\qquad \qquad
%
\mbox{
\begin{picture}(100,100)(0,0)
\put(0,30){\line(0,-1){5}} \put(45,15){\line(0,-1){5}}
\put(90,30){\line(0,-1){5}} \put(45,90){\line(0,1){5}}
\put(45,94){\vector(0,1){}} \put(0,25){\vector(0,-1){}}
\put(90,60){\line(-3,2){45}} \put(0,30){\line(1,1){30}}
\put(60,30){\line(-1,-1){15}}
\put(60,60){\line(1,-1){30}} \put(60,30){\line(1,1){13}}
\put(90,60){\line(-1,-1){13}}
\put(45,60){\oval(30,30)[t]}
\end{picture} }
$$
In this case $\quad$
\begin{picture}(40,50)(-5,-15)
\put(30,0){\line(-1,1){13}} \put(0,0){\line(1,1){30}}
\put(0,30){\line(1,-1){13}} \put(33,35){$\imath$}
\put(33,-10){$\jmath$} \put(0,30){\vector(-1,1){0}}
\put(0,0){\vector(-1,-1){0}}
\end{picture}
$\quad$ in ${\bf T}$ is replaced by $\quad$
\begin{picture}(40,50)(-5,-15)
\put(30,30){\line(-1,-1){13}} \put(0,0){\line(1,1){13}}
\put(0,30){\line(1,-1){30}} \put(-7,-10){$\jmath$}
\put(33,-10){$\imath$} \put(0,30){\vector(-1,1){0}}
\put(30,30){\vector(1,1){0}}
\end{picture}
$\quad$ in ${\bf T}'$. Observe that $u'_{\sf d}$ and $u'_{\sf u}$
agree with $u_{\sf d}$ and $u_{\sf u}$ respectively with the
exception $u'_{\sf d}(\imath) = u_{\sf d}(\imath) + 1$. Since
\begin{eqnarray*}
\lefteqn{b^{-1}(T_{\sf d}^{u'_{\sf d}(\imath)} {\circ} T_{\sf u}^{u'_{\sf u}(\imath)}
(c_{(\imath)}), T_{\sf d}^{u'_{\sf d}(\jmath)} {\circ} T_{\sf u}^{u'_{\sf u}
(\jmath)}(c_{(\jmath)}) ) } \\
& = & b^{-1}(T_{\sf d}^{u_{\sf d}(\imath) {+}1} {\circ} T_{\sf
u}^{u_{\sf u}(\imath)}(c_{(\imath)}), T_{\sf d}^{u_{\sf
d}(\jmath)} {\circ} T_{\sf u}^{u_{\sf u}(\jmath)}(c_{(\jmath)})
\end{eqnarray*}
and the contributions which the other crossings of ${\bf T}$ make
to ${\bf Inv}_C({\bf T})(c)$ are unaffected by the replacement of
the figure on the left in (\ref{EqSildeMove}) with the right on
the right, ${\bf Inv}_C({\bf T})(c) = {\bf Inv}_C({\bf T}')(c)$ in
this case.

The calculations in the other two cases, which are Cases M.4.1 and
M.4.2 with orientations reversed, are similar to those in Cases
M.4.1 and M.4.2 respectively. Thus ${\bf Inv}_C({\bf T})(c) = {\bf
Inv}_C({\bf T}')(c)$ when ${\bf T}$ is altered according to
(\ref{EqSildeMove}).

By a similar argument it follows that  ${\bf Inv}_C({\bf T})(c) =
{\bf Inv}_C({\bf T}')(c)$ when ${\bf T}$ is altered according to
the second equivalence of (M.4) and, since (qc.2) holds for
$b^{-1}$ also,  ${\bf Inv}_C({\bf T})(c) = {\bf Inv}_C({\bf
T}')(c)$ when ${\bf T}$ is altered according to (M.4rev).

It is clear that  ${\bf Inv}_C({\bf T})(c) = {\bf Inv}_C({\bf
T}')(c)$ when ${\bf T}$ is altered according to (M.1). The
following non-standard notation for the coproduct
\begin{eqnarray*}
\Delta^{(m-1)}(c) & = & c_{(1)} {\otimes} \cdots {\otimes} c_{(\ell)}
{\otimes} c_{(\ell {+} 1)} {\otimes} \cdots {\otimes} c_{(\ell')}
{\otimes} c_{(\ell' {+} 1)}\cdots {\otimes} c_{(m)} \\
& = & c_{(1)} {\otimes} \cdots {\otimes} c_{(\ell)(1)} {\otimes}
c_{(\ell)(2)} {\otimes} \cdots {\otimes} c_{(\ell')(1)} {\otimes}
c_{(\ell')(2)}\cdots {\otimes} c_{(m)}
\end{eqnarray*}
will be very useful in our analysis of (M.2). To emphasize, in the
second expression for $\Delta^{(m-1)}$ differs from the first {\em
only} in that the subscripts $(\ell), (\ell {+} 1)$ are replaced
by $(\ell)(1), (\ell)(2)$ and that $(\ell'), (\ell' {+} 1)$ are
replaced by $(\ell')(1), (\ell')(2)$. Likewise the non-standard
notation
\begin{eqnarray*}
\lefteqn{\Delta^{(m-1)}(c) } \\
& = & c_{(1)} {\otimes} \cdots {\otimes} c_{(\ell)}{\otimes} c_{(\ell {+} 1)}
\cdots {\otimes} c_{(\ell')} {\otimes} c_{(\ell' {+} 1)} \cdots
{\otimes}c_{(\ell'')} {\otimes}c_{(\ell'' {+} 1)}  \cdots {\otimes} c_{(m)} \\
& = & c_{(1)} {\otimes} \cdots {\otimes} c_{(\ell)(1)} {\otimes}
c_{(\ell)(2)} \cdots {\otimes} c_{(\ell')(1)} {\otimes}
c_{(\ell')(2)} \cdots  {\otimes}c_{(\ell'')(1)}
{\otimes}c_{(\ell'' )(2)}  \cdots {\otimes} c_{(m)}
\end{eqnarray*}
will be very useful in our analysis of our analysis of (M.3).
These manipulations with the subscripts are justified the
coassociativity of the coproduct.

We consider (M.2) next. There are four cases to analyze.
\medskip

\noindent {\em Case} M.2.1:
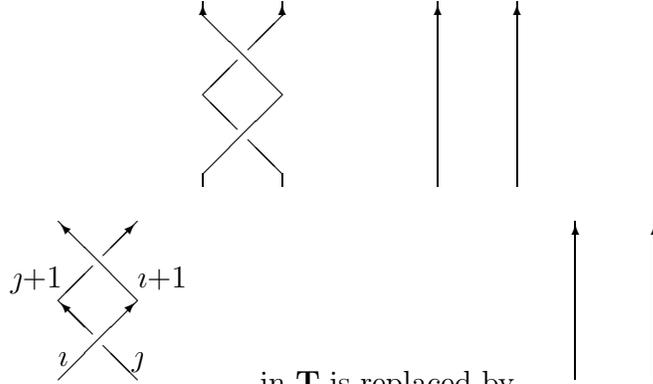
\begin{center}
\mbox{
\begin{picture}(30,70)(0,-5)
\put(0,0){\line(0,-1){5}} \put(0,60){\line(0,1){5}}
\put(30,0){\line(0,-1){5}} \put(30,60){\line(0,1){5}} \put(0,
63){\vector(0, 1){}} \put(30, 63){\vector(0, 1){}}
\put(0,0){\line(1,1){30}} \put(0,30){\line(1,-1){13}}
\put(30,0){\line(-1,1){13}}
\put(0,30){\line(1,1){13}} \put(30,60){\line(-1,-1){13}}
\put(0,60){\line(1,-1){30}}
\end{picture} }
\qquad \qquad \mbox{
\begin{picture}(30,70)(0,-5)
\put(0,-5){\line(0,1){70}} \put(30,-5){\line(0,1){70}} \put(0,
63){\vector(0, 1){}} \put(30, 63){\vector(0, 1){}}
\end{picture} }
\end{center}
In this case $\quad$
\mbox{
\begin{picture}(30,70)(0,-5)
\put(2, 58){\vector(-1,1){}} \put(28, 58){\vector(1,1){}} \put(2,
28){\vector(-1,1){}} \put(28, 28){\vector(1,1){}}
\put(0,0){\line(1,1){30}} \put(0,30){\line(1,-1){13}}
\put(30,0){\line(-1,1){13}}
\put(0,30){\line(1,1){13}} \put(30,60){\line(-1,-1){13}}
\put(0,60){\line(1,-1){30}}
\put(0,5){$\imath$} \put(-18,35){$\jmath {+} 1$}
\put(28,5){$\jmath$} \put(30,35){$\imath {+} 1$}
\end{picture} }
$\qquad \quad$ in ${\bf T}$ is replaced by $\quad$
\mbox{
\begin{picture}(30,70)(0,-5)
\put(0,0){\line(0,1){60}} \put(30,0){\line(0,1){60}} \put(0,
58){\vector(0,1){}} \put(30, 58){\vector(0,1){}}
\end{picture} }
$\quad$ in ${\bf T}'$. Since $b^{-1}$ is right inverse of $b$ and
$T_{\sf d}, T_{\sf u}$ are commuting coalgebra automorphisms of
$C$ with respect to $\{ b, b^{-1}\}$, it follows by part a) of
Lemma \ref{TmAuto} that the contribution which the two crossings
above make to the calculation of ${\bf Inv}_C({\bf T})(c)$ is
\begin{eqnarray*}
\lefteqn{b(T_{\sf d}^{u_{\sf d}(\imath)}{\circ} T_{\sf u}^{u_{\sf
u}(\imath)} (c_{(\imath)}), T_{\sf d}^{u_{\sf d} (\jmath)} {\circ}
T_{\sf u}^{u_{\sf u} (\jmath)}(c_{(\jmath)})) b^{-1}(T_{\sf
d}^{u_{\sf d}(\imath)}{\circ} T_{\sf u}^{u_{\sf u}(\imath)}
(c_{(\imath{+}1)}), T_{\sf d}^{u_{\sf d} (\jmath)} {\circ} T_{\sf
u}^{u_{\sf u}
(\jmath)} (c_{(\jmath {+} 1)})) } \\
& = & b(T_{\sf d}^{u_{\sf d}(\imath)}{\circ} T_{\sf u}^{u_{\sf
u}(\imath)} (c_{(\imath)(1)}), T_{\sf d}^{u_{\sf d} (\jmath)}
{\circ} T_{\sf u}^{u_{\sf u}
(\jmath)}(c_{(\jmath)(1)}))b^{-1}(T_{\sf d}^{u_{\sf
d}(\imath)}{\circ} T_{\sf u}^{u_{\sf u}(\imath)}
(c_{(\imath)(2)}), T_{\sf d}^{u_{\sf d} (\jmath)} {\circ} T_{\sf
u}^{u_{\sf u} (\jmath)}
(c_{(\jmath)(2)}))   \\
& = & b(T_{\sf d}^{u_{\sf d}(\imath)}{\circ} T_{\sf u}^{u_{\sf
u}(\imath)} (c_{(\imath)})_{(1)}, T_{\sf d}^{u_{\sf d} (\jmath)}
{\circ} T_{\sf u}^{u_{\sf u}
(\jmath)}(c_{(\jmath)})_{(1)})b^{-1}(T_{\sf d}^{u_{\sf
d}(\imath)}{\circ} T_{\sf u}^{u_{\sf u}(\imath)}
(c_{(\imath)})_{(2)}, T_{\sf d}^{u_{\sf d} (\jmath)} {\circ}
T_{\sf u}^{u_{\sf u} (\jmath)} (c_{(\jmath)})_{(2)})  \\
& = & \epsilon (T_{\sf d}^{u_{\sf d}(\imath)} {\circ} T_{\sf u}^{u_{\sf u}
(\imath)}(c_{(\imath)})) \epsilon (T_{\sf d}^{u_{\sf d}(\jmath)} {\circ}
T_{\sf u}^{u_{\sf u}(\jmath)}(c_{(\imath')})) \\
& = & \epsilon (c_{(\imath)}) \epsilon (c_{(\jmath)}) \\
& = & \epsilon (c_{(\imath)(1)}) \epsilon (c_{(\imath)(2)}) \epsilon
(c_{(\jmath)(1)}) \epsilon (c_{(\jmath)(2)}) \\
& = & \epsilon (c_{(\imath)}) \epsilon (c_{(\imath {+} 1)})
\epsilon (c_{(\jmath)}) \epsilon (c_{(\jmath {+} 1)}).
\end{eqnarray*}
Since the contributions which the other crossings of ${\bf T}$
make to ${\bf Inv}_C({\bf T})(c)$ are unaffected by the
replacement of the figure on the left in Case M.2.1 with the right
on the right, ${\bf Inv}_C({\bf T})(c) = {\bf Inv}_C({\bf T}')(c)$
in this case.
\medskip

\noindent {\em Case} M.2.2:
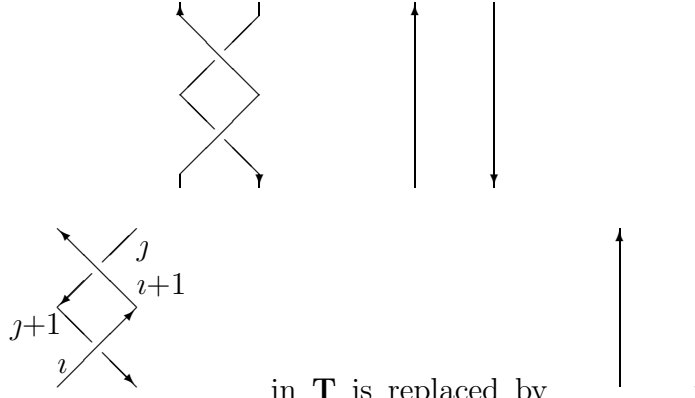
\begin{center}
\mbox{
\begin{picture}(30,70)(0,-5)
\put(0,0){\line(0,-1){5}} \put(0,60){\line(0,1){5}}
\put(30,0){\line(0,-1){5}} \put(30,60){\line(0,1){5}} \put(0,
63){\vector(0, 1){}} \put(30, -3){\vector(0, -1){}}
\put(0,0){\line(1,1){30}} \put(0,30){\line(1,-1){13}}
\put(30,0){\line(-1,1){13}}
\put(0,30){\line(1,1){13}} \put(30,60){\line(-1,-1){13}}
\put(0,60){\line(1,-1){30}}
\end{picture} }
\qquad \qquad \mbox{
\begin{picture}(30,70)(0,-5)
\put(0,-5){\line(0,1){70}} \put(30,-5){\line(0,1){70}} \put(0,
63){\vector(0, 1){}} \put(30, -3){\vector(0, -1){}}
\end{picture} }
\end{center}
In this case $\quad$
\mbox{
\begin{picture}(30,70)(0,-5)
\put(28, 28){\vector(1,1){}} \put(2, 32){\vector(-1,-1){}} \put(2,
58){\vector(-1,1){}} \put(28, 2){\vector(1,-1){}}
\put(0,0){\line(1,1){30}} \put(0,30){\line(1,-1){13}}
\put(30,0){\line(-1,1){13}}
\put(0,30){\line(1,1){13}} \put(30,60){\line(-1,-1){13}}
\put(0,60){\line(1,-1){30}}
\put(0,5){$\imath$} \put(30,50){$\jmath $} \put(-18, 20){$\jmath
{+} 1$} \put(30,35){$\imath {+} 1$}
\end{picture} }
$\qquad \quad$ in ${\bf T}$ is replaced by $\quad$
\mbox{
\begin{picture}(30,70)(0,-5)
\put(0,0){\line(0,1){60}} \put(30,0){\line(0,1){60}} \put(0,
58){\vector(0,1){}} \put(30, -3){\vector(0,-1){}}
\end{picture} }
$\quad$ in ${\bf T}'$. Since $T_{\sf d}, T_{\sf u}$ are commuting
coalgebra automorphisms of $C$ with respect to $\{ b, b^{-1}\}$,
it follows by part b) of Lemma \ref{TmAuto} and the second
equation of (qc.1) that the contribution which the two crossings
above make to the calculation of ${\bf Inv}_C({\bf T})(c)$ is
\begin{eqnarray*}
\lefteqn{b^{-1}(T_{\sf d}^{u_{\sf d}(\imath)} {\circ} T_{\sf
u}^{u_{\sf u}(\imath){+}1}(c_{(\imath)}), T_{\sf d}^{u_{\sf
d}(\jmath)} {\circ} T_{\sf u}^{u_{\sf u}(\jmath)}(c_{(\jmath)}))
b(T_{\sf d}^{u_{\sf d}(\imath)} {\circ} T_{\sf u}^{u_{\sf u}(\imath)}
(c_{(\imath{+}1)}), T_{\sf d}^{u_{\sf d}(\jmath)} {\circ} T_{\sf u}^{u_{\sf u}
(\jmath){+}1}(c_{(\jmath)})) } \\
& = & b^{-1}(T_{\sf d}^{u_{\sf d}(\imath)} {\circ} T_{\sf
u}^{u_{\sf u}(\imath){+}1}(c_{(\imath)(1)}), T_{\sf d}^{u_{\sf
d}(\jmath)} {\circ} T_{\sf u}^{u_{\sf
u}(\jmath)}(c_{(\jmath)(2)})) \times \\
& & \phantom{aaa} b(T_{\sf d}^{u_{\sf d}(\imath)} {\circ} T_{\sf
u}^{u_{\sf u}(\imath)} (c_{(\imath)(2)}), T_{\sf d}^{u_{\sf
d}(\jmath)} {\circ} T_{\sf u}^{u_{\sf u}
(\jmath){+}1}(c_{(\jmath)(1)}))  \\
& = & b^{-1}(T_{\sf u}(T_{\sf d}^{u_{\sf d}(\imath)} {\circ}
T_{\sf u}^{u_{\sf u}(\imath)}(c_{(\imath)})_{(1)}), T_{\sf
d}^{u_{\sf d}(\jmath)} {\circ} T_{\sf u}^{u_{\sf
u}(\jmath)}(c_{(\jmath)})_{(2)}) \times \\
& & \phantom{aaa} b(T_{\sf d}^{u_{\sf d}(\imath)} {\circ} T_{\sf
u}^{u_{\sf u}(\imath)} (c_{(\imath)})_{(2)}, T^{\sf u}(T_{\sf
d}^{u_{\sf d}(\jmath)}
{\circ} T_{\sf u}^{u_{\sf u}(\jmath)}(c_{(\jmath)})_{(1)})) ) \\
& = & \epsilon (T^{u(\imath)}(c_{(\imath)})) \epsilon (T^{u(\imath')}(c_{(\imath')})) \\
& = & \epsilon (T_{\sf d}^{u_{\sf d}(\imath)} {\circ} T_{\sf u}^{u_{\sf u}(\imath)}
(c_{(\imath)})) \epsilon (T_{\sf d}^{u_{\sf d}(\jmath)}
{\circ} T_{\sf u}^{u_{\sf u}(\jmath)}(c_{(\imath')})) \\
& = & \epsilon (c_{(\imath)}) \epsilon (c_{(\jmath)}) \\
& = & \epsilon (c_{(\imath)(1)}) \epsilon (c_{(\imath)(2)})
\epsilon (c_{(\jmath)(1)}) \epsilon (c_{(\jmath)(2)}) \\
& = & \epsilon (c_{(\imath)}) \epsilon (c_{(\imath {+} 1)})
\epsilon (c_{(\jmath)}) \epsilon (c_{(\jmath {+} 1)}).
\end{eqnarray*}
Since the contributions which the other crossings of ${\bf T}$
make to ${\bf Inv}_C({\bf T})(c)$ are unaffected by the
replacement of the figure on the left in Case M.2.1 with the right
on the right, ${\bf Inv}_C({\bf T})(c) = {\bf Inv}_C({\bf T}')(c)$
in this case.

Using the fact that $b^{-1}$ is a left inverse for $b$ the
argument for Cases M.2.1 is easily modified to show that ${\bf
Inv}_C({\bf T})(c) = {\bf Inv}_C({\bf T}')(c)$ in Case M.2.3,
which is Case M.2.1 with orientations reversed. Using the first
equation of (qc.1) the argument for Cases M.2.2 is easily modified
to show that ${\bf Inv}_C({\bf T})(c) = {\bf Inv}_C({\bf T}')(c)$
in Case M.2.4, which is Case M.2.2 with orientations reversed.

It remains to analyze (M.3). We consider the possible orientations
of the lines of the figures described in (M.3), reading left to
right.
\medskip

\noindent {\em Case} M.3.1
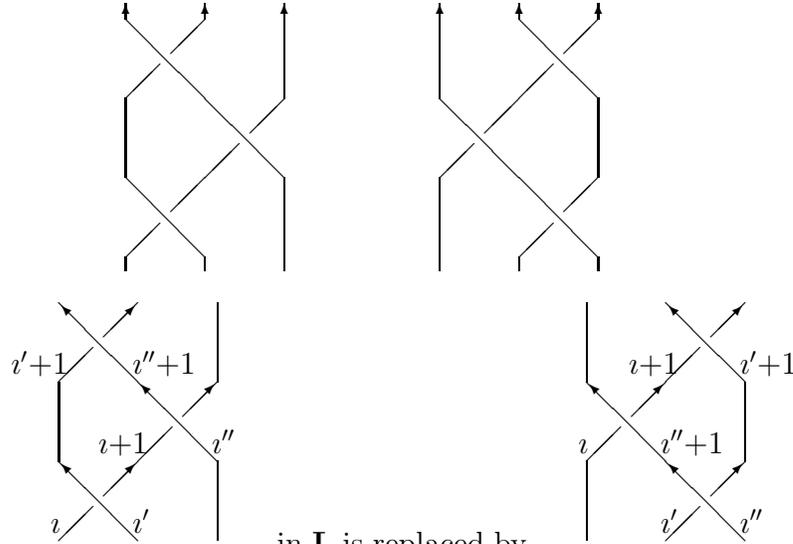
\begin{center}
\mbox{
\begin{picture}(60,100)(0,-5)
\put(0,0){\line(0,-1){5}} \put(30,0){\line(0,-1){5}}
\put(60,0){\line(0,-1){5}}
\put(0,90){\line(0,1){5}} \put(30,90){\line(0,1){5}}
\put(60,90){\line(0,1){5}} \put(0,95){\vector(0,1){}}
\put(30,95){\vector(0,1){}} \put(60,95){\vector(0,1){}}
\put(0,30){\line(0,1){30}} \put(60,0){\line(0,1){30}}
\put(60,60){\line(0,1){30}}
\put(0,0){\line(1,1){13}} \put(30,30){\line(-1,-1){13}}
\put(0,30){\line(1,-1){30}}
\put(0,60){\line(1,1){13}} \put(30,90){\line(-1,-1){13}}
\put(0,90){\line(1,-1){30}}
\put(30,30){\line(1,1){13}} \put(60,60){\line(-1,-1){13}}
\put(30,60){\line(1,-1){30}}
\end{picture} }
%
\qquad \qquad
\mbox{
\begin{picture}(60,100)(0,-5)
\put(0,0){\line(0,-1){5}} \put(30,0){\line(0,-1){5}}
\put(60,0){\line(0,-1){5}}
\put(0,90){\line(0,1){5}} \put(30,90){\line(0,1){5}}
\put(60,90){\line(0,1){5}} \put(0,95){\vector(0,1){}}
\put(30,95){\vector(0,1){}} \put(60,95){\vector(0,1){}}
\put(0,0){\line(0,1){30}} \put(0,60){\line(0,1){30}}
\put(60,30){\line(0,1){30}}
\put(0,30){\line(1,1){13}} \put(30,60){\line(-1,-1){13}}
\put(30,30){\line(-1,1){30}}
\put(30,0){\line(1,1){13}} \put(60,30){\line(-1,-1){13}}
\put(60,0){\line(-1,1){30}}
\put(30,60){\line(1,1){13}} \put(60,90){\line(-1,-1){13}}
\put(60,60){\line(-1,1){30}}
\end{picture} }
\end{center}
In this case $\quad$ \mbox{
\begin{picture}(60,100)(0,-5)
\put(0,30){\line(0,1){30}} \put(60,0){\line(0,1){30}}
\put(60,60){\line(0,1){30}}
\put(0,0){\line(1,1){13}} \put(30,30){\line(-1,-1){13}}
\put(0,30){\line(1,-1){30}} \put(2,28){\vector(-1,1){}}
\put(28,28){\vector(1,1){}} \put(-3,3){$\imath$}
\put(-18,63){$\imath' {+}1$}
\put(0,60){\line(1,1){13}} \put(30,90){\line(-1,-1){13}}
\put(0,90){\line(1,-1){30}} \put(2,88){\vector(-1,1){}}
\put(28,88){\vector(1,1){}} \put(28,3){$\imath'$}
\put(15,33){$\imath{+}1$} \put(28,63){$\imath''{+}1$}
\put(30,30){\line(1,1){13}} \put(60,60){\line(-1,-1){13}}
\put(30,60){\line(1,-1){30}} \put(32,58){\vector(-1,1){}}
\put(58,58){\vector(1,1){}} \put(58,33){$\imath''$}
\end{picture} }
$\quad$ in ${\bf L}$ is replaced by $\quad$ \mbox{
\begin{picture}(60,100)(0,-5)
\put(0,0){\line(0,1){30}} \put(0,60){\line(0,1){30}}
\put(60,30){\line(0,1){30}}
\put(0,30){\line(1,1){13}} \put(30,60){\line(-1,-1){13}}
\put(30,30){\line(-1,1){30}} \put(2,58){\vector(-1,1){}}
\put(28,58){\vector(1,1){}} \put(-3,33){$\imath$}
\put(30,0){\line(1,1){13}} \put(60,30){\line(-1,-1){13}}
\put(60,0){\line(-1,1){30}} \put(32,28){\vector(-1,1){}}
\put(58,28){\vector(1,1){}} \put(28,3){$\imath'$}
\put(28,33){$\imath''{+}1$} \put(16,63){$\imath{+}1$}
\put(30,60){\line(1,1){13}} \put(60,90){\line(-1,-1){13}}
\put(60,60){\line(-1,1){30}} \put(32,88){\vector(-1,1){}}
\put(58,88){\vector(1,1){}} \put(58,3){$\imath''$}
\put(58,63){$\imath'{+}1$}
\end{picture} }
$\quad$ in ${\bf L}'$. Since $u(\imath {+} 1) = u(\imath)$,
$u(\imath' {+} 1) = u(\imath')$ and $u(\imath'' {+} 1) =
u(\imath'')$, the contribution which the figure on the left above
makes to the calculation of ${\bf Inv}_C({\bf T})(c)$ is
$$
b^{-1}(T^{u(\imath'')}(c_{(2)}), T^{u(\imath')}(d_{(2)}))
b^{-1}(T^{u(\imath'')}(c_{(1)}), T^{u(\imath)}(e_{(2)}))
b^{-1}(T^{u(\imath')}(d_{(1)}), T^{u(\imath)}(e_{(1)}))
$$
and the contribution which the figure on the right above makes to
the calculation of ${\bf Inv}_C({\bf T}')(c)$ is
$$
b^{-1}(T^{u(\imath'')}(c_{(1)}), T^{u(\imath')}(d_{(1)}))
b^{-1}(T^{u(\imath'')} c_{(2)}, T^{u(\imath)}) (e_{(1)}))
b^{-1}(T^{u(\imath')} d_{(2)}, T^{u(\imath)}) (e_{(2)}))
$$
where $c = c_{(\imath'')}$, $d = c_{(\imath'')}$, and $e =
c_{(\imath)}$. By part a) of Lemma \ref{TmAuto} the two
contributions are the same if
\begin{eqnarray*}
\lefteqn{b^{-1}(c_{(2)}, T^v (d_{(2)})) b^{-1}(c_{(1)}, e_{(1)})
b^{-1}(T^v (d_{(1)}), e_{(2)}) } \\
& = & b^{-1}(c_{(1)}, T^v (d_{(1)})) b^{-1}(c_{(2)}, e_{(1)})
b^{-1}(T^v (d_{(2)}), e_{(2)})
\end{eqnarray*}
for all $c ,d, e \in C$. Since $b^{-1}$ satisfies (qc.2), using
part a) of Lemma \ref{TmAuto} again we see that this last equation
holds if and only if
\begin{eqnarray*}
\lefteqn{b^{-1}(T^{-v}(c_{(2)}), d_{(2)}) b^{-1}(T^{-v} (c_{(1)}),
T^{-v} (e_{(1)}))
b^{-1}(d_{(1)}, T^{-v} (e_{(2)})) } \\
& = & b^{-1}(T^{-v}(c_{(1)}), d_{(1)}) b^{-1}(T^{-v}(c_{(2)}),
T^{-v}(e_{(1)})) b^{-1}(T^v (d_{(2)}), e_{(2)})
\end{eqnarray*}
holds for all $c, d, e \in C$ if and only if
\begin{eqnarray*}
\lefteqn{b^{-1}(T^{-v}(c)_{(2)}, d_{(2)}) b^{-1}(T^{-v} (c)_{(1)},
T^{-v} (e)_{(1)})
b^{-1}(d_{(1)}, T^{-v} (e)_{(2)}) } \\
& = & b^{-1}(T^{-v}(c)_{(1)}, d_{(1)}) b^{-1}(T^{-v}(c)_{(2)},
T^{-v}(e)_{(2)}) b^{-1}(d_{(2)}, T^{-v}(e)_{(1)}))
\end{eqnarray*}
holds for all $c, d, e \in C$ which in turn holds if and only if
\begin{eqnarray*}
\lefteqn{b^{-1}(c_{(2)}, d_{(2)}) b^{-1}(c_{(1)}, e_{(1)})
b^{-1}(d_{(1)}, e_{(2)}) } \\
& = & b^{-1}(c_{(1)}, d_{(1)}) b^{-1}(c_{(2)}, e_{(1)})
b^{-1}(d_{(2)}, e_{(2)})
\end{eqnarray*}
holds for all $c, d, e \in C$. The last equation is (qc.3) for
$b^{-1}$ which holds since (qc.3) holds for $b$. Thus the two
contributions are the same. Since the contributions which the
other crossings of ${\bf T}$ make to ${\bf Inv}_C({\bf T})(c)$ are
unaffected by the replacement of the figure on the left in Case
M.3.1 with the right on the right, ${\bf Inv}_C({\bf T})(c) = {\bf
Inv}_C({\bf T}')(c)$ in Case M.3.1.

Using similar arguments one can show that ${\bf Inv}_C({\bf T})(c)
= {\bf Inv}_C({\bf T}')(c)$ for all $c \in C$ in Case M.3.2 (up up
down) if
\begin{eqnarray}
\lefteqn{b^{-1}(c_{(2)}, d_{(2)}) b(c_{(1)}, e_{(1)}) b(d_{(1)}, e_{(2)}) \nonumber } \\
& = & b^{-1}(c_{(1)}, d_{(1)}) b(c_{(2)}, e_{(2)}) b(d_{(2)},
e_{(1)}) \label{EqCaseM32}
\end{eqnarray}
for all $c, d, e \in C$, in Case M.3.3 (up down up) if
\begin{eqnarray}
\lefteqn{b(c_{(2)}, d_{(1)}) b^{-1}(c_{(1)}, e_{(2)}) b(d_{(2)}, T^2(e_{(1)})) \nonumber}  \\
&  = & b(c_{(1)}, d_{(2)}) b^{-1}(c_{(2)}, e_{(1)}) b(d_{(1)},
T^2(e_{(2)})) \label{EqCaseM33}
\end{eqnarray}
for all $c, d, e \in C$, and in Case M.3.4 (up down down) if
\begin{eqnarray}
\lefteqn{b(c_{(2)}, d_{(1)}) b(c_{(1)}, e_{(1)}) b^{-1}(d_{(2)}, e_{(2)}) \nonumber } \\
& = & b(c_{(1)}, d_{(2)}) b(c_{(2)}, e_{(2)}) b^{-1}(d_{(1)},
e_{(1)}) \label{EqCaseM34}
\end{eqnarray}
for all $c, d, e \in C$. Cases M.2.5--M.3.8, which are Cases
M.3.1--M.3.4 with orientations reversed, reduce to Cases
M.3.1--M.3.4. Cases M.3rev.1--M.3rev.8 are Cases M.3.1--M.3.8 for
the oriented quantum coalgebra $(C, b^{-1}, T^{-1})$. Thus to
complete the proof of the theorem we need only establish
(\ref{EqCaseM32})--(\ref{EqCaseM34}).

To establish (\ref{EqCaseM32}) we define linear maps $\ell, r, u :
C {\otimes} C {\otimes} C \la k$ by
$$
\ell (c {\otimes} d {\otimes} e) = b^{-1}(c_{(2)}, d_{(2)})
b(c_{(1)}, e_{(1)}) b(d_{(1)}, e_{(2)}),
$$
$$
r(c {\otimes} d {\otimes} e) = b^{-1}(c_{(1)}, d_{(1)}) b(c_{(2)},
e_{(2)}) b(d_{(2)}, e_{(1)})
$$
and
$$
u(c {\otimes} d {\otimes} e) =  b(c, d) \epsilon (e)
$$
for all $c, d, e \in C$. Since $u$ is invertible in the dual
algebra $(C {\otimes} C {\otimes} C)^*$ and $u \ell u = uru$, we
conclude that $\ell = r$, which is to say that (\ref{EqCaseM32})
holds.

To establish (\ref{EqCaseM34}) we define linear maps $\ell, r, u :
C {\otimes} C {\otimes} C \la k$ by
$$
\ell (c {\otimes} d {\otimes} e) = b(c_{(2)}, d_{(1)}) b(c_{(1)},
e_{(1)}) b^{-1}(d_{(2)}, e_{(2)}),
$$
$$
r(c {\otimes} d {\otimes} e) = b(c_{(1)}, d_{(1)}) b(c_{(2)},
e_{(2)}) b^{-1}(d_{(1)}, e_{(1)})
$$
and
$$
u(c {\otimes} d {\otimes} e) =  \epsilon (e) b(d, e)
$$
for all $c, d, e \in C$. Again, $u$ is invertible in the dual
algebra $(C {\otimes} C {\otimes} C)^*$ and again $u \ell u =
uru$. Thus $\ell = r$, or equivalently (\ref{EqCaseM32}) holds.

Equation (\ref{EqCaseM33}) is perhaps the most interesting of
(\ref{EqCaseM32})--(\ref{EqCaseM34}). Since $T^{-1}$ is a
coalgebra automorphism of $C$ with respect to $\{ b, b^{-1}\}$ by
part a) of Lemma \ref{TmAuto}, the equations of (qc.1) can be
reformulated
$$
b(c_{(1)}, T^2(d_{(2)}))b^{-1}(c_{(2)}, d_{(1)}) = \epsilon (c)
\epsilon (d)
$$
and
$$
b^{-1}(c_{(1)}, d_{(2)})b(c_{(2)}, T^2(d_{(1)})) = \epsilon (c)
\epsilon (d)
$$
for all $c. d \in C$. Thus the right hand side of the equation of
(\ref{EqCaseM33})
\begin{eqnarray*}
\lefteqn{b(c_{(1)}, d_{(2)}) b^{-1}(c_{(2)}, e_{(1)}) b(d_{(1)},
T^2(e_{(2)}))
} \\
& = &
b^{-1}(c_{(1)}, e_{(4)})b(c_{(3)}, d_{(2)}) b^{-1}(c_{(4)}, e_{(1)}) b(d_{(1)}, T^2(e_{(2)})b(c_{(2)}, T^2(e_{(3)})) \\
& = &
b^{-1}(c_{(1)}, e_{(4)})b(c_{(2)}, d_{(1)}) b^{-1}(c_{(4)}, e_{(1)}) b(d_{(2)}, T^2(e_{(3)})b(c_{(3)}, T^2(e_{(2)})) \\
& = &
b^{-1}(c_{(1)}, e_{(2)})b(c_{(2)}, d_{(1)})b(d_{(2)}, T^2(e_{(1)})) \\
& = & b(c_{(2)}, d_{(1)})b^{-1}(c_{(1)}, e_{(2)})b(d_{(2)},
T^2(e_{(1)}))
\end{eqnarray*}
is equal to the left hand side. We have established
(\ref{EqCaseM33}) which completes the proof of the theorem. \qed
\medskip

Apropos of the proof of Theorem \ref{fCRegIso}, observe that
(\ref{EqCaseM32}) and (\ref{EqCaseM34}) reduce to the
finite-dimensional case since $C$ is the sum of its
finite-dimensional subcoalgebras. For suppose that $C$ is
finite-dimensional, $b : C {\otimes} C \la k$ is a bilinear form
which is invertible and satisfies (qc.3). Let $R \in C^* {\otimes}
C^*$ be defined by $b(c, d) = R(c {\otimes} d)$ for all $c, d \in
C$. Then $R$ is invertible and $R_{1\, 2} R_{1\, 3} R_{2\, 3}  =
R_{2\, 3} R_{1\, 3}R_{1\, 2}$. Equations (\ref{EqCaseM32}) and
(\ref{EqCaseM34}) translate to $R^{-1}_{1\, 2} R_{2\, 3} R_{1\, 3}
= R_{1\, 3} R_{2\, 3}R^{-1}_{1\, 2}$ and  $R^{-1}_{2\, 3} R_{1\,
2} R_{1\, 3}  = R_{1\, 3} R_{1\, 2}R^{-1}_{2\, 3}$ respectively
which are consequences of the preceding equation. Also, once ${\bf
Inv}_C({\bf T})(c) = {\bf Inv}_C({\bf T}')(c)$ for all $c \in C$
is established in the cases for (M.2), (M.4) and in Case M.3.1,
necessarily ${\bf Inv}_C({\bf T})(c) = {\bf Inv}_C({\bf T}')(c)$
for all $c \in C$ in Case M.4.3 for topological reasons.

To calculate the invariant ${\bf Inv}_C$ we need only consider
standard oriented quantum coalgebras.
\begin{TH}\label{StdTang}
Let $(C, b, T_{\sf d}, T_{\sf u})$ be an oriented quantum
coalgebra over the field $k$. Then ${\bf Inv}{(C, b, T_{\sf d},
T_{\sf u})}({\bf T}) = {\bf Inv}_{(C, b, 1_C, T_{\sf d} {\circ}
T_{\sf u})}({\bf T})$ for all ${\bf T} \in {\bf Tang}$.
\end{TH}

\pf We may assume ${\bf T}$ has a crossing and that all of of its
crossings are oriented upward. The result now follows as $u_{\sf
u}(\imath) = u_{\sf d}(\imath)$ for all lines $\imath$ of ${\bf
T}$. See the discussion preceding \cite[Proposition 3]{KRoalg2}.
\qed
\medskip

The invariants ${\bf Inv}_C$ and ${\bf Inv}_{C^{cop}}$ have a very
natural relationship.
\begin{LE}\label{InvOP}
Let $(C, b, T_{\sf d}, T_{\sf u})$ be an oriented quantum
coalgebra over $k$. Then ${\bf Inv}_C({\bf T}^{op}) = {\bf
Inv}_{C^{cop}}({\bf T})$ for all ${\bf T} \in {\bf Tang}$.
\end{LE}

\pf We may assume that ${\bf T}$ has $n \geq 1$ crossings. A
crossing line of ${\bf T}$ which has label $\imath$ in ${\bf
T}^{op}$ has label $n + \imath - 1$ in ${\bf T}$. Let $s$ be the
sum of the local extrema of ${\bf T}$ which are oriented counter
clockwise minus the number oriented clockwise. Then $s = u(n +
\imath -1) - u^{op}(\imath)$, or equivalently $u^{op}(\imath) =
u(n + \imath -1) - s$, for all $1 \leq \imath \leq n$. Thus
\begin{eqnarray*}
{\bf Inv}_C({\bf T}^{op})(c)
&  = & \ldots b^r(T^{u^{op}}(c_{(\imath)}), T^{u^{op}}(c_{(\jmath)})) \ldots \\
& = & \ldots b^r(T^{u(n + \imath - 1) - s}(c_{(\imath)}),
T^{u(n + \jmath - 1) - s}(c_{(\jmath)})) \ldots \\
& = & \ldots b^r(T^{u(n + \imath - 1)}(c_{(\imath)}),
T^{u(n + \jmath - 1)}(c_{(\jmath)})) \ldots \\
& = & {\bf f}_{C^{cop}}({\bf T})(c)
\end{eqnarray*}
for all $c \in C$, where $r = \pm 1$. \qed
\section{Oriented $1$--$1$ Tangle Invariants Arising from
Cocommutative Oriented Quantum Coalgebras}\label{Sec11TINV}
Let $(C, b, T_{\sf d}, T_{\sf u})$ be an oriented quantum
coalgebra over $k$ and suppose that $C$ is cocommutative. To
compute ${\bf Inv}_C$ we may assume that $(C, b, T_{\sf d}, T_{\sf
u}) = (C, b, 1_C, T)$ is standard by Theorem \ref{StdTang}. Since
$C$ is cocommutative it follows by (qc.1) that $b$ and the
bilinear form $b' : C {\times} C \la k$ defined by $b'(c, d) =
b(c, T(d))$ for all $c, d \in C$ are both inverses for $b^{-1}$.
Therefore $b' = b$, and using (qc.2) we deduce
\begin{equation}\label{EqbCocomm}
b(c, T(d)) = b(c, d) = b(T(c), d)
\end{equation}
for all $c, d \in C$. Since $T$ is a coalgebra automorphism of $C$
with respect to $\{ b, b^{-1}\}$ it follows by (\ref{EqbCocomm})
that (\ref{EqbCocomm}) holds for $b^{-1}$ and $T$ as well;
therefore
\begin{equation}\label{EqbCocomm2}
b^r(T^u(c), T^v(d)) = b(c, d)
\end{equation}
for all integers $u, v$ and $c, d \in C$, where $r = \pm 1$.

Suppose that ${\bf T} \in {\bf Tang}$ is an oriented $1$--$1$
tangle diagram with $n \geq 1$ crossings. Let $c \in C$. Since $C$
is cocommutative $\Delta^{(2n -1)}(c) = c_{(\imath_1)} {\otimes}
\cdots {\otimes} c_{(\imath_{2n})}$, where $\imath_1, \ldots,
\imath_{2n}$ is any arrangement of $1, \ldots, 2n$; see
\cite[Section 7.4]{QC} for example. This last equation and
(\ref{EqbCocomm2}) show that any crossing $\chi$ of ${\bf T}$ with
crossing lines labeled $\imath$ and $\jmath$ contributes the
factor $b^{{\rm sign}\, \chi}(c_{(\imath)}, c_{(\jmath)})$ to the
formulation of ${\bf Inv}_C({\bf T})(c)$.

We follow \cite{KNOTS} in our convention for the sign of an
oriented crossing. The sign of an oriented crossing is $1$ if as
the under crossing line is traversed in the direction of
orientation the direction of the over crossing line is to the
right, otherwise the sign of the crossing is $-1$. The writhe of
an oriented $1$--$1$ tangle diagram, denoted by ${\rm writhe}\,
{\bf T}$, is $0$ if the tangle has no crossings; otherwise the
writhe is defined to be the sum of the signs of the crossings.

Now let $\chi_1, \ldots, \chi_n$ be the crossings of ${\bf T}$.
Using the cocommutativity of $C$ again we may thus write
\begin{eqnarray*}
{\bf Inv}_C({\bf T})(c) & = & b^{{\rm sign}\, \chi_1}(c_{(1)(1)},
c_{(2)(1)}) \cdots
b^{{\rm sign}\, \chi_n}(c_{(1)(n)}, c_{(2)(n)}) \\
& = & b_{(\ell)}^{{\rm sign}\, \chi_1}(c_{(1)(1)})(c_{(2)(1)})
\cdots
b_{(\ell)}^{{\rm sign}\, \chi_n}(c_{(1)(n)})(c_{(2)(n)}) \\
& = & \left(b_{(\ell)}^{{\rm sign}\, \chi_1}(c_{(1)(1)}) \cdots
b_{(\ell)}^{{\rm sign}\, \chi_n}(c_{(1)(n)})\right)(c_{(2)}) \\
& = & b_{(\ell)}^{{\rm sign}\, \chi_1 + \cdots + {\rm sign}\, \chi_n}(c_{(1)})(c_{(2)}) \\
& = & b_{(\ell)}^{{\rm writhe}\, {\bf T}}(c_{(1)})(c_{(2)}) \\
& = & b^{{\rm writhe}\, {\bf T}}(c_{(1)},  c_{(2)}).
\end{eqnarray*}
With the convention $b^0(c, d) = \epsilon (c) \epsilon (d)$ for
all $c, d \in C$, we conclude that
\begin{equation}\label{EqINVCocomm}
{\bf Inv}_C({\bf T})(c) = b^{{\rm writhe}\, {\bf T}}(c_{(1)},
c_{(2)})
\end{equation}
for all ${\bf T} \in {\bf Tang}$ and $c \in C$. Thus the regular
isotopy invariant ${\rm writhe}$ of oriented $1$--$1$ tangle
diagrams dominates ${\bf Inv}_C$, meaning that whenever ${\bf T},
{\bf T}' \in {\bf Tang}$ satisfy ${\rm writhe} \, {\bf T} = {\rm
writhe} \, {\bf T}'$ then ${\bf Inv}_C({\bf T}) = {\bf Inv}_C({\bf
T}') $.
\section{Regular Isotopy Invariants of Oriented Knots and Links
Which Arise from a Twist Oriented Quantum
Coalgebra}\label{SecKnotINV} Throughout this section $(C, b,
T_{\sf d}, T_{\sf u}, G)$ is a twist oriented quantum coalgebra
over $k$; that is $(C, b, T_{\sf d}, T_{\sf u})$ is a strict
oriented quantum coalgebra over $k$ and $G \in C^*$ in an
invertible element which satisfies $T_{\sf d}^*(G) = T_{\sf
u}^*(G) = G$ and $T_{\sf d}{\circ}T_{\sf u}(c) = G^{-1}
\rightharpoonup c \leftharpoonup G$ for all $c \in C$. The notion
of twist quantum coalgebra is introduced in \cite[Section
4]{RKinv}. Note that $(C^{cop}, b, T_{\sf d}, T_{\sf u}, G^{-1})$
is a twist oriented quantum coalgebra as well.

We represent oriented knots and links as diagrams in the plane
with respect to the vertical direction. Let ${\cal K}$ be the set
of oriented knot diagrams and ${\cal L}$ be the set of oriented
link diagrams with respect to the vertical direction. We will show
$T$-invariant cocommutative elements ${\sf c} \in C$ give rise to
scalar valued functions ${\bf f}_{C, \, {\sf c}} : {\cal L} \la k$
which are constant on the regular isotopy classes of oriented link
diagrams (and thus ${\bf f}_{C, \, {\sf c}}$ defines a regular
isotopy invariant of oriented knots and links). The function ${\bf
f}_{C, \, {\sf c}}$ restricted to the set of oriented knot
diagrams ${\cal K}$ is closely related to the function ${\bf
Inv}_C$ of Section \ref{SecORINV}.

A very important example of a cocommutative element is the trace
function ${\rm Tr} : {\rm M}_n(k) \la k$ which we regard as an
element of ${\rm C}_n(k) = {\rm M}_n(k)^*$. Since any algebra
automorphism $t$ of ${\rm M}_n(k)$ is described by $t(x) =
GxG^{-1}$ for all $x \in {\rm M}_n(k)$, where $G \in {\rm M}_n(k)$
is invertible, it follows that ${\rm Tr}$ is $T_{\sf d}, T_{\sf
u}$-invariant for all twist oriented quantum coalgebra structures
$({\rm C}_n(k), b, T_{\sf d}, T_{\sf u})$ on ${\rm C}_n(k)$. See
the corollary to \cite[Theorem 4.3.1]{HER}.
\subsection{The Function ${\bf  f}_{C, \, {\sf c}}$ Defined on Oriented
Knot Diagrams}\label{SubSecfCcK} Let ${\sf c}$ be a $T_{\sf
d}{\circ} T_{\sf u}$-invariant cocommutative element of $C$ and
suppose that ${\bf K} \in {\cal K}$. To define the scalar ${\bf
f}_{C, \, {\sf c}}({\bf K})$ we first construct a functional ${\bf
f} \in C^*$ as follows. If ${\bf K}$ has no crossings set ${\bf f}
= \epsilon$.

Suppose that ${\bf K}$ has $n \geq 1$ crossings. Choose a point
$P$ on a vertical line in the knot diagram ${\bf K}$. (There is no
harm, under regular isotopy, in inserting a vertical line at the
end of a crossing line or local extrema -- thus we may assume that
${\bf K}$ has a vertical line.) We refer to our chosen point $P$
as the starting point.

Traverse the knot diagram ${\bf K}$, starting at $P$ and moving in
the direction of the orientation, labelling the crossing lines $1,
\ldots, 2n$ in the order encountered. For $c \in C$ let ${\bf
f}(c)$ be a sum of products, where each crossing contributes a
factor by the same algorithm which was used to describe ${\bf
Inv}_C({\bf T})(c)$ in Section \ref{SecORINV}. The proof of
Theorem \ref{fCRegIso} can be repeated verbatim to show that ${\bf
f}(c)$ is unaffected by the replacement of local parts of the knot
diagram ${\bf K}$ by their equivalents according to (M.1)--(M.5)
and (M.2rev)--(M.5rev).

Let $d$ be the Whitney degree of the oriented knot diagram ${\bf
K}$. Then $2d$ is the number of local extrema with clockwise
orientation minus the number of extrema with counterclockwise
orientation. We will show that the scalar
$$
(G^d{\bf f})({\sf c}) = G^d({\sf c}_{(1)}){\bf f}({\sf c}_{(2)})
$$
does not depend on the starting point $P$.  Observe to calculate
${\bf f}$ we may assume that all crossings are oriented in the
upright position. Altering ${\bf K}$ to achieve this will not
change the Whitney degree. Thus we may assume that all crossings
are oriented in the upright position. In light of the proof of
Theorem \ref{StdTang} we may also assume that $(C, b, T_{\sf d},
T_{\sf u})$ is standard.  Set $T = T_{\sf u}$.

Consider a new starting point $P_{new}$ which precedes $P$ in the
orientation of ${\bf K}$ and has the property that traversal of
the portion of the diagram ${\bf K}$ form $P_{new}$ to $P$ in the
direction of the orientation passes through exactly one local
extremum. Let ${\bf f}_{new}$ be the analog of ${\bf f}$
constructed for $P_{new}$ and let $m+1, \ldots, 2n$ be the labels
of the crossing lines between $P_{new}$ and $P$. Set $r = 1$ if
the extremum which precedes $P$ has clockwise orientation and set
$r = -1$ otherwise. Then
$$
G^d({\sf c}_{(1)}){\bf f}({\sf c}_{(2)}) = G^d({\sf c}_{(1)}){\bf
f}_{new}({\sf c}_{(2)})
$$
if
\begin{eqnarray}
\lefteqn{G^d({\sf c}_{(1)}) T^{\ell_1}({\sf c}_{(2)(1)}) {\otimes}
\cdots {\otimes}
T^{\ell_m}({\sf c}_{(2)(m)}) {\otimes} {\sf c}_{(2)(m+1)}
{\otimes} \cdots {\otimes} {\sf c}_{(2)(2n)} \nonumber} \\
& = & G^d({\sf c}_{(1)}) T^{\ell_1+r}({\sf c}_{(2)(2n-m+1)})
{\otimes} \cdots {\otimes}
T^{\ell_m +r}({\sf c}_{(2)(2n)}) {\otimes} \nonumber \\
& & \qquad T^{-2d + r}({\sf c}_{(2)(1)}) {\otimes} \cdots
{\otimes}
 T^{-2d + r}({\sf c}_{(2)(2n -  m)}) \label{EqfInvK}
\end{eqnarray}
for all integers $\ell_1, \ldots, \ell_m$. We will establish
(\ref{EqfInvK}) by showing for all $a_1, \ldots, a_{2n} \in C^*$
that $a_1 {\otimes} \cdots {\otimes} a_{2n}$ applied to both sides
of the equation of (\ref{EqfInvK}) gives the same result.

Now $t = T^*$ is an algebra automorphism of $C^*$ since $T$ is a
coalgebra automorphism of $C$. The axioms $T^*(G) = G$ and $T(c) =
G^{-1} {\rightharpoonup} c {\leftharpoonup} G$ for all $c \in C$
translate to $t(G) = G$ and $t(a) = GaG^{-1}$ for all $a \in C^*$.
Since ${\sf c}$ is cocommutative $ab({\sf c}) = ba({\sf c})$ for
all $a, b \in C^*$. Let $a_1, \ldots, a_{2n} \in C^*$. Applying
$a_1 {\otimes} \cdots {\otimes} a_{2n}$ to the righthand side of
the equation of (\ref{EqfInvK}) gives
\begin{eqnarray*}
\lefteqn{
G^dt^{-2d + r}(a_{m+1}) \cdots t^{-2d + r}(a_{2n}) t^{\ell_1 + r}(a_1)
\cdots t^{\ell_m + r}(a_m) ({\sf c}) }\\
& = & t^r(G^dt^{-2d}(a_{m+1} \cdots a_{2n})  t^{\ell_1}(a_1) \cdots
t^{\ell_m }(a_m) )({\sf c}) \\
& = & G^dG^{-d}a_{m+1} \cdots a_{2n}G^d t^{\ell_1}(a_1) \cdots
t^{\ell_m }(a_m) (T^r({\sf c})) \\
& = & a_{m+1} \cdots a_{2n}G^d  t^{\ell_1}(a_1) \cdots t^{\ell_m }(a_m)
({\sf c}) \\
& = & G^d  t^{\ell_1}(a_1) \cdots t^{\ell_m }(a_m)(a_{m+1} \cdots
a_{2n} )({\sf c})
\end{eqnarray*}
which is $a_1 {\otimes} \cdots {\otimes} a_{2n}$ applied to the
left hand side of the equation of (\ref{EqfInvK}). We have
established (\ref{EqfInvK}).

Set
\begin{equation}\label{EqBff}
{\bf f}_{C, \, {\sf c}}({\bf K}) = (G^d{\bf f})({\sf c}).
\end{equation}
The preceding calculations show that (\ref{EqBff}) describes a
well-defined function on ${\bf K}$, which by abuse of notation we
will refer to as ${\bf f}_{C, \, {\sf c}} : {\cal K} \la k$.

Observe that the oriented knot diagram ${\bf K}$ is regularly
isotopic to an oriented knot diagram ${\bf K}({\bf T})$, where
${\bf K}({\bf T})$ is
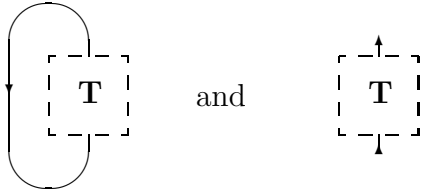
\begin{center}
\mbox{
\begin{picture}(60,80)(-5,-25)
\put(0,-5){\line(0,1){40}} \put(15,-5){\oval(30, 30)[b]}
\put(15,35){\oval(30, 30)[t]} \put(30,0){\line(0,-1){5}}
\put(30,30){\line(0,1){5}} \put(15,0){\dashbox{5}(30,30)}
\put(26,12){$\bf T$} \put(0,15){\vector(0, -1){0}}
\end{picture}}
\raisebox{7ex}{\quad and \qquad} \mbox{
\begin{picture}(30,80)(0,-25)
\put(15,0){\line(0,-1){5}} \put(15,-3){\vector(0, 1){0}}
\put(15,30){\line(0,1){5}} \put(15,38){\vector(0, 1){0}}
\put(0,0){\dashbox{5}(30,30)} \put(11,12){$\bf T$}
\end{picture}}
\end{center}
\noindent is an oriented $1$--$1$ tangle diagram, or ${\bf K}({\bf
T})$ is
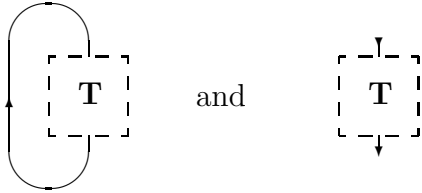
\begin{center}
\mbox{
\begin{picture}(60,80)(-5,-25)
\put(0,-5){\line(0,1){40}} \put(15,-5){\oval(30, 30)[b]}
\put(15,35){\oval(30, 30)[t]} \put(30,0){\line(0,-1){5}}
\put(30,30){\line(0,1){5}} \put(15,0){\dashbox{5}(30,30)}
\put(26,12){$\bf T$} \put(0,15){\vector(0, 1){0}}
\end{picture}}
\raisebox{7ex}{\quad and \qquad} \mbox{
\begin{picture}(30,80)(0,-25)
\put(15,0){\line(0,-1){5}} \put(15,-8){\vector(0, -1){0}}
\put(15,30){\line(0,1){5}} \put(15,33){\vector(0, -1){0}}
\put(0,0){\dashbox{5}(30,30)} \put(11,12){$\bf T$}
\end{picture}}
\end{center}
is an oriented $1$--$1$ tangle diagram. Since the Whitney degree
is a regular isotopy invariant of oriented knot diagrams, the
Whitney degrees of ${\bf K}$ and ${\bf K}({\bf T})$ are the same.
\begin{TH}\label{fCcKINV}
Let $(C, b, T_{\sf d}, T_{\sf u}, G)$ be a twist oriented quantum
coalgebra over the field $k$, let ${\sf c}$ be a $T_{\sf
d}{\circ}T_{\sf u}$-invariant cocommutative element of $C$, and
let ${\bf f}_{C, \, {\sf c}} : {\cal K} \la k$ be the function
defined by (\ref{EqBff}).
\begin{enumerate}
\item[{\rm a)}] Suppose that ${\bf K}, {\bf K}' \in {\cal K}$ are
regularly isotopic. Then ${\bf f}_{C, \, {\sf c}}({\bf K}) = {\bf
f}_{C, \, {\sf c}}({\bf K}')$. \item[{\rm b)}] Suppose that ${\bf
K} \in {\cal K}$ and that ${\bf K}$ is regularly isotopic to ${\bf
K}({\bf T})$ for some ${\bf T} \in {\bf Tang}$. Then
$$
{\bf f}_{C, \, {\sf c}}({\bf K}) = (G^d{\bf T}_C)({\sf c}),
$$
where $d$ is the Whitney degree of ${\bf K}$. \item[{\rm c)}]
${\bf f}_{C, \, {\sf c}}({\bf K}^{op}) = {\bf f}_{C^{cop}, \, {\sf
c}}({\bf K})$ for all ${\bf K} \in {\cal K}_{\rm knots}$.
\end{enumerate}
\end{TH}
\qed
\medskip

Observe that the formula in part b) of the preceding theorem may
be written
$$
{\bf f}_{C, \,  {\sf c}}({\bf K}) = (G^d{\bf T}_C)({\sf c}) = {\bf
T}_C(G^d({\sf c}_{(1)}) {\sf c}_{(2)}) = {\bf T}_C({\sf
c}{\leftharpoonup}G^d).
$$
Part c) of the preceding theorem follows with this observation
together with the fact that we may assume ${\bf K} = {\bf K}({\bf
T})$ for some ${\bf T} \in {\bf Tang}^o_{\rm tangles}$.

Observe that ${\bf f}_{C, \,  {\sf c}}({\bf K}) = G^d({\sf c})$
when ${\bf K}$ has no crossings.
\subsection{The function ${\bf  f}_{C, \, {\sf c}}$ Defined for
Oriented Link Diagrams}\label{SubSecfCcL} Let ${\bf L} \in {\cal
L}$ be an oriented link diagram with components ${\bf L}_1,
\ldots, {\bf L}_r$ and suppose that ${\sf c}$ a $T_{\sf
d}{\circ}T_{\sf u}$-invariant cocommutative element of the twist
oriented quantum coalgebra $C$. To construct the scalar ${\bf
f}_{C, \, {\sf c}}({\bf L})$ we modify the procedure for the
construction of ${\bf f}_{C, \, {\sf c}}({\bf K})$, where ${\bf K}
\in {\cal K}$ is an oriented knot diagram, described in the
preceding section.

For each $1 \leq \ell \leq r$ let $d_\ell$ denote the Whitney
degree of the component ${\bf L}_\ell$, let
$$
{\sf c}(\ell) = {\sf c} {\leftharpoonup} G^{d_\ell} =
G^{d_\ell}({\sf c}_{(1)}){\sf c}_{(2)}
$$
and choose a point on a vertical line of ${\bf L}_\ell$. We refer
to this point as a starting point. (As in the case of knot
diagrams we can always assume that each component of ${\bf L}$ has
a vertical line.) Traverse the component ${\bf L}_\ell$, beginning
at the starting point and moving in the direction of the
orientation,  labelling the crossing lines contained in ${\bf
L}_\ell$ by $(\ell {:}1), (\ell {:}2), \ldots$ in the order
encountered. Let $u(\ell {:} \imath)$ denote the number of local
extrema which are traversed in the counterclockwise direction
minus the number of  local extrema which are traversed in the
clockwise direction during the portion of traversal of the link
component from line labelled $(\ell {:} \imath)$ to the starting
point.

Next we construct a scalar ${\bf f}'_{C, \, {\sf c}}({\bf L})$. If
${\bf L}$ has no crossings we set ${\bf f}'_{C, \, {\sf c}}({\bf
L}) = 1$. Suppose that ${\bf L}$ has at least one crossing. Then
we define ${\bf f}'_{C, \, {\sf c}}({\bf L})$ to be a sum of
products, where each crossing contributes a factor of the form
$$
\cdots b^\pm (T_{\sf d}^\bullet{\circ}T_{\sf u}^\bullet(\bullet),
T_{\sf d}^\bullet{\circ}T_{\sf u}^\bullet(\bullet)) \cdots
$$
according to the conventions of Section \ref{SecORINV}, where
$(\ell {:} \imath)$ replaces $\imath$, $(\ell' {:} \imath')$
replaces $\jmath$, and then ${\sf c}(\ell)_{(\imath)}$ replaces
$c_{(\ell {:} \imath)}$ and ${\sf c}(\ell')_{(\imath')}$ replaces
$c_{(\ell' {:} \imath')}$.

We define
\begin{equation}\label{EqLINV}
{\bf f}_{C, \, {\sf c}}({\bf L}) = \omega {\bf f}'_{C, \, {\sf
c}}({\bf L}),
\end{equation}
where $\omega$ is the product of the $G^{d_\ell}({\sf c})$'s such
that the component ${\bf L}_\ell$ has no crossing lines. The
reader is left with the exercise of showing that ${\bf f}_{C, \,
{\sf c}}({\bf L})$ does not depend on the particular starting
points and is not affected by the replacement of local parts of
the diagram ${\bf L}$ by their equivalents according to
(M.1)--(M.5) and (M.2rev)--(M.5rev). The proof of Theorem
\ref{fCRegIso} provides a blueprint for the latter. Collecting
results:
\begin{TH}\label{THLINV}
Let $(C, b, T_{\sf d}, T_{\sf u}, G)$ be a twist oriented quantum
coalgebra over the field $k$, suppose that ${\sf c}$ is a $T_{\sf
d}{\circ}T_{\sf u}$-invariant cocommutative element of $C$ and let
${\bf f}_{C, \, {\sf c}} : {\cal L} \la k$ be the function
described by (\ref{EqLINV}). If ${\bf L}, {\bf L}' \in {\cal L}$
are regularly isotopic then ${\bf f}_{C, \, {\sf c}}({\bf L}) =
{\bf f}_{C, \, {\sf c}}({\bf L}')$.
\end{TH}
\qed
\medskip

Observe that ${\bf f}_{C, \, {\sf c}}$ restricted to ${\cal K}$ is
the function described in (\ref{EqBff}). By virtue of the
preceding theorem the function ${\bf f}_{C, \, {\sf c}}$
determines a regularly isotopy invariant of oriented links. When
$C$ is the dual twist quantum oriented coalgebra of a
finite-dimensional twist oriented quantum algebra $A$ over $k$
then the scalar ${\bf f}_{C, \, {\sf c}}({\bf L})$ is the
invariant $K(L)$ of \cite{RADKAUFFinv} defined for $A$. See also
\cite{kauffqalg}.

We end this section with two examples, the Hopf link and the
Borromean rings. Consider the oriented Hopf link  ${\bf L}_{Hopf}$
depicted below left with components ${\bf L}_1$ and ${\bf L}_2$,
reading left to right. The symbol  $\circ$ denotes a starting
point.
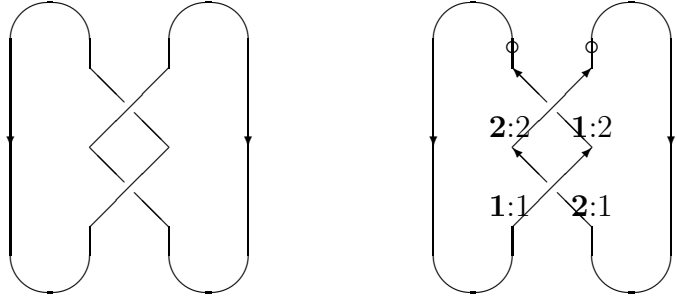
\begin{center}
\begin{picture}(280,140)(0,0)
\put(170,30){\line(0,1){80}} \put(260,30){\line(0,1){80}}
\put(200,30){\line(0,1){10}} \put(200,100){\line(0,1){10}}
\put(230,30){\line(0,1){10}} \put(230,100){\line(0,1){10}}
\put(170,70){\vector(0, -1){0}} \put(260,70){\vector(0, -1){0}}
\put(185,110){\oval(30,30)[t]} \put(245,110){\oval(30,30)[t]}
\put(185,30){\oval(30,30)[b]} \put(245,30){\oval(30,30)[b]}
\put(200,40){\line(1,1){30}} \put(200,70){\line(1,1){30}}
\put(230,70){\vector(1, 1){0}} \put(200,70){\line(1,-1){13}}
\put(200,70){\vector(-1, 1){0}} \put(230,100){\vector(1, 1){0}}
\put(200,100){\vector(-1, 1){0}} \put(200,100){\line(1,-1){13}}
\put(230,40){\line(-1,1){13}} \put(230,70){\line(-1,1){13}}
\put(10,30){\line(0,1){80}} \put(100,30){\line(0,1){80}}
\put(40,30){\line(0,1){10}} \put(40,100){\line(0,1){10}}
\put(70,30){\line(0,1){10}} \put(70,100){\line(0,1){10}}
\put(10,70){\vector(0, -1){0}} \put(100,70){\vector(0, -1){0}}
\put(25,110){\oval(30,30)[t]} \put(85,110){\oval(30,30)[t]}
\put(25,30){\oval(30,30)[b]} \put(85,30){\oval(30,30)[b]}
\put(40,40){\line(1,1){30}} \put(40,70){\line(1,-1){13}}
\put(70,40){\line(-1,1){13}} \put(40,70){\line(1,1){30}}
\put(40,100){\line(1,-1){13}} \put(70,70){\line(-1,1){13}}
\put(191, 44){${\bf 1} {:} 1$} \put(222, 44){${\bf 2} {:} 1$}
\put(191, 74){${\bf 2} {:} 2$} \put(222, 74){${\bf 1} {:} 2$}
\put(197, 105){$\circ$} \put(227, 105){$\circ$}
\end{picture}
\end{center}

Observe that $d_1 = -1$, $d_2 = 1$ and
$$
{\bf f}_{C, \, {\sf c}}({\bf L}_{Hopf}) = b(d_{(1)},
e_{(1)})b(e_{(2)}, d_{(2)}),
$$
where $d = c{\leftharpoonup}G^{-1}$ and $e = c{\leftharpoonup}G$.

Suppose ${\bf L}_{\rm Borro}$ is the Borromean rings with the
orientation given in the diagram below left and let ${\bf L}_1,
{\bf L}_2, {\bf L}_3$ be the components of ${\bf L}_{\rm Borro}$,
reading from left to right.
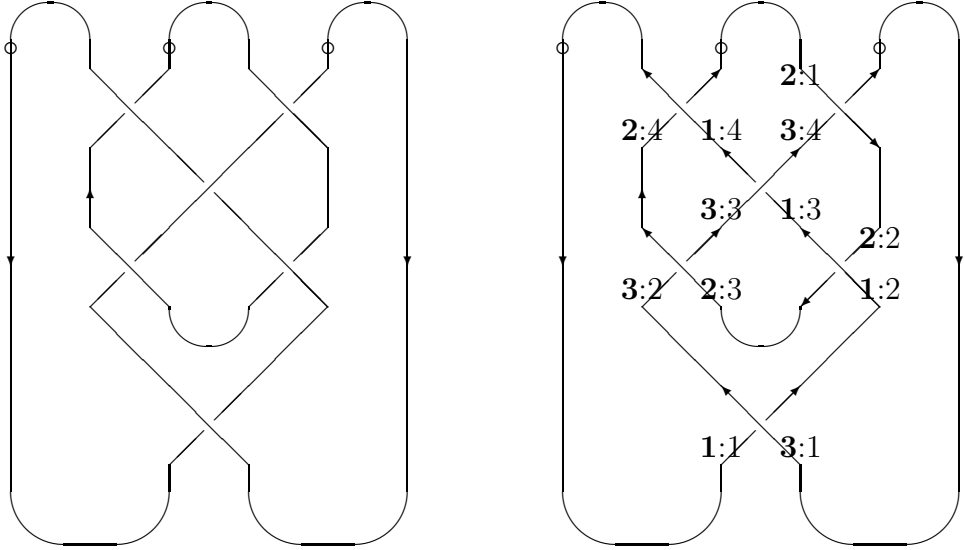
\begin{center}
\mbox{
\begin{picture}(150,200)(0,-30)
\put(30, 0){\oval(60,60)[b]} \put(120, 0){\oval(60,60)[b]}
\put(15, 160){\oval(30,30)[t]} \put(75, 160){\oval(30,30)[t]}
\put(135, 160){\oval(30,30)[t]} \put(75, 60){\oval(30,30)[b]}
\put(0,0){\line(0,1){160}} \put(150,0){\line(0,1){160}}
\put(30,90){\line(0,1){30}} \put(120,90){\line(0,1){30}}
\put(30,150){\line(0,1){10}} \put(60,150){\line(0,1){10}}
\put(90,150){\line(0,1){10}} \put(120,150){\line(0,1){10}}
\put(30,60){\line(1,-1){30}} \put(30,90){\line(1,-1){30}}
\put(30,150){\line(1,-1){30}}
\put(30,60){\line(1,1){13}} \put(30,120){\line(1,1){13}}
\put(60,90){\line(-1,-1){13}} \put(60,150){\line(-1,-1){13}}
\put(60,90){\line(1,1){30}} \put(60,30){\line(1,-1){30}}
\put(60,0){\line(1,1){13}} \put(60,120){\line(1,-1){13}}
\put(90,30){\line(-1,-1){13}} \put(90,90){\line(-1,1){13}}
\put(90,30){\line(1,1){30}} \put(90,90){\line(1,-1){30}}
\put(90,150){\line(1,-1){30}} \put(90,60){\line(1,1){13}}
\put(90,120){\line(1,1){13}} \put(120,60){\line(-1,1){13}}
\put(120,120){\line(-1,1){13}} \put(120,90){\line(-1,-1){13}}
\put(120,150){\line(-1,-1){13}}
\put(-3, 155){$\circ$} \put(57, 155){$\circ$} \put(117,
155){$\circ$}
\put(0, 75){\vector(0, -1){0}} \put(150, 75){\vector(0, -1){0}}
\put(30, 105){\vector(0, 1){0}}
\end{picture} }
\qquad \qquad
%
%
\mbox{
\begin{picture}(150,200)(0,-30)
\put(30, 0){\oval(60,60)[b]} \put(120, 0){\oval(60,60)[b]}
\put(15, 160){\oval(30,30)[t]} \put(75, 160){\oval(30,30)[t]}
\put(135, 160){\oval(30,30)[t]} \put(75, 60){\oval(30,30)[b]}
\put(0,0){\line(0,1){160}} \put(150,0){\line(0,1){160}}
\put(30,90){\line(0,1){30}} \put(120,90){\line(0,1){30}}
\put(30,150){\line(0,1){10}} \put(60,150){\line(0,1){10}}
\put(90,150){\line(0,1){10}} \put(120,150){\line(0,1){10}}
\put(30,60){\line(1,-1){30}} \put(30,90){\line(1,-1){30}}
\put(30,150){\line(1,-1){30}}
\put(30,60){\line(1,1){13}} \put(30,120){\line(1,1){13}}
\put(60,90){\line(-1,-1){13}} \put(60,150){\line(-1,-1){13}}
\put(60,90){\line(1,1){30}} \put(60,30){\line(1,-1){30}}
\put(60,0){\line(1,1){13}} \put(60,120){\line(1,-1){13}}
\put(90,30){\line(-1,-1){13}} \put(90,90){\line(-1,1){13}}
\put(90,30){\line(1,1){30}} \put(90,90){\line(1,-1){30}}
\put(90,150){\line(1,-1){30}} \put(90,60){\line(1,1){13}}
\put(90,120){\line(1,1){13}} \put(120,60){\line(-1,1){13}}
\put(120,120){\line(-1,1){13}} \put(120,90){\line(-1,-1){13}}
\put(120,150){\line(-1,-1){13}}
\put(-3, 155){$\circ$} \put(57, 155){$\circ$} \put(117,
155){$\circ$}
\put(0, 75){\vector(0, -1){0}} \put(150, 75){\vector(0, -1){0}}
\put(30, 105){\vector(0, 1){0}}
%
%
\put(60, 90){\vector(1, 1){0}} \put(60, 150){\vector(1, 1){0}}
\put(30, 90){\vector(-1, 1){0}} \put(30, 150){\vector(-1, 1){0}}
\put(90, 30){\vector(1, 1){0}} \put(90, 120){\vector(1, 1){0}}
\put(60, 30){\vector(-1, 1){0}} \put(60, 120){\vector(-1, 1){0}}
\put(90, 60){\vector(-1, -1){0}} \put(90, 90){\vector(-1, 1){0}}
\put(120, 120){\vector(1, -1){0}} \put(120, 150){\vector(1, 1){0}}
\put(22, 63){${\bf 3} {:} 2$} \put(22, 123){${\bf 2} {:} 4$}
\put(52, 3){${\bf 1} {:} 1$} \put(52, 63){${\bf 2} {:} 3$}
\put(52, 93){${\bf 3} {:} 3$} \put(52, 123){${\bf 1} {:} 4$}
\put(82, 3){${\bf 3} {:} 1$} \put(82, 93){${\bf 1} {:} 3$}
\put(82, 123){${\bf 3} {:} 4$} \put(82, 144){${\bf 2} {:} 1$}
\put(112, 63){${\bf 1} {:} 2$} \put(112, 82){${\bf 2} {:} 2$}
\end{picture} }
\end{center}

Observe that $d_1 = -1$, $d_2 = 1 = d_3$ and
\begin{eqnarray*}
\lefteqn{{\bf f}_{C, \, {\sf c}}({\bf L}_{\rm Borro}) =
b^{-1}(e_{(1)}, c_{(1)}) b^{-1}(T^2(c_{(2)}), d_{(2)}) b(e_{(3)}, c_{(3)}) {\times} } \\
& & \qquad  \qquad \qquad b^{-1}(c_{(4)}, d_{(4)}) b^{-1}(d_{(3)},
e_{(2)}) b^{-1}(d_{(1)}, e_{(4)})
\end{eqnarray*}


\begin{thebibliography}{99}

\bibitem{DRIN} V. G. Drinfel'd, Quantum Groups. {\em Proceedings of the
International Congress of Mathematicians,}
Berkeley, California, USA (1987), 798--820.

\bibitem{HER} I. N. Herstein, ``Noncommutative rings." The Carus
Mathematical Monographs, No. 15. Published by The Mathematical
Association of America; distributed by John Wiley \& Sons, Inc.,
New York 1968.

\bibitem{kauffqalg} Louis Kauffman, Gauss codes, quantum groups
and ribbon Hopf algebras. {\em Reviews in Math. Physics} {\bf 5}
(1993), 735--773.

\bibitem{KNOTS} Louis Kauffman, ``Knots and Physics.'' World Scientific,
Singapore/New Jersey/London/Hong Kong, 1991, 1994.


\bibitem{RKQCS} Louis H. Kauffman and David E. Radford, A separation
result for quantum coalgebras with an application to pointed
quantum coalgebras of small dimension. {\em  J. Algebra} {\bf 225}
(2000), 162--200.

\bibitem{RADKAUFFinv} Louis H. Kauffman and David E. Radford, Invariants
of $3$-manifolds derived from finite-dimensional Hopf algebras.
{\em  J. Knot Theory Ramifications} {\bf  4} (1995), 131--162.

\bibitem{RKinv} Louis H. Kauffman and David E. Radford, On invariants
of knots and links which arise from quantum coalgebras. {\em in
preparation}.

\bibitem{KRoalg2} Louis H. Kauffman and David E. Radford, Oriented quantum
algebras and invariants of knots and links. {\em   J. Algebra}
{\bf 246} (2001), 253--291.

\bibitem{RKO} Louis H. Kauffman and David E. Radford, Oriented quantum algebras,
categories, knots and links. {\em   J. Knot Theory Ramifications}
{\bf 10}  (2001), 1047--1084.

\bibitem{RK2DIM} Louis H. Kauffman and David E. Radford, Quantum algebra structures on
$n\times n$ matrices. {\em J. Algebra} {\bf 213} (1999), 405--436.

\bibitem{QC}  Louis H. Kauffman and David E. Radford, Quantum algebras,
quantum coalgebras, invariants of $1$--$1$ tangles and knots. {\em
Comm. Algebra}  {\bf 28}  (2000),  5101--5156.

\bibitem{LRBOOK}  Larry A. Lambe and David E. Radford, ``Introduction to the
Quantum Yang--Baxter Equation and Quantum Groups: An algebraic
Approach.'' Kluwer Academic Publishers, Boston/Dordrecht/London,
1997.

\bibitem{Mont} S. Montgomery, ``Hopf Algebras and their actions on rings.''
{\bf 82}, {\em Regional Conference Series in Mathematics}, AMS, Providence, RI, 1993.

\bibitem{RMIN} D. E. Radford, Minimal quasitriangular Hopf
algebras. {\em J. Algebra} {\bf 157} (1993), 285--315.

\bibitem{RBracket} David E. Radford, On parameterized family of twist
quantum coalgebras and the bracket polynomial. {\em  J. Algebra}
{\bf 225}  (2000), 93--123.

\bibitem{MSRI} David E. Radford, On quantum algebras and coalgebras,
oriented quantum algebras and coalgebras, invariants of $1$-$1$
tangles, knots and links. {\em New directions in Hopf algebras},
263--319, Math. Sci. Res. Inst. Publ., {\bf 43}, {\em Cambridge
Univ. Press}, Cambridge,  2002.

\bibitem{SBK} M. E. Sweedler, ``Hopf Algebras.'' Benjamin, New York, 1969.
\end{thebibliography}
\end{document}